\pdfoutput=1

\newif\ifsiam

\siamtrue
\siamfalse

\newif\ifarxiv
\ifsiam \arxivfalse \else \arxivtrue \fi

\ifsiam
\documentclass[review]{siamonline0516}

\ifsiam \usepackage{amsfonts,amssymb,amsopn} \fi 

\usepackage{mathtools}
\usepackage{graphicx}
\usepackage[caption=false]{subfig}
\usepackage{enumitem}
\usepackage{xspace}
\usepackage[nocompress]{cite}

\ifsiam
\newsiamremark{remark}{Remark}
\crefname{remark}{Remark}{Remarks}
\fi

\ifpdf
\DeclareGraphicsExtensions{.eps,.pdf,.png,.jpg}
\else
\DeclareGraphicsExtensions{.eps}
\fi

\newcommand{\TheTitle}{Switching to nonhyperbolic cycles from codimension two bifurcations of equilibria of delay differential equations} 
\newcommand{\TheAuthors}{M.M. Bosschaert, S.G. Janssens, and Yu.A. Kuznetsov}

\ifsiam
\headers{\TheTitle}{\TheAuthors}
\fi

\ifsiam \title{{\TheTitle}\thanks{Submitted to the editors DATE.}} \fi
\ifarxiv \title{\TheTitle} \fi

\author{ 
  M.M. Bosschaert\thanks{Department of Mathematics, Hasselt University,
    Diepenbeek Campus, Agoralaan Gebouw D, 3590 Diepenbeek, Belgium
    (\email{maikel.bosschaert@uhasselt.be}).}
  \and
  S.G. Janssens\thanks{Department of Mathematics, Utrecht University, 
    Budapestlaan 6, 3508 TA Utrecht, The Netherlands (\email{s.g.janssens@uu.nl}, \email{sj@dydx.nl}).}
  \and
  Yu.A. Kuznetsov\thanks{Department of Mathematics, Utrecht University, 
    Budapestlaan 6, 3508 TA Utrecht, The Netherlands and\newline 
    Department of Applied Mathematics, University of Twente, Zilverling Building, 
    7500AE Enschede, The Netherlands (\email{I.A.Kouznetsov@uu.nl}).}
}

\ifarxiv \date{March 19, 2019} \fi

\newcommand{\blist}[1]{\mbox{\lstinline!#1!}}

\definecolor{var}{rgb}{0,0.25,0.25}
\definecolor{comment}{rgb}{0,0.5,0}
\definecolor{kw}{rgb}{0,0,0.5}
\definecolor{str}{rgb}{0.5,0,0}
\definecolor{darkblue}{cmyk}{1,0,0,0.8}
\definecolor{darkred}{cmyk}{0,1,0,0.7}
\definecolor{orange}{cmyk}{0,0.5,1,0}
\definecolor{royalblue}{rgb}{0.00000,0.44700,0.74100}

\definecolor{matlabYellow}{HTML}{FCFCDC}
\definecolor{matlabrulecolor}{HTML}{A8A8A8}
\definecolor{matlabGreen}{HTML}{379634}
\definecolor{matlabLilas}{RGB}{170,55,241}

\usepackage{textcomp} 

\usepackage{listings}

\lstdefinestyle{customMatlab}{
  breaklines,showstringspaces=false,
  frame=top,frame=bottom,prebreak=...,upquote=true,
  backgroundcolor = \color{matlabYellow},
  keywordstyle=\color{blue},
  morekeywords=[2]{1}, keywordstyle=[2]{\color{black}},
  identifierstyle=\color{black},
  stringstyle=\color{matlabLilas},
  commentstyle=\color{matlabGreen},
  rulecolor=\color{matlabrulecolor},
  emph=[1]{elseif,isempty,sqrt,acos,tanh,linspace,load,clf,plot,xlabel,ylabel,zlabel,hold,exist,clear,close,all,length,figure,eps,
    colormap,real,size,axis,legend,text,pi,imag,beta,title,zeros,grid,view,strcmp,
    get, gca, set, flipud},%
  emphstyle=[1]\color{black} 
}

\lstdefinestyle{customBash}{language=bash,
  basicstyle=\small\color{black}\ttfamily,
  identifierstyle=\color{black}, 
  keywordstyle=\color{black}, 
  stringstyle=\color{black}, 
  commentstyle=\color{black},
  breaklines,
  moredelim=[is][\color{blue}\bfseries\underbar]{?}{?},
}

\lstdefinestyle{matlabConsole}{language=bash,
  basicstyle=\small\color{black}\ttfamily,
  identifierstyle=\color{black}, 
  keywordstyle=\color{black}, 
  stringstyle=\color{black}, 
  commentstyle=\color{black},
  breaklines,
  moredelim=[is][\color{blue}\bfseries\underbar]{?}{?},
}

\newcommand{\matlabrule}{\vspace{-\baselineskip}\textcolor{matlabrulecolor}{\hrule}\vspace{2pt}}

\crefname{lstlisting}{\lstlistingname}{\lstlistingname s}
\renewcommand\lstlistingname{Listing}

\usepackage{pgfplots}
\pgfplotsset{compat=newest,
  ylabsh/.style={every axis y label/.style={at={(0,0.5)}, xshift=#1, rotate=90}}} 

\usetikzlibrary{external} 
\tikzsetexternalprefix{images/} 
\tikzset{external/system call={lualatex -shell-escape -halt-on-error -interaction=batchmode -jobname "\image" "\texsource"}}
\usetikzlibrary{shapes,arrows,calc}
\usepgfplotslibrary{groupplots}

\newcommand{\PAIR}[3][]{\langle #2,#3 \rangle_{#1}}
\newcommand{\STAR}[1]{#1^{\star}}
\newcommand{\SUN}[1]{#1^{\odot}}
\newcommand{\SUNSTAR}[1]{#1^{\odot\star}}
\newcommand{\SUNSUN}[1]{#1^{\odot\odot}}
\newcommand{\STARSTAR}[1]{#1^{\star\star}}
\newcommand{\SUNSTARSTAR}[1]{#1^{\odot\star\star}}
\newcommand{\rss}{\SUNSTAR{r}}
\newcommand{\rsswp}[1]{\left[ #1 \right] \rss} 

\newcommand{\LHS}[1]{\left(#1-\SUNSTAR{A}\right)}
\newcommand{\LHSZ}{-\SUNSTAR{A}}
\newcommand{\WSTAR}{\ensuremath{\text{weak}^{\star}}}
\newcommand{\WSTARLY}{\ensuremath{\text{weakly}^{\star}}}

\newcommand{\BAR}[1]{\bar{#1}}
\newcommand{\NN}{\mathbb{N}}
\newcommand{\RR}{\mathbb{R}}
\newcommand{\RRR}[1]{\RR^{#1\star}} 
\newcommand{\CC}{\mathbb{C}}
\newcommand{\CCC}[1]{\CC^{#1\star}} 
\newcommand{\KK}{\mathbb{K}}
\newcommand{\KKK}[1]{\KK^{#1\star}} 

\let\Re\relax
\let\Im\relax
\DeclareMathOperator{\Re}{Re}
\DeclareMathOperator{\Im}{Im}

\newcommand{\DEF}{\coloneqq}

\newcommand{\DOM}{\mathcal{D}}
\newcommand{\INV}[1]{#1^{\mathrm{INV}}}
\newcommand{\BINV}[1]{B^{\mathrm{INV}}_{#1}}
\DeclareMathOperator{\DET}{det}
\DeclareMathOperator{\DIM}{dim}
\DeclareMathOperator{\diag}{diag}
\DeclareMathOperator*{\RES}{Res} 

\newcommand{\BM}[1]{\boldsymbol{#1}} 
\renewcommand{\phi}{\varphi}

\newcommand{\CM}{\mathcal{W}^c_{\text{loc}}}
\newcommand{\CMBOLD}{\BM{\mathcal{W}}^c_{\text{loc}}}

\newcommand{\MATLAB}{\texttt{MATLAB}\xspace}
\newcommand{\OCTAVE}{\texttt{GNU Octave}\xspace}
\newcommand{\DDEBIFTOOL}{\texttt{DDE-BifTool}\xspace}
\newcommand{\MATCONT}{\texttt{MatCont}\xspace}
\newcommand{\PYDELAY}{\texttt{pydelay}\xspace}

\fi

\ifarxiv
\documentclass[english,final,twoside,pdftex]{myarticle}

\fi

\ifpdf
\hypersetup{
  pdftitle={\TheTitle},
  pdfauthor={\TheAuthors}
}
\fi

\newif\ifcompileimages
\compileimagesfalse

\begin{document}

\ifarxiv \pdfbookmark[0]{Main Text}{maintext} \fi

\maketitle

\begin{abstract}
In this paper we perform the parameter-dependent center manifold reduction near the generalized Hopf (Bautin), fold-Hopf, Hopf-Hopf and transcritical-Hopf bifurcations in delay differential equations (DDEs). This allows us to initialize the continuation of codimension one equilibria and cycle bifurcations emanating from these codimension two bifurcation points. The normal form coefficients are derived in the functional analytic perturbation framework for dual semigroups (sun-star calculus) using a normalization technique based on the Fredholm alternative. The obtained expressions give explicit formulas which have been implemented in the freely available numerical software package \DDEBIFTOOL. While our theoretical results are proven to apply more generally, the software implementation and examples focus on DDEs with finitely many discrete delays. Together with the continuation capabilities of \DDEBIFTOOL, this provides a powerful tool to study the dynamics near equilibria of such DDEs. The effectiveness is demonstrated on various models.
\end{abstract}



\pagestyle{myheadings}
\thispagestyle{plain}
\markboth{\TheAuthors}{SWITCHING TO NONHYPERBOLIC CYCLES IN DDES}

\section{Introduction}\label{sec:introduction}
Great interest has recently been shown in the analysis of degenerate Hopf bifurcations in delay differential equations (DDEs), see e.g. \cite{MR2296886, MR3020901, Xu2010, MR2819829, Wang2010Hopftranscritical, Ma2011, MR3178278, MR2775253, qesmi2014HH, Agrawal2016, MR2889930, MR3047823, MR3342118, Peng2013, MR3146341, MR3430930, Song2009}.  In the simplest case, often encountered in applications, such DDEs have the form
\begin{equation}
  \label{Eq:FiniteDDE}
  \dot{x}(t)=f(x(t),x(t-\tau_1),\ldots,x(t-\tau_m),\alpha), \qquad t \geq 0,
\end{equation}
where  $x(t) \in \RR^n,\ \alpha \in \RR^p$, $f : \RR^{n \times (m+1)} \times \RR^p \to \RR^n$ is a smooth  mapping and the delays $0 < \tau_1 < \cdots <\tau_m$ are constant. They are known as \emph{discrete} DDEs or DDEs of \emph{point type}.

Using the framework of perturbation theory for dual semigroups developed in \cite{Clement1987, Clement1988, Clement1989, Clement1989b} the existence of a finite dimensional smooth center manifold for DDEs can be rigorously established \cite{diekmann1995delay}. As a consequence the normalization method for local bifurcations of ODEs developed in \cite{Kuznetsov1999} can be lifted \cite{Janssens:Thesis} rather easily to the infinite dimensional setting of DDEs. One of the advantages of this normalization technique is that the center manifold reduction and the calculation of the normal form coefficients are performed simultaneously by solving the so-called \emph{homological equation}. The method gives explicit expressions for the coefficients rather than a procedure as developed in \cite{Faria1995201, Faria1995}. The critical normal form coefficients for all five generic codimension two bifurcations of equilibria of DDEs have been derived \cite{Janssens:Thesis} and partially implemented \cite{Wage:Thesis:2014} into the fully \OCTAVE compatible \MATLAB package \DDEBIFTOOL \cite{DDEBIFTOOL,2014arXiv1406.7144S}.

In this paper we will perform the parameter-dependent center manifold reduction and normalization for three codimension two Hopf cases: the \emph{generalized Hopf, fold-Hopf} and \emph{Hopf-Hopf} bifurcations. This will allow us to initialize the continuation of codimension one bifurcation curves of nonhyperbolic equilibria and cycles emanating from the codimension two points. These are the only codimension two bifurcation points of equilibria in generic DDEs where codimension one bifurcation curves of nonhyperbolic cycles could originate. We also treat the more special transcritical-Hopf bifurcation which is frequently found in applications.

The center manifold theorem for parameter-dependent DDEs as presented in \cite{diekmann1995delay} assumes explicitly that the steady state exists for all nearby parameter values. However, for a generic fold-Hopf bifurcation this assumption is not satisfied. An attempt to deal with this complication has been made in \cite{GuoMan2011parCM}, where it is discussed how to reduce a parameter-dependent DDE to a DDE without parameters by appending the trivial equation $\dot{\alpha}=0$. However, the reduction in \cite{GuoMan2011parCM} is based on the formal adjoint approach \cite{hale1969functional} and applies specifically to DDEs, while at times it lacks consistency. Therefore we demonstrate in this paper how the reduction to the parameter-independent case can be done in the sun-star framework, enabling a rigorous approach to the existence of parameter-dependent center manifolds for a class of evolution equations that includes DDEs. This allows us to treat bifurcations of equilibria with zero eigenvalues in generic DDEs while at the same time achieving applicability of our results to other classes of delay equations.

This paper is organized as follows. In \cref{sec:sunstar} we offer a concise review of perturbation theory for dual semigroups (also called sun-star calculus), both on an abstract level as well as in application to the analysis of classical DDEs as dynamical systems. We also recall from \cite{Janssens:Thesis} various results that are needed for the normalization.

In \cref{sec:pd} we show how the theory from the previous section also applies to parameter-dependent classical DDEs by converting them into a parameter-\emph{in}dependent system on a product state space. We again present the material in two stages: Results are first established at a more abstract semigroup level and next applied to classical DDEs depending on parameters. In particular, we define the parameter-dependent local center manifold and give an explicit ODE for solutions that are confined to it.

In \cref{sec:normal-forms} we describe the general technique used to derive expressions for the normal form coefficients on the parameter-dependent center manifold in the infinite dimensional setting of classical DDEs.

Before we apply this technique to the previously mentioned codimension two bifurcations, we summarize in \cref{sec:predictors} relevant smooth normal forms and we list asymptotics for the codimension one cycle bifurcation curves emanating from the codimension two points as well as for the corresponding nonhyperbolic equilibria.

In \cref{sec:Coefficients-of-parameter} the method is then applied to the generalized Hopf (Bautin), fold-Hopf, Hopf-Hopf and transcritical-Hopf bifurcations in classical DDEs. We provide explicit expressions for all normal form coefficients necessary for the predictors of codimension one bifurcation curves. While most of the critical normal form coefficients for these bifurcations were obtained in \cite{Janssens:Thesis}, we briefly re-derive them to ensure readability.

In \cref{sec:Implement} we provide explicit computational formulas for the evaluation of the linear and multilinear forms used in the normal form coefficients and predictors for the simplest subclass consisting of DDEs \cref{Eq:FiniteDDE} with finite many discrete delays. These formulas are actually implemented in \DDEBIFTOOL.

In \cref{sec:Examples} we employ our implementation in \DDEBIFTOOL to illustrate the accuracy of the codimension one bifurcation curve predictors through various example models displaying all aforementioned degenerate Hopf cases. A complete step-by-step walk-through of the examples, including all code to reproduce the obtained results, is provided in the \hyperlink{mysupplement}{Supplement}.

\section{Dual perturbation theory and classical DDEs}\label{sec:sunstar}
We begin by presenting those general elements of perturbation theory for dual semigroups that are useful for the study of classical DDEs as dynamical systems. Throughout we assume sun-reflexivity - a term that will be introduced in \cref{sec:duality}. From \cref{sec:ddecase} onward, we then explain how the general results apply to classical DDEs. The standard reference for this entire section is \cite{diekmann1995delay}, while for the underlying theory of semigroups of linear operators we recommend \cite{Engel2000, Engel2006}.

\subsection{Duality structure and linear perturbation}\label{sec:duality}
The starting point is a $\mathcal{C}_0$-semigroup $T_0$ on a real or complex Banach space $X$. Let $A_0$ with domain $\DOM(A_0)$ be the infinitesimal generator (or: generator, for short) of $T_0$.  We denote by $\STAR{X}$ the topological dual space of $X$ and we use the prefix notation for the pairing between $\STAR{x} \in \STAR{X}$ and $x \in X$,
\[
\PAIR{\STAR{x}}{x} \DEF \STAR{x}(x).
\]
If $X$ is not reflexive then the adjoint semigroup $\STAR{T_0}$ is in general only $\WSTAR$ continuous on $\STAR{X}$ and $\STAR{A_0}$ generates $\STAR{T_0}$ only in the $\WSTAR$ sense. The maximal subspace of strong continuity
\[
\SUN{X} \DEF \left\{\STAR{x} \in \STAR{X}\,:\, t \mapsto \STAR{T_0}(t)\STAR{x} \text{ is norm-continuous on } \RR_+\right\}
\]
is invariant under $\STAR{T_0}$ and we have the characterization
\[
  \SUN{X} = \overline{\DOM(\STAR{A_0})}
\]
where the bar denotes the norm closure in $\STAR{X}$. By construction the restriction of $\STAR{T_0}$ to $\SUN{X}$ is a $\mathcal{C}_0$-semigroup that we denote by $\SUN{T_0}$. Its generator $\SUN{A_0}$ is the \emph{part} of $\STAR{A_0}$ in $\SUN{X}$,
\[
  \DOM(\SUN{A_0}) = \left\{\SUN{x} \in \DOM(\STAR{A_0}) \,:\, \STAR{A_0}\SUN{x} \in \SUN{X}\right\}, \qquad \SUN{A_0}\SUN{x} = \STAR{A_0}\SUN{x}.
\]
At this stage we again have a $\mathcal{C}_0$-semigroup $\SUN{T_0}$ with generator $\SUN{A_0}$ on a Banach space $\SUN{X}$ so we can iterate the above construction. On the dual space $\SUNSTAR{X}$ we obtain the adjoint semigroup $\SUNSTAR{T_0}$ with $\WSTAR$ generator $\SUNSTAR{A_0}$. By restriction to the maximal subspace of strong continuity $\SUNSUN{X} = \overline{\DOM(\SUNSTAR{A_0})}$ we end up with the $\mathcal{C}_0$-semigroup $\SUNSUN{T_0}$. Its generator $\SUNSUN{A_0}$ is the part of $\SUNSTAR{A_0}$ in $\SUNSUN{X}$.

The canonical injection $j : X \to \SUNSTAR{X}$ defined by
\begin{equation}
  \label{eq:j}
  \PAIR{jx}{\SUN{x}} \DEF \PAIR{\SUN{x}}{x}
\end{equation}
maps $X$ into $\SUNSUN{X}$. If $j$ maps $X$ \emph{onto} $\SUNSUN{X}$ then $X$ is called \emph{$\odot$-reflexive} (pronounce: sun-reflexive) with respect to $T_0$. One may define an equivalent norm on $X$ with respect to which $j$ becomes an isometry, but this need not be assumed. However, \emph{sun-reflexivity of $X$ with respect to $T_0$ will be assumed throughout}.

With the abstract duality structure in place, we next turn our attention to perturbation. Let $L : X \to \SUNSTAR{X}$ be a bounded linear operator. Then there exists a unique $\mathcal{C}_0$-semigroup $T$ on $X$ that satisfies the linear integral equation
\begin{equation}
  \label{eq:T_T0_AIE}
  T(t)x = T_0(t)x + j^{-1} \int_0^t \SUNSTAR{T_0}(t-\tau) L T(\tau)x \, d\tau, \qquad t \ge 0,\,x \in X,
\end{equation}
where the $\WSTAR$ Riemann integral takes values in $\SUNSUN{X}$ and the running assumption of sun-reflexivity justifies the application of $j^{-1}$. By using \cref{eq:T_T0_AIE} to express the difference $T - T_0$ of the perturbed and the unperturbed semigroups, one proves that the maximal subspaces of strong continuity $\SUN{X}$ and $\SUNSUN{X}$ are \emph{the same} for $T$ and $T_0$, so there is no need to distinguish them with a subscript. In particular, $X$ is sun-reflexive also with respect to $T$. On $\SUNSTAR{X}$ the perturbation $L$ appears additively in the action of $\SUNSTAR{A}$,
\begin{equation}
  \label{eq:A_sunstar}
  \DOM(\SUNSTAR{A}) = \DOM(\SUNSTAR{A_0}), \qquad \SUNSTAR{A} = \SUNSTAR{A_0} + Lj^{-1}.
\end{equation}
We recover the generator $A$ of $T$ by considering the part of $\SUNSTAR{A}$ in $\SUNSUN{X}$. As a consequence $L$ moves into the domain and we find
\[
  \DOM(A) = \left\{x \in X\,:\, jx \in \DOM(\SUNSTAR{A_0}) \text{ and } \SUNSTAR{A_0}jx + Lx \in \SUNSUN{X}\right\}, \quad Ax = j^{-1}(\SUNSTAR{A_0}jx + Lx).
\]
For proofs of the statements so far, see \cite[Appendix II.3 and Chapter III]{diekmann1995delay}.

\subsection{Nonlinear perturbation and linearization}\label{sec:nonlinear}
The $\mathcal{C}_0$-semigroup $T$ arose as a linear perturbation of the original $\mathcal{C}_0$-semigroup $T_0$, so the next step is to introduce a nonlinear perturbation of $T$ itself. In keeping with the tradition for nonlinear problems \cite[Sections VII.1 and VIII.1]{diekmann1995delay} we only regard the case that $X$ is a \emph{real} Banach space, also see \cref{rem:complex} below. Let $R : X \to \SUNSTAR{X}$ be a $C^k$-operator for some $k \ge 1$ such that
\[
  R(0) = 0, \qquad DR(0) = 0,
\]
and consider the nonlinear integral equation
\begin{equation}
  \label{eq:aie}
  u(t) = T(t)x + j^{-1} \int_{0}^{t}\SUNSTAR{T}(t-\tau)R(u(\tau))\,d\tau, \qquad t \ge 0,\, x \in X.
\end{equation}
Due to the nonlinearity, for a given initial condition $x \in X$ one can at most guarantee existence of a \emph{maximal solution} $u_x : I_x \to X$ of \cref{eq:aie} on a forward time interval $I_x \DEF [0,t_x)$ for some $0 < t_x \le \infty$ \cite[Chapter VII]{diekmann1995delay}. The family of all such maximal solutions defines a nonlinear semiflow $\Sigma : \DOM(\Sigma) \to X$,
\
\begin{equation}
  \label{eq:semiflow}
  \DOM(\Sigma) \DEF \{(t,x) \in [0,\infty) \times X\,:\, t \in I_x\}, \qquad \Sigma(t, x) \DEF u_x(t),
\end{equation}
that may in addition depend on parameters \cite[Defs. VII.2.1 and VII.2.9]{diekmann1995delay}. (For reasons discussed in \cref{sec:pd}, we will treat parameter dependence differently and separately. Until then, the reader can consider all parameters to be held fixed and absent in the notation.) The domain of $\Sigma$ is open in $[0,\infty) \times X$ and $0 \in X$ is a stationary point of $\Sigma$,
\[
  I_0 = [0,\infty), \qquad \Sigma(t,0) = 0, \qquad \forall\,t \ge 0.
\]
The semiflow $\Sigma$ is (in fact, uniformly) differentiable with respect to the state at $(t,0) \in \DOM(\Sigma)$, with the partial derivative
\begin{equation}\label{eq:linsigma}
D_2\Sigma(t,0) = T(t), \qquad \forall\,t \ge 0,
\end{equation}
where $T$ is the $\mathcal{C}_0$-semigroup that satisfies \cref{eq:T_T0_AIE}.

\begin{remark} \label{rem:complex}
  For nonlinear problems it is customary to work on a real Banach space $X$. The reason is that these problems often come from concrete equations with nonlinear right-hand sides for which it is unclear if and how they can be extended to complex arguments. Consequently, if we want to analyze the linearization of $\Sigma$ at $0 \in X$ using spectral theory, then it becomes necessary to \emph{complexify} $X$ and the linear operators acting on $X$ \textup{\cite[Section III.7 and last part of Section IV.2]{diekmann1995delay}, \cite[Section 1.3]{Ruston1986}}. In particular, by the spectrum of $A$ we mean the spectrum of its complexification on the complexified Banach space.   
\end{remark}

\subsection{Critical local center manifolds}\label{sec:criticalcm}
As in \cref{sec:nonlinear} we continue to assume that $T_0$ is a $\mathcal{C}_0$-semigroup on a \emph{real} Banach space $X$ that is sun-reflexive with respect to $T_0$. In addition we assume that $T_0$ is eventually compact and $L$ is a compact operator. This implies that the perturbed semigroup $T$ defined by \cref{eq:T_T0_AIE} is eventually compact as well \cite[Theorem 2.8]{diekmann2007stability}.

When considering solutions that exist for all (positive and negative) time - such as periodic orbits - it is useful to write \cref{eq:aie} in the translation invariant form
\begin{equation}
  \label{eq:AIE-st}
  u(t) = T(t-s)u(s) + j^{-1} \int_{s}^{t} \SUNSTAR{T}(t-\tau)R(u(\tau))\,d\tau, \qquad -\infty < s \leq t < \infty.
\end{equation}
A \emph{solution} of \cref{eq:AIE-st} is a continuous function $u : I \to X$ on some nondegenerate -- possibly unbounded -- interval $I \subseteq \RR$ that satisfies \cref{eq:AIE-st} for all $s, t \in I$ with $s \le t$. Naturally, $u$ is a solution of \cref{eq:AIE-st} if and only if
\[
  t - s \in I_{u(s)}, \qquad u(t) = \Sigma(t - s, u(s)), \qquad \forall\,s,t \in I \text{ with } s \le t,
\]
where $\Sigma : \DOM(\Sigma) \to X$ is the nonlinear semiflow from \cref{eq:semiflow}. The interval $I$ is often left implicit.

The general center manifold theorems from \cite[Chapter IX]{diekmann1995delay} for equations of the type \cref{eq:AIE-st} apply to the particular case where $T$ is an eventually compact $\mathcal{C}_0$-semigroup on a real, sun-reflexive Banach space. Let us therefore suppose that $0 \in X$ is a nonhyperbolic equilibrium of $\Sigma$, so the generator $A$ of $T$ possesses $1 \le n_0 < \infty$ purely imaginary eigenvalues, counting algebraic multiplicities - see \cref{rem:complex}. Let $X_0 \subseteq X$ be the \emph{real} center eigenspace corresponding to these eigenvalues. Then there exists a $C^k$-smooth $n_0$-dimensional \emph{local} center manifold $\CM$ that is tangent to $X_0$ at the origin. Any solution $u : I \to X$ of \cref{eq:AIE-st} that lies on $\CM$ is differentiable on $I$ and satisfies
\begin{equation}
  \label{eq:aode}
  j\dot{u}(t) = \SUNSTAR{A}j u(t) + R(u(t)), \qquad \forall\,t \in I,
\end{equation}
where $\SUNSTAR{A}$ is the {\WSTAR} generator of $\SUNSTAR{T}$. We note that \cref{eq:aode} is an identity in $\SUNSTAR{X}$.

\subsection{The special case of classical DDEs}\label{sec:ddecase}
It will now be explained how the general results from \cref{sec:duality,sec:nonlinear,sec:criticalcm} apply to classical DDEs. We choose the nonreflexive Banach space $X \DEF C([-h,0],\RR^n)$ as the state space, introduce a $C^k$-smooth operator $F:X \to \RR^n$, and consider an equation with a finite delay $0 < h < \infty$ of the form
\begin{equation}
  \label{eq:DDE}\tag{DDE}
  \dot{x}(t)=F(x_t), \qquad t \ge 0,
\end{equation}
with an initial condition
\begin{equation}
  \label{eq:DDE-ic}\tag{IC}
  x_0 = \phi \in X.
\end{equation}
For each $t \ge 0$, the function $x_t : [-h,0] \to \RR^n$ defined by
\[
  x_t(\theta) \DEF x(t + \theta), \qquad \forall\,\theta \in [-h,0],
\]
is called the \emph{history} of the unknown function $x$ at time $t$. Equations of the type \cref{eq:DDE} will be called \emph{classical} DDEs. Note that \cref{Eq:FiniteDDE} is quite literally a case in point. By a \emph{solution of the initial value problem} \crefrange{eq:DDE}{eq:DDE-ic} we mean a continuous function $x : [-h,t_+) \to \RR^n$ for some $0 < t_+ \le \infty$ that is differentiable on $[0,t_+)$ and satisfies \cref{eq:DDE,eq:DDE-ic}. When $t_+ = \infty$ we call $x$ a \emph{global solution}.

We want to study \cref{eq:DDE} near an equilibrium at the origin, so assume that $F(0)=0$ and split $F$ into its linear and nonlinear parts,
\[
  F(\phi)=\int_0^h d\zeta(\theta) \phi(-\theta) + G(\phi), \qquad \phi \in X.
\]
Here $\zeta : [0,h] \to \RR^{n \times n}$ is a matrix-valued function of bounded variation, normalized by the requirement that $\zeta(0) = 0$ and $\zeta$ is right-continuous on the open interval $(0,h)$. The integral is of the Riemann-Stieltjes type, and $G : X \to \RR^n$ is a $C^k$-smooth nonlinear operator. It is common to denote the linear part more succinctly as
\begin{equation}
  \label{eq:lindde_shorthand}
  \PAIR{\zeta}{\phi} \DEF \int_0^h{d\zeta(\theta)\phi(-\theta}),
\end{equation}
so that
\begin{equation}
  \label{eq:DDE-RHS}
  F(\phi) = \PAIR{\zeta}{\phi} + G(\phi),  \qquad \phi \in X.
\end{equation}
We first consider the case $G = 0$, whence \cref{eq:DDE} reduces to the linear equation
\begin{equation}
  \label{eq:lindde}
  \dot{x}(t) = \PAIR{\zeta}{x_t}, \qquad t \ge 0.
\end{equation}
In order to understand the relationship between \cref{eq:lindde} and \cref{eq:T_T0_AIE} we begin by observing that the trivial DDE
\begin{equation}
  \label{eq:trivialdde}
  \dot{x}(t) = 0, \qquad t \ge 0,
\end{equation}
with initial condition \cref{eq:DDE-ic} has the obvious solution
\[
  x^{\phi}(t) =
  \begin{cases}
    \phi(t),& t \in [-h,0],\\
    \phi(0),& t > 0.
  \end{cases}
\]
Using this solution, we define the strongly continuous \emph{shift semigroup} $T_0$ on $X$ by
\begin{equation}
  \label{eq:LDDE}
 (T_0(t)\phi)(\theta) \DEF x^{\phi}(t + \theta) =
  \begin{cases}
    \phi(t + \theta),& t + \theta \in [-h,0],\\
    \phi(0),& t + \theta > 0.
  \end{cases}
\end{equation}
We note that $T_0(h)$ is a compact operator, so $T_0$ is eventually compact. For this particular combination of $X$ and $T_0$ the abstract duality structure from \cref{sec:duality} can be constructed systematically and explicitly \cite[Section II.5]{diekmann1995delay}. We only  summarize the few facts that will be used in the sequel.

\begin{remark}[Notation]\label{rem:dotnotation}
For $\KK \in \{\RR, \CC\}$ let $\KK^n$ be the linear space of column vectors and let $\KKK{n}$ be the linear space of row vectors, both over $\KK$. Elements of $\KK^n$ are denoted by $q = (q_1,q_2,\ldots,q_n)$ - commas between the entries - while elements in $\KKK{n}$ are denoted by $p = (p_1 ~ p_2 ~ \cdots ~ p_n)$ - no commas between the entries. We sometimes use the pairing defined by the row-column matrix multiplication:
\[
p \cdot q \DEF pq = \sum_{j=1}^n p_jq_j, \qquad p \in \KKK{n},\, q \in \KK^n.
\]
Note that the standard Hermitian inner product between two vectors $p^T,q \in \CC^n$ should be written as $\bar{p} \cdot q$ and {\em not} as $p\cdot q$.
\end{remark}

\begin{description}[wide]
\item[On $\SUN{X}$:]
  The maximal domain of strong continuity of $\STAR{T_0}$ has the representation
  \begin{equation}
    \label{eq:xsun_dde}
    \SUN{X} = \RRR{n} \times L^1([0,h],\RRR{n}),
  \end{equation}
  and the duality pairing between $\SUN{\phi} = (c,g) \in \SUN{X}$ and $\phi \in X$ is
  \begin{equation}
    \label{eq:pairing_X_sun_X}
    \PAIR{\SUN{\phi}}{\phi} = c\phi(0)+ \int_{0}^{h}g(\theta)\,\phi(-\theta)\,d\theta.
  \end{equation}
\item[On $\SUNSTAR{X}$:]
  Switching to the dual space of \cref{eq:xsun_dde} yields the representation
  \[
    \SUNSTAR{X} = \RR^n \times L^\infty([-h,0], \RR^n),
  \]
  and the duality pairing between $\SUNSTAR{\phi} = (a,\psi) \in \SUNSTAR{X}$ and $\SUN{\phi} = (c,g) \in \SUN{X}$ is
  \begin{equation}
    \label{eq:pairing_X_sun_star_X_sun}
    \PAIR{\SUNSTAR{\phi}}{\SUN{\phi}} = ca + \int_{0}^{h}g(\theta)\,\psi(-\theta)\,d\theta.
  \end{equation}
  The canonical injection \cref{eq:j} sends $\phi \in X$ to $j\phi = (\phi(0), \phi)$,  mapping $X$ \emph{onto} $\SUNSUN{X}$. Therefore $X$ is sun-reflexive with respect to the shift semigroup $T_0$.
\end{description}

Next, we specify the linear and nonlinear perturbations $L$ and $R$ in \cref{eq:T_T0_AIE,eq:aie}, respectively, and to relate these two abstract integral equations in $X$ to the linear and nonlinear initial value problems for \cref{eq:DDE}. For $i = 1,\ldots,n$ we denote $\rss_i \DEF (e_i, 0)$ where $e_i$ is the $i$th standard basis vector of $\RR^n$. It is conventional and convenient to introduce the shorthand
\[
  w \rss \DEF \sum_{i=1}^n{w_i\rss_i}, \qquad \forall\,w = (w_1,\ldots,w_n) \in \RR^n,
\]
and we note that $w\rss = (w, 0) \in \SUNSTAR{X}$. First we define the compact linear perturbation in \cref{eq:T_T0_AIE} as
\begin{equation}
  \label{eq:L}
  L\phi \DEF \PAIR{\zeta}{\phi}\rss,
\end{equation}
where the pairing in the right-hand side is given by \cref{eq:lindde_shorthand}. Now \cref{eq:lindde} with \cref{eq:DDE-ic} is equivalent to \cref{eq:T_T0_AIE} with \cref{eq:L} in the following sense: If $T$ is the unique $\mathcal{C}_0$-semigroup on $X$ satisfying \cref{eq:T_T0_AIE} with \cref{eq:L} then $x^{\phi} : [-h,\infty) \to \RR^n$ defined by
\[
  x^{\phi}_0 \DEF \phi, \qquad x^{\phi}(t) \DEF (T(t)\phi)(0), \qquad \forall\,t \ge 0,
\]
is the unique global solution of \cref{eq:lindde} with \cref{eq:DDE-ic} and
\[
  x^{\phi}_t = T(t)\phi, \qquad \forall\,t \ge 0.
\]
It remains to specify the nonlinear perturbation $R$ in \cref{eq:aie} as
\begin{equation}
  \label{eq:Rss}
  R(\phi) \DEF G(\phi) \rss,
\end{equation}
where $G$ is the nonlinear operator appearing in the splitting \cref{eq:DDE-RHS}. Let $\Sigma$ as in \cref{eq:semiflow} be the nonlinear semiflow generated by the family of maximal solutions of \cref{eq:aie} with \cref{eq:Rss}. The equivalence between \crefrange{eq:DDE}{eq:DDE-ic} and \cref{eq:aie} with \cref{eq:Rss} can be formulated as follows \cite[Prop. VII.6.1]{diekmann1995delay}. The function $x^{\phi} : [-h, t_{\phi}) \to \RR^n$ defined by
\[
  x^{\phi}_0 \DEF \phi, \qquad x^{\phi}(t) \DEF \Sigma(t,\phi)(0), \qquad \forall\,t \in I_{\phi},
\]
is the \emph{maximal solution} of \crefrange{eq:DDE}{eq:DDE-ic}, in the sense that any other solution necessarily exists only on a subinterval $[-h,t_+)$ for some $0 < t_+ \le t_{\phi}$ and coincides with $x^{\phi}$ there. Moreover,
\[
  x^{\phi}_t = \Sigma(t,\phi), \qquad \forall\,t \in I_{\phi}.
\]
It is the content of \cref{eq:linsigma} that generation and linearization commute: Starting with \cref{eq:DDE}, linearization of the semiflow $\Sigma$ at the stationary point $0 \in X$ yields precisely the eventually compact $\mathcal{C}_0$-semigroup $T$ corresponding to the linearized DDE \cref{eq:lindde}.

\subsection{Spectral computations for classical DDEs}
 The eventual compactness of $T$ implies that the spectrum of its generator $A$ - see \cref{rem:complex} - consists entirely of isolated eigenvalues of finite algebraic multiplicity. It is clear from \cref{eq:L} that $L$ is not just compact, but actually of finite rank. This implies that all spectral information about $A$ is contained in a holomorphic \emph{characteristic matrix function} $\Delta : \CC \to \CC^{n \times n}$ defined by
\begin{equation}
\label{eq:CharMatrix}
  \Delta(z) \DEF zI - \hat{\zeta}(z) \qquad \text{with} \qquad \hat{\zeta}(z) \DEF \int_0^h{e^{-z\theta}\,d\zeta(\theta)},
\end{equation}
where $\zeta$ is the \emph{real} kernel from \cref{eq:lindde_shorthand} \cite[Sections IV.4 and IV.5]{diekmann1995delay}. In particular, the eigenvalues of $A$ are the roots of the \emph{characteristic equation}
\begin{equation}
  \label{eq:main:det_delta}
\DET{\Delta(z)} = 0,
\end{equation}
and the algebraic multiplicity of an eigenvalue equals its order as a root of \cref{eq:main:det_delta}.

We will be concerned exclusively with \emph{simple} eigenvalues, for which the geometric and algebraic multiplicities are both equal to one. Let $\lambda \in \CC$ be such a simple eigenvalue of $A$. There exist nonzero right and left null vectors  $q \in \CC^n$ and $p \in \CCC{n}$ of $\Delta(\lambda)$,
\[
  \Delta(\lambda)q = 0, \qquad p\Delta(\lambda) = 0.
\]
The second equation is of course equivalent to $p^T$ being a nonzero right null vector of $\Delta^T(\lambda)$. The one-dimensional eigenspaces of $A$ and $\STAR{A}$ corresponding to $\lambda$ are spanned by eigenfunctions $\phi$ and $\SUN{\phi}$, respectively, with
\begin{equation}
  \label{eq:eigenfunction}
  \phi(\theta) = e^{\lambda\theta}q, \quad \theta \in [-h,0],
\end{equation}
and
\begin{equation}
  \label{eq:eigenfunction1}
  \SUN{\phi} = \left(p, \theta \mapsto p\int_{\theta}^{h}e^{\lambda(\theta-\tau)}\,d\zeta(\tau)\right), \quad \theta \in [0,h].
\end{equation}
We note that we have implicitly used - and will use consistently - the complexifications of $X$ and of the representation \cref{eq:xsun_dde} of $\SUN{X}$. For a simple eigenvalue $\lambda$,
\[
  \PAIR{\SUN{\phi}}{\phi} \neq 0,
\]
where the duality pairing is understood to be the complexification of \cref{eq:pairing_X_sun_X}. This nonequality implies that the eigenfunctions can be normalized to satisfy $\PAIR{\SUN{\phi}}{\phi} = 1$. In fact, from \cref{eq:pairing_X_sun_X,eq:eigenfunction} one computes
\begin{equation}
  \label{eq:normpair}
  \PAIR{\SUN{\phi}}{\phi} = p\Delta'(\lambda)q,
\end{equation}
so this normalization can be effectuated by scaling $p$ and $q$ such that $p\Delta'(\lambda)q = 1$. Finally, it is easily seen that if $\mu \neq \lambda$ is another simple eigenvalue of $A$ with eigenvector $\psi$ and adjoint eigenvector $\SUN{\psi}$, then
\begin{equation}
  \label{eq:zeropair}
  \PAIR{\SUN{\phi}}{\psi} = 0, \qquad \PAIR{\SUN{\psi}}{\phi} = 0.
\end{equation}

\subsection{Solvability of linear operator equations} \label{sec:solvability}
When computing the normal form coefficients in \cref{sec:Coefficients-of-parameter} using \cref{eq:homological_equation}, we will encounter linear operator equations of the form
\begin{equation}
  \label{eq:general_system_sunstar}
  (z I-\SUNSTAR{A})(v_0,v) = (w_0,w),
\end{equation}
where $z$ is a complex number, $(w_0,w) \in \SUNSTAR{X}$ is given and $(v_0,v) \in D(\SUNSTAR{A})$ is the unknown. In general, both $z$ and the right-hand side will have a nontrivial imaginary part, so here and from here onward, it is necessary to regard systems of the form \cref{eq:general_system_sunstar} as the complexification of the original operator equations. We will however not attach additional subscripts to the operator symbols, hoping that this omission will not cause confusion.

Since $\sigma(A)$ consists exclusively of point spectrum, there are two situations to consider depending on whether or not $z$ is an eigenvalue. If $z$ is \emph{not} an eigenvalue of $A$ then $z$ belongs to the resolvent set $\rho(A)$ of $A$ and \cref{eq:general_system_sunstar} admits a unique solution,
\[
  (v_0,v) =  \left(z I-\SUNSTAR{A}\right)^{-1}(w_0,w).
\]
In order to actually find this solution, one needs a representation of the resolvent operator of $\SUNSTAR{A}$. The general result can be found in \cite[Corollary IV.5.4]{diekmann1995delay}, but here we only require a special case.
\begin{lemma}\label{lem:regular_solution}
  Suppose that $z$ is not an eigenvalue of $A$, so \cref{eq:general_system_sunstar} has a unique solution $(v_0, v)$. If the right-hand side is represented by
  \[
    (w_0, w) = \left(w_0, \theta \mapsto  e^{z\theta}\Delta^{-1}(z)\eta\right),
  \]
  for some fixed vector $\eta \in \CC^n$, then this solution has the representation
  \[
    v_0 = v(0), \qquad v(\theta) = \Delta^{-1}(z)\left(e^{z\theta}w_0  + \left(\Delta'(z) - I -\theta\Delta(z)\right)w(\theta) \right).
  \]
\end{lemma}
\begin{proof}
 Write $(w_0, w) = (w_0, 0) + (0, \theta \mapsto  e^{z\theta}\Delta^{-1}(z)\eta)$, use the linearity of $(zI - \SUNSTAR{A})^{-1}$ and apply both cases of \cite[Corollary 3.4]{Janssens:Thesis}.
\end{proof}
\par
On the other hand, suppose that $z = \lambda$ is an eigenvalue. Then \cref{eq:general_system_sunstar} need not be consistent. In fact, a solution exists if and only if
\begin{equation}\tag{FSC}
  \label{eq:FSC}
  \PAIR{(w_0,w)}{\SUN{\phi}} = 0, \qquad \forall\,\SUN{\phi} \in \mathcal{N}(\lambda I - \STAR{A}),
\end{equation}
see \cite[Lemma 3.2]{Janssens:Thesis}). This condition is often referred to as the \emph{Fredholm solvability condition}. We note that the duality pairing in \cref{eq:FSC} may be evaluated in concrete cases using \cref{eq:eigenfunction} and the complexification of \cref{eq:pairing_X_sun_star_X_sun}. This will be done many times in \cref{sec:Coefficients-of-parameter} when we apply \cref{eq:FSC} to specific operator equations.
\par
If $z = \lambda$ is an eigenvalue and \cref{eq:general_system_sunstar} is consistent, then clearly its solutions are not unique. The bordered \emph{operator} inverse
\[
\INV{(\lambda I - \SUNSTAR{A})} : \mathcal{R}(\lambda I - \SUNSTAR{A}) \to \DOM(\SUNSTAR{A}),
\]
is used to select a particular solution in a systematic and convenient way. For the case that $\lambda$ is a \emph{simple} eigenvalue, it assigns the unique solution of the extended linear system
\begin{equation}
  \label{eq:bordop}
(\lambda I - \SUNSTAR{A})(v_0,v) = (w_0,w), \qquad \PAIR{(v_0,v)}{\SUN{\phi}} = 0,
\end{equation}
to every $(w_0,w)$ for which \cref{eq:general_system_sunstar} is consistent. The following lemma gives an explicit representation for a special case \cite[Proposition 3.6 and Corollary 3.7]{Janssens:Thesis}.
\begin{lemma}\label{lem:bordered}
  Let $z = \lambda$ be a simple eigenvalue with eigenvector $\phi$ and adjoint eigenvector $\SUN{\phi}$ as in \cref{eq:eigenfunction}, normalized to $\PAIR{\SUN{\phi}}{\phi} = 1$. Suppose \cref{eq:general_system_sunstar} is consistent for a given right-hand side of the form
  \[
    (w_0,w) = (\eta,0) + \kappa (q, \phi),
  \]
  where $\eta \in \CC^n$ and $\kappa \in \CC$. Then the unique solution $(v_0,v)$ of \cref{eq:bordop} is given by
  \[
    v_0 = \xi + \gamma q, \qquad v(\theta) = e^{\lambda\theta}(v_0 - \kappa \theta q),
  \]
  with $\xi = \INV{\Delta}(\lambda)(\eta + \kappa\Delta'(\lambda)q)$ and $\gamma = -p\Delta'(\lambda)\xi + \frac{1}{2}\kappa p \Delta''(\lambda)q$.
\end{lemma}
In \cref{sec:Coefficients-of-parameter} we will use the shorthand notation
\[
  v = \BINV{\lambda}(\eta,\kappa),
\]
for the solution in \cref{lem:bordered}. We observe that the expression for $\xi$ itself involves a bordered {\em matrix} inverse,
\[
\INV{\Delta}(\lambda) : \mathcal{R}(\Delta(\lambda)) \to \CC^n,
\]
which assigns the unique solution of the extended linear system
\[
  \Delta(\lambda)x = y, \qquad p\cdot x = 0,
\]
to every $y \in \CC^n$ for which the system $\Delta(\lambda)x = y$ is consistent - also see \cref{rem:dotnotation} for the notation. In practice, $x=\INV{\Delta}(\lambda)y$ can be obtained by solving the nonsingular bordered \emph{matrix} system
\[
\left(\begin{array}{cc}
\Delta(\lambda) & q\\
p & 0
\end{array}\right)\left(\begin{array}{c}
x\\
s
\end{array}\right)=\left(\begin{array}{c}
y\\
0
\end{array}\right),
\]
for the unknown $(x,s) \in \CC^{n+1}$ that necessarily satisfies $s = 0$. The properties of (finite dimensional) bordered linear systems and their role in numerical bifurcation analysis are discussed more extensively in \cite{Keller1987Numerical} and \cite[Chapter 3]{govaerts2000numerical}.

\section{Parameter dependence and classical DDEs}\label{sec:pd}
In \cref{sec:pd:motivation} we motivate our approach by explaining why the standard literature result does not apply to the problem at hand. This is most easily done at the concrete level of classical DDEs. The structure of the remaining subsections parallels that of \cref{sec:sunstar}. Namely, we first solve the problem of parameter dependence at the more abstract level of dual perturbation theory. In the final \cref{sec:pd:ddes} we then return to classical DDEs to see how the general results apply in this special case.

\subsection{Motivation}\label{sec:pd:motivation}
We are concerned with the situation where the right-hand side of \cref{eq:DDE} depends explicitly on parameters. Specifically, we consider
\begin{equation}
  \label{eq:pd-DDE}
  \dot{x}(t)= F(x_t, \alpha), \qquad t \ge 0,
\end{equation}
where $F: X \times \RR^p \to \RR^n$ is $C^k$-smooth for some $k \ge 1$ with $F(0,0) = 0$. We assume that at the critical parameter value $\alpha = 0$ the linearization of \cref{eq:pd-DDE} has $1 \le n_0 < \infty$ purely imaginary eigenvalues, counting multiplicities. The goal of \cref{sec:pd} is to obtain a parameter-dependent family of local center manifolds for a class of evolution equations that includes \cref{eq:pd-DDE}.

In \cite[Section IX.9.1]{diekmann1995delay} this problem is approached as follows. One augments \cref{eq:pd-DDE} with a trivial equation for the constant parameter dynamics. This gives the system
\begin{equation}
  \label{eq:pd:augmented}
  \left\{
    \begin{aligned}
      \dot{x}(t) &= F(x_t,\mu(t)),\\
      \dot{\mu}(t) &= 0,
    \end{aligned}
  \right.
  \qquad t \ge 0,
\end{equation}
on the state space $\BM{X} \DEF X \times \RR^p$, with $X = C([-h,0],\RR^n)$ as before in \cref{sec:ddecase}. Then the right-hand side of the first equation of \cref{eq:pd:augmented} is split as
\begin{equation}
  \label{eq:pd:splitting_book}
  F(\phi,\alpha) = D_1F(0,0)\phi + \tilde{G}(\phi,\alpha),
\end{equation}
which defines $\tilde{G} : \BM{X} \to \RR^n$, cf. \cite[(9.7) in Section IX.9.1]{diekmann1995delay}. The first term on the right of \cref{eq:pd:splitting_book} acts only on the $X$-component of the state in $\BM{X}$, so the semigroup $\tilde{\BM{T}}$ on $\BM{X}$ obtained by perturbing the shift-semigroup $\BM{T}_0$ is \emph{diagonal}.

However, there is an obstruction. In order to satisfy the hypotheses of the parameter-\emph{in}dependent center manifold theorem, $\tilde{G}$ must be a pure nonlinearity on $\BM{X}$, i.e.
\[
\tilde{G}(0,0) = 0, \qquad D_1\tilde{G}(0,0) = 0, \qquad D_2\tilde{G}(0,0) = 0.
\]
The first two of these conditions are clearly fulfilled, but in general there is no reason for the third condition to be met. It \emph{does} hold when $\tilde{G}(0,\alpha) = 0$ for all $\alpha \in \RR^p$ in a neighborhood of zero, i.e. when the zero equilibrium of \cref{eq:pd-DDE} persists under small parameter variations. For a generic fold-Hopf bifurcation - as well as for a generic Bogdanov-Takens bifurcation that we do not discuss here - \emph{there is no such persistence}.

In this article, the above difficulty is resolved by considering instead of \cref{eq:pd:splitting_book} the splitting
\begin{equation}
  \label{eq:pd:splitting}
  F(\phi,\alpha) = D_1F(0,0)\phi + D_2F(0,0)\alpha + G(\phi,\alpha).
\end{equation}
Using this splitting, \cref{eq:pd:augmented} is written as
\begin{equation}
  \label{eq:pd-DDE12}
  \left\{
    \begin{aligned}
      \dot{x}(t) &= D_1F(0,0)x_t + D_2F(0,0)\mu(t) + G(x_t, \mu(t)),\\
      \dot{\mu}(t) &= 0,
    \end{aligned}
  \right.
  \qquad t \ge 0.
\end{equation}
Now \emph{both} $D_1F(0,0)$ \emph{as well as} $D_2F(0,0)$ appear in the perturbation of $\BM{T}_0$. As a consequence the perturbed semigroup $\BM{T}$ is no longer diagonal, but still simple enough for a complete analysis. Moreover,
\begin{equation}
  \label{eq:G_conditions}
  G(0,0) = 0, \qquad D_1G(0,0) = 0, \qquad D_2G(0,0) = 0,
\end{equation}
so the  parameter-independent center manifold theorem can be applied without having to assume equilibrium persistence. Of course $G = \tilde{G}$ and $\BM{T} = \tilde{\BM{T}}$ whenever $D_2F(0,0) = 0$.

\begin{remark}\label{rem:more_general}
 In a first attempt we regarded the augmented system \cref{eq:pd:augmented} as a DDE on the state space $C([-h,0], \RR^{n + p})$, also see \textup{\cite{GuoMan2011parCM}}, but we found this approach to be a bit unsatisfactory: The proofs in \cref{sec:pd:duality,sec:spectral_bold,sec:nonlinear_bold,sec:cm_bold} do not depend on the details of the class of equations under consideration, while those same details sometimes lead to notation that is more complicated than necessary.
\end{remark}

\subsection{Duality structure and linear perturbation}\label{sec:pd:duality}
We work in the setting of \cref{sec:duality}. Namely, let $T_0$ be a $\mathcal{C}_0$-semigroup on a real or complex Banach space $X$ that is sun-reflexive with respect to $T_0$. We write $\KK \in \{\RR, \CC\}$ for the underlying scalar field - as in \cref{rem:dotnotation}. Define $\BM{T}_0$ on $\BM{X}$ by
\begin{equation}
  \label{eq:T0_bold}
  \BM{T}_0(t) \DEF \diag{(T_0(t), I_p)}.
\end{equation}
The procedure of taking adjoints and restrictions (twice) then yields semigroups $\STAR{\BM{T}_0}$, $\SUN{\BM{T}_0}$, $\SUNSTAR{\BM{T}_0}$ and $\SUNSUN{\BM{T}_0}$ on $\STAR{\BM{X}} \simeq \STAR{X} \times \KK^p$, $\SUN{\BM{X}} \simeq \SUN{X} \times \KK^p$, $\SUNSTAR{\BM{X}} \simeq \SUNSTAR{X} \times \KK^p$ and $\SUNSUN{\BM{X}} \simeq \SUNSUN{X} \times \KK^p$. (The symbol $\simeq$ indicates an identification via a natural isometric isomorphism.) It is straightforward to check that on $\SUNSTAR{\BM{X}}$ we have
\begin{equation}
  \label{eq:T0_sun_star_bold}
  \SUNSTAR{\BM{T}_0}(t) \ = \diag{(\SUNSTAR{T_0}(t), I_p)},
\end{equation}
and that the canonical injection $\BM{j} : \BM{X} \to \SUNSTAR{\BM{X}}$ has the form
\begin{equation}
  \label{eq:j_bold}
  \BM{j} = \diag{(j, I_p)},
\end{equation}
where $j : X \to \SUNSTAR{X}$ is the canonical injection from \cref{eq:j}. In particular, $\BM{X}$ is sun-reflexive with respect to $\BM{T}_0$.

As in \cref{sec:duality} we now introduce a bounded linear perturbation $\BM{L} : \BM{X} \to \SUNSTAR{\BM{X}}$ of $\BM{T}_0$. We let it be of the form
\begin{equation}
  \label{eq:B_bold}
  \BM{L} =
  \begin{pmatrix}
    L& L_p\\
    0& 0
  \end{pmatrix},
\end{equation}
with $L : X \to \SUNSTAR{X}$ and $L_p : \KK^p \to \SUNSTAR{X}$ bounded linear operators. Perturbing $T_0$ by $L$ and $\BM{T}_0$ by $\BM{L}$ yields $\mathcal{C}_0$-semigroups $T$ on $X$ and $\BM{T}$ on $\BM{X}$, respectively. Let $A$ and $\BM{A}$ be their generators.

\begin{remark}
  There are at least two equivalent ways to compute $\BM{T}$ and $\BM{A}$ on $\BM{X}$ and their {\WSTAR} counterparts $\SUNSTAR{\BM{T}}$ and $\SUNSTAR{\BM{A}}$ on $\SUNSTAR{\BM{X}}$. One approach - suggested to us by Odo Diekmann - uses integrated semigroup theory to calculate first $\BM{T}$ and next $\SUNSTAR{\BM{T}}$. Then $\SUNSTAR{\BM{A}}$ and $\BM{A}$ are calculated, in that order.

  Here we go the other way around: We start by calculating $\SUNSTAR{\BM{A}}$ and use it to obtain $\SUNSTAR{\BM{T}}$. If desired $\BM{A}$ and $\BM{T}$ can then be found by application of \cref{eq:j_bold} and its inverse. This approach is more elementary - we use only theory that was already introduced in \cref{sec:duality} - and it yields the same outcome, as it should.
\end{remark}

\begin{proposition}\label{prop:A_sun_star_bold}
  The {\WSTAR} generator $\SUNSTAR{\BM{A}}$ of $\SUNSTAR{\BM{T}}$ has the representation
  \[
    \DOM(\SUNSTAR{\BM{A}}) = \DOM(\SUNSTAR{A})\times \KK^p, \qquad \SUNSTAR{\BM{A}} =
    \begin{pmatrix}
      \SUNSTAR{A}& L_p\\
      0& 0
    \end{pmatrix},
  \]
  where $\SUNSTAR{A}$ is the {\WSTAR} generator of $\SUNSTAR{T}$.
\end{proposition}
\begin{proof}
  According to the general theory of \cref{sec:duality} and \cref{eq:A_sunstar} in particular, we have
  \[
    \DOM(\SUNSTAR{\BM{A}}) = \DOM(\SUNSTAR{\BM{A}_0}), \qquad \SUNSTAR{\BM{A}} = \SUNSTAR{\BM{A}_0} + \BM{L}\BM{j}^{-1},
  \]
  From \cref{eq:T0_sun_star_bold} we see that
  \[
    \DOM(\SUNSTAR{\BM{A}_0}) = \DOM(\SUNSTAR{A_0})\times \KK^p, \qquad \SUNSTAR{\BM{A}_0} = \diag{(\SUNSTAR{A_0}, 0)}.
  \]
  Using \cref{eq:j_bold,eq:B_bold} we calculate
  \[
    \SUNSTAR{\BM{A}_0} + \BM{L}\BM{j}^{-1} =
    \begin{pmatrix}
      \SUNSTAR{A_0}& 0\\
      0& 0
    \end{pmatrix} +
    \begin{pmatrix}
      L& L_p\\
      0& 0
    \end{pmatrix}
    \begin{pmatrix}
      j^{-1}& 0\\
      0& I_p
    \end{pmatrix},
  \]
  and the result follows.
\end{proof}

\begin{lemma}\label{lem:Lp}
  Let $\SUN{L_p} : \SUN{X} \to \KK^p$ be the restriction of $\STAR{L_p}$ to $\SUN{X}$. Then $\SUNSTAR{L_p} = L_p$.
\end{lemma}
\begin{proof}
  We begin by noting that - strictly speaking - this statement involves two canonical identifications. Namely, let $i : \SUN{X} \to \SUNSTARSTAR{X}$ and $i_p : \KK^p \to \STARSTAR{{\KK^p}}$ be the canonical injection and bijection, respectively. Then $\SUN{L_p} \DEF \STAR{L_p} i$ and we need to prove that $\SUNSTAR{L_p}i_p = L_p$. For this it is not difficult to show that
  \[
    \PAIR{\SUNSTAR{L_p}i_p\alpha}{\SUN{\phi}} = \PAIR{L_p\alpha}{\SUN{\phi}},
  \]
  for all $\alpha \in \KK^p$ and for all $\SUN{\phi} \in \SUN{X}$.
\end{proof}

For the purpose of notation, we define the \emph{integrated semigroup} $\SUNSTAR{W}$ for  $\SUNSTAR{T}$ as
\begin{equation}\label{eq:integrated_semigroup}
  \SUNSTAR{W}(t)\SUNSTAR{\phi} \DEF \int_0^t{\SUNSTAR{T}(\tau)\SUNSTAR{\phi}\,d\tau}, \qquad t \ge 0.
\end{equation}
with on the right a {\WSTAR} Riemann integral of the same type as the integral in \cref{eq:T_T0_AIE}.

\begin{proposition}\label{prop:T_bold_sun_star}
  The semigroup $\SUNSTAR{\BM{T}}$ that is {\WSTARLY} generated by $\SUNSTAR{\BM{A}}$ has the representation
  \begin{equation}\label{eq:T_bold_sun_star}
    \SUNSTAR{\BM{T}}(t) =
    \begin{pmatrix}
      \SUNSTAR{T}(t)& \SUNSTAR{W}(t)L_p\\
      0& I_p
    \end{pmatrix}, \qquad t \ge 0.
  \end{equation}
\end{proposition}
\begin{proof}
  We define a one-parameter family $\BM{S}$ of bounded linear operators on $\SUN{\BM{X}}$ by
  \[
    \BM{S}(t) =
    \begin{pmatrix}
      \SUN{T}(t)& 0\\
      \SUN{L_p}\SUN{W}(t)& I_p
    \end{pmatrix}, \qquad t \ge 0,
  \]
  where
  \[
    \SUN{W}(t)\SUN{\phi} \DEF \int_0^t{\SUN{T}(\tau)\SUN{\phi}\,d\tau}, \qquad t \ge 0.
  \]
  It is easy to check that $\BM{S}$ is a $\mathcal{C}_0$-semigroup on $\SUN{\BM{X}}$. By \cref{lem:Lp} the adjoint semigroup $\STAR{\BM{S}}(t)$ equals the right-hand side of \cref{eq:T_bold_sun_star} for all $t \ge 0$. We will show that the {\WSTAR} generator $\SUNSTAR{\BM{A}}$ of $\SUNSTAR{\BM{T}}$ is also the {\WSTAR} generator of $\STAR{\BM{S}}$. This will then imply that $\SUNSTAR{\BM{T}} = \STAR{\BM{S}}$.
  \par
We use \cref{prop:A_sun_star_bold}. Let $\BM{C}$ be the generator of $\BM{S}$, so $\STAR{\BM{C}}$ is the {\WSTAR} generator of $\STAR{\BM{S}}$. For any $(\SUNSTAR{\phi}, \alpha)$ in $\SUNSTAR{\BM{X}}$ and any $t > 0$ we have
\[
  \frac{1}{t}\left(\STAR{\BM{S}}(t)(\SUNSTAR{\phi}, \alpha) - (\SUNSTAR{\phi}, \alpha) \right) = \frac{1}{t}
  \begin{pmatrix}
    \SUNSTAR{T}(t)\SUNSTAR{\phi} - \SUNSTAR{\phi}\\
    0
  \end{pmatrix} +
  \frac{1}{t}
  \begin{pmatrix}
    \SUNSTAR{W}(t)L_p\alpha\\
    0
  \end{pmatrix}.
\]
We note that $t^{-1}\SUNSTAR{W}(t)L_p\alpha \to L_p\alpha$ {\WSTARLY} as $t \downarrow 0$. It follows that the right-hand side converges {\WSTARLY} if and only if $\SUNSTAR{\phi} \in D(\SUNSTAR{A})$ and in that case the {\WSTAR}-limit equals $(\SUNSTAR{A}\SUNSTAR{\phi} + L_p\alpha, 0) = \SUNSTAR{\BM{A}}(\SUNSTAR{\phi}, \alpha)$. We conclude that $\STAR{\BM{C}} = \SUNSTAR{\BM{A}}$.
\end{proof}

\subsection{Spectral theory and the center eigenspace}\label{sec:spectral_bold}
Let $T_0$ be a $\mathcal{C}_0$-semigroup on a complex Banach space $X$ that is sun-reflexive with respect to $T_0$. For the purpose of spectral theory, we explicitly take $\CC$ as the underlying scalar field. In examples, $X$ will often be a complexification of a real Banach space, see \cref{rem:complex}.

We are interested in a description of the spectrum and the corresponding (generalized) eigenspaces of the generator $\BM{A}$ of $\BM{T}$. In particular, \cref{prop:compact_bold,prop:real_center_subspace_bold} below guarantee, respectively, the existence and smooth parametrization of the parameter-dependent local center manifold in \cref{sec:cm_bold}. 

\begin{proposition}\label{prop:spectrum_bold}
  The spectrum $\sigma(\SUNSTAR{\BM{A}}) = \sigma(\SUNSTAR{A}) \cup \{0\}$ with resolvent operator
   \begin{equation}
    \label{eq:resolvent_bold}
    R_z(\SUNSTAR{\BM{A}}) =
    \begin{pmatrix}
      R_z(\SUNSTAR{A})& z^{-1}R_z(\SUNSTAR{A})L_p\\
      0 & z^{-1}I_p
    \end{pmatrix},
  \end{equation}
  for every $z$ in the resolvent set $\rho(\SUNSTAR{\BM{A}})$.
\end{proposition}
\begin{proof}
  From \cref{prop:A_sun_star_bold} we have
  \[
    z I - \SUNSTAR{\BM{A}} =
    \begin{pmatrix}
      z I - \SUNSTAR{A}& -L_p\\
      0& z I_p
    \end{pmatrix}.
  \]
  This upper triangular operator matrix has a bounded inverse if and only if \emph{both} entries on its diagonal have bounded inverses, which happens if and only if $z \in \rho(\SUNSTAR{A})$ \emph{and} $z \neq 0$. In that case, the inverse is given precisely by the stated expression for $R_z(\SUNSTAR{\BM{A}}) \DEF (z I - \SUNSTAR{\BM{A}})^{-1}$.
\end{proof}

\emph{In addition we assume that $T_0$ is eventually compact and the perturbation $L$ in \cref{eq:B_bold} is compact.} As a consequence, the spectral analysis of $\SUNSTAR{\BM{A}}$ reduces to an analysis of the poles of its resolvent operator \cite[Corollary V.3.2]{Engel2000}, \cite[Section V.10]{Taylor1980}.

\begin{proposition}\label{prop:compact_bold}
  $\BM{T}$ is an eventually compact $\mathcal{C}_0$-semigroup.
\end{proposition}
\begin{proof}
The eventual compactness of $T_0$, the finite rank of $I_p$ and \cref{eq:T0_bold} together imply that $\BM{T}_0$ is eventually compact. Since $L_p$ has finite rank and $L$ is compact by assumption, it follows from \cref{eq:B_bold} that $\BM{L}$ is compact, so $\BM{T}$ is eventually compact by \cite[Theorem 2.8]{diekmann2007stability}.
\end{proof}

\begin{theorem}
  \label{thm:eigspaces_bold}
  The generalized eigenspace corresponding to $\lambda \in \sigma(\SUNSTAR{\BM{A}})$ is given by
  \[
    \mathcal{M}_{\lambda}(\SUNSTAR{\BM{A}}) =
    \begin{cases}
      \mathcal{M}_{\lambda}(\SUNSTAR{A}) \times \{0\}, &\text{if } \lambda \neq 0,\\
      \mathcal{M}_0(\SUNSTAR{A}) \times \{0\} \oplus
    \left\{ \left(\Gamma_0L_p\alpha, \alpha\right) \,:\, \alpha \in \CC^p \right\}, &\text{if } \lambda = 0,
    \end{cases}
  \]
  where $\Gamma_0$ is a bounded linear operator on $\SUNSTAR{X}$ mapping into $\SUNSUN{X}$.
\end{theorem}
\begin{proof}
  Let $\lambda \in \sigma(\SUNSTAR{\BM{A}})$ be arbitrary. Taking residues at $z = \lambda$ in \cref{eq:resolvent_bold}, we obtain
  \begin{equation}
    \label{eq:P_bold}
    \SUNSTAR{\BM{P}_{\lambda}} =
    \begin{pmatrix}
      \SUNSTAR{P_{\lambda}}& \Gamma_{\lambda}L_p\\
      0& I_p\delta_{\lambda}
    \end{pmatrix}, \qquad \Gamma_{\lambda} \DEF \RES_{z=\lambda}{z^{-1}R_z(\SUNSTAR{A})},
  \end{equation}
  where $\delta_{\lambda} \DEF \delta_{\lambda,0}$ is the Kronecker delta and $\SUNSTAR{\BM{P}_{\lambda}}$ and $\SUNSTAR{P_{\lambda}}$ are the spectral projectors corresponding to $\lambda$ for $\SUNSTAR{\BM{A}}$ and $\SUNSTAR{A}$. (If $\lambda$ is in the resolvent set of the respective operator, then the residue - hence the spectral projector - is identically zero.) $\Gamma_{\lambda}$ is pointwise equal to a contour integral with an integrand in the closed subspace $\SUNSUN{X}$ of $\SUNSTAR{X}$, so $\Gamma_{\lambda}$ maps into $\SUNSUN{X}$.

  We will now calculate the range of $\SUNSTAR{\BM{P}_{\lambda}}$ from \cref{eq:P_bold}. In general,
  \begin{equation}
    \label{eq:gen_eigenspace_bold}
    \mathcal{M}_{\lambda}(\SUNSTAR{\BM{A}}) =
    \left\{
      \begin{pmatrix}
        \SUNSTAR{P_{\lambda}}\SUNSTAR{\phi}\\
        0
      \end{pmatrix}
      +
      \begin{pmatrix}
        \Gamma_{\lambda}L_p\alpha\\
        \alpha\delta_{\lambda}
      \end{pmatrix}
      \,:\, (\SUNSTAR{\phi},\alpha) \in \SUNSTAR{\BM{X}}
    \right\}.
  \end{equation}
  First we assume that $\lambda \neq 0$, so $\delta_{\lambda} = 0$. We are going to show that
   \begin{equation}\label{eq:lambda_not_zero}
    \left\{\SUNSTAR{P_{\lambda}}\SUNSTAR{\phi} + \Gamma_{\lambda}L_p\alpha \,:\, (\SUNSTAR{\phi}, \alpha) \in \SUNSTAR{\BM{X}}\right\} = \mathcal{M}_{\lambda}(\SUNSTAR{A}).
  \end{equation}
   Together with \cref{eq:gen_eigenspace_bold} this will then prove the theorem for $\mathcal{M}_{\lambda}(\SUNSTAR{\BM{A}})$. To verify \cref{eq:lambda_not_zero} let $p \in \NN$ be the order of $\lambda$ as a pole of $R_z(\SUNSTAR{A})$. For $n = 1,\ldots,p$ let $B_n$ be the coefficient of $(z - \lambda)^{-n}$ in the Laurent series for $R_z(\SUNSTAR{A})$. A small computation shows that
  \begin{equation}
    \label{eq:Gamma}
    \Gamma_{\lambda} = \sum_{k=1}^p{(-1)^{k+1}\lambda^{-k}B_k}.
  \end{equation}
  From \cite[Section V.10]{Taylor1980} we recall the relation $B_{n+1} = (\SUNSTAR{A} - \lambda I)^n B_1$ for all $n \in \NN$. Since $B_1 = \SUNSTAR{P_{\lambda}}$ and its range $\mathcal{M}_{\lambda}(\SUNSTAR{A})$ is an invariant subspace of $\SUNSTAR{A}$, this relation implies that $B_k$ takes values in $\mathcal{M}_{\lambda}(\SUNSTAR{A})$ for all $k = 1,\ldots,p$, so the same is true for $\Gamma_{\lambda}$ by \cref{eq:Gamma}. From this it follows that \cref{eq:lambda_not_zero} holds.
  \par
  For the remaining case $\lambda = 0$ we have $\delta_{\lambda} = 1$, so from \cref{eq:gen_eigenspace_bold} we get the direct sum
  \[
    \mathcal{M}_0(\SUNSTAR{\BM{A}}) =
    \left\{
      \begin{pmatrix}
        \SUNSTAR{P_0}\SUNSTAR{\phi}\\
        0
      \end{pmatrix}
      \,:\, \SUNSTAR{\phi} \in \SUNSTAR{X}
    \right\}
    \oplus
    \left\{
      \begin{pmatrix}
        \Gamma_0L_p\alpha\\
        \alpha
      \end{pmatrix}
      \,:\, \alpha \in \CC^p
    \right\}.
  \]
  The first summand equals $\mathcal{M}_0(\SUNSTAR{A}) \times \{0\}$ and this gives the result.
\end{proof}

\begin{corollary}
  \label{cor:center_subspace_bold}
  The center eigenspace $\BM{X}_0$ corresponding to the purely imaginary eigenvalues of $\BM{A}$ is given by
  \[
    \BM{X}_0 = X_0 \times \{0\} \oplus \left\{ \left(j^{-1}\Gamma_0L_p\alpha, \alpha \right) \,:\, \alpha \in \CC^p \right\},
  \]
  with $\DIM{\BM{X}_0} = \DIM{X_0} + p$.
\end{corollary}
\begin{proof}
  By \cref{prop:spectrum_bold} we have the disjoint union $\sigma(\SUNSTAR{\BM{A}}) = (\sigma(\SUNSTAR{A}) \setminus \{0\}) \cup \{0\}$. Using this and \cref{thm:eigspaces_bold} we first compute the center eigenspace for $\SUNSTAR{\BM{A}}$ as
  \[
    \SUNSTAR{\BM{X}_0} = \SUNSTAR{X_0} \times \{0\} \oplus \{(\Gamma_0L_p\alpha,\alpha)\,:\,\alpha \in \CC^p\},
  \]
  and then we apply $\BM{j}^{-1}$ from \cref{eq:j_bold} to both sides of this equality.
\end{proof}

In \cref{sec:nonlinear_bold,sec:cm_bold,sec:pd:ddes} we will consider nonlinear problems on a real Banach space. In this case spectral analysis must be preceded by complexification, see \cref{rem:complex} and in particular \cite[last part of Section IV.2]{diekmann1995delay}.

\begin{proposition}
  \label{prop:real_center_subspace_bold}
  Suppose that $X = Y_{\CC}$ is a complexification of a real Banach space $Y$ and let $Y_0 \subseteq Y$ be the \emph{real} center eigenspace associated with $X_0$. Then the \emph{real} center eigenspace $\BM{Y}_0$ associated with $\BM{X}_0$ is
  \begin{equation}
    \label{eq:Y0bold}
    \BM{Y}_0 = Y_0 \times \{0\} \oplus \left\{\left(Q\alpha, \alpha \right) \,:\, \alpha \in \RR^p \right\} \subseteq Y \times \RR^{p},
  \end{equation}
  where $Q : \RR^p \to Y$ is a bounded linear operator. Furthermore, $\iota : \BM{Y}_0 \to Y_0 \times \RR^p$ defined by $\iota(\psi, \alpha) \DEF (\psi - Q\alpha, \alpha)$ is a linear isomorphism.
\end{proposition}
\begin{proof}
  $\BM{X}$ is naturally identified with $\BM{Y}_{\CC}$ where $\BM{Y} = Y \times \RR^p$. Let $\BM{P}_{\Lambda}$ with range $\BM{X}_0$ be the spectral projector on $\BM{X}$ for the spectral set $\Lambda = \overline{\Lambda}$ of all purely imaginary eigenvalues of $\BM{A}_{\CC}$. A direct generalization of \cite[Corollary IV.2.19]{diekmann1995delay} implies that $\BM{P}_{\Lambda}$ is the complexification of a projector $\BM{P}_{\Lambda}^Y$ on $\BM{Y}$ and the range $\BM{Y}_0$ of $\BM{P}_{\Lambda}^Y$ - identified with a subspace of $\BM{Y}$ - is the real center eigenspace for $\BM{A}$. Also, $\Gamma_0$ on $\SUNSTAR{X}$ is self-conjugate by \cref{eq:P_bold}. Together with \cref{cor:center_subspace_bold} this implies \cref{eq:Y0bold}. It is easily verified that the linear operator $\iota$ is an isomorphism.
\end{proof}

\begin{remark}\label{rem:X0Y0}
We will not make a \emph{notational} distinction between the real and complex center eigenspaces, indicating both $X_0$ and $Y_0$ with $X_0$ and both $\BM{X}_0$ and $\BM{Y}_0$ with $\BM{X}_0$, respectively. We hope that the underlying scalar field will be clear from the immediate context.
\end{remark}

\subsection{Nonlinear perturbation}\label{sec:nonlinear_bold}
Let $T_0$ be a $\mathcal{C}_0$-semigroup on a \emph{real} Banach space $X$ that is sun-reflexive with respect to $T_0$. Introduce a nonlinear perturbation $\BM{R} : \BM{X} \to \SUNSTAR{\BM{X}}$ of the form
\begin{equation}\label{eq:R_bold}
  \BM{R}(\phi,\alpha) = (R(\phi,\alpha), 0),
\end{equation}
where $R : \BM{X} \to \SUNSTAR{X}$ is $C^k$-smooth, satisfying
\begin{equation}\label{eq:R_conditions}
  R(0,0) = 0, \qquad D_1{R}(0,0) = 0, \qquad D_2{R}(0,0) = 0.
\end{equation}
We associate with $\BM{T}$ and $\BM{R}$ the integral equation
\begin{equation}
  \label{eq:aie_bold}
  \BM{u}(t) = \BM{T}(t - s)\BM{u}(s) + \BM{j}^{-1}\int_s^t{\SUNSTAR{\BM{T}}(t - \tau)\BM{R}(\BM{u}(\tau))\,d\tau}, \qquad -\infty < s \le t < \infty.
\end{equation}
We expect all nontrivial dynamics to be contained in the first component, and this is indeed the case:
\begin{proposition}\label{prop:aie_bold}
  The function $\BM{u} = (u, u_p) : I \to \BM{X}$ is a solution of \cref{eq:aie_bold} if and only if $u_p$ is constant on $I$ and $u : I \to X$ is a solution of
  \begin{equation}\label{eq:aie_bold_1st}
    u(t) = T(t - s)u(s) + j^{-1}\int_s^t{\SUNSTAR{T}(t - \tau)(L_p\alpha + R(u(\tau),\alpha))\,d\tau}, \qquad -\infty < s \le t < \infty,
  \end{equation}
  where $\alpha \in \RR^p$ denotes the constant value of $u_p$.
\end{proposition}
\begin{proof}
  We use \cref{prop:T_bold_sun_star,eq:j_bold}. For any continuous function $\BM{u} = (u, u_p) : I \to \BM{X}$ we compute
  \[
    \BM{T}(t - s)\BM{u}(s) = \BM{j}^{-1}\SUNSTAR{\BM{T}}(t - s)\BM{j}\BM{u}(s) =
    \begin{pmatrix}
      T(t - s)u(s) + j^{-1}\SUNSTAR{W}(t - s)L_p u_p(s)\\
      u_p(s)
    \end{pmatrix},
  \]
  while another computation shows that
  \[
    \BM{j}^{-1}\int_s^t{\SUNSTAR{\BM{T}}(t - \tau)\BM{R}(\BM{u}(\tau))\,d\tau} =
    \begin{pmatrix}
      j^{-1}\int_s^t{\SUNSTAR{T}(t - \tau)R(\BM{u}(\tau))\,d\tau}\\
      0
    \end{pmatrix}.
  \]
  From the above together with the definition \cref{eq:integrated_semigroup} of $\SUNSTAR{W}(t)$ we see that \cref{eq:aie_bold} is equivalent to the system
  \[
    \left\{
      \begin{aligned}
        u(t) &= T(t - s)u(s) + j^{-1}\int_s^t{\SUNSTAR{T}(t - \tau)(L_p u_p(s) + R(u(\tau),u_p(\tau)))\,d\tau},\\
        u_p(t) &= u_p(s),
      \end{aligned}
    \right.
  \]
   for $-\infty < s \le t < \infty$. The statement now follows.
\end{proof}

\subsection{Parameter-dependent local center manifolds}\label{sec:cm_bold}
We consider again a $\mathcal{C}_0$-semigroup $T_0$ on a \emph{real} Banach space $X$ that is sun-reflexive with respect to $T_0$. We also assume that $T_0$ is eventually compact and $L$ in \cref{eq:B_bold} is compact, so \cref{prop:compact_bold} implies that $\BM{T}$ is eventually compact. If furthermore the nonlinearity $\BM{R}$ satisfies \cref{eq:R_conditions}, then all conditions are fulfilled for the application of the center manifold theory from \cite[Chapter IX]{diekmann1995delay} to \cref{eq:aie_bold}.

Therefore, if the generator $A$ of $T$ has $1 \le n_0 < \infty$ purely imaginary eigenvalues, counting algebraic multiplicities, then by \cref{prop:real_center_subspace_bold} the \emph{real} center eigenspace $\BM{X}_0$ has dimension $n_0 + p$. There exists a $C^k$-smooth local center manifold $\CMBOLD$ in $\BM{X}$ that is tangent at the origin to $\BM{X}_0$. In fact, \cref{prop:real_center_subspace_bold} implies that $\CMBOLD$ is the image of a $C^k$-smooth map
\[
  \BM{\mathcal{C}} : U \times U_p \subseteq X_0 \times \RR^p \to \BM{X},
\]
where $U \subseteq X_0$ and $U_p \subseteq \RR^p$ are neighborhoods of the origin. Since \cref{eq:R_bold} has a zero in the second component, it follows from \cite[(5.1) in Section IX.5]{diekmann1995delay} that $\BM{\mathcal{C}}$ has the form
\begin{equation}
  \label{eq:cmmap_bold}
\BM{\mathcal{C}}(\phi, \alpha) = (\mathcal{C}(\phi, \alpha), \alpha), \qquad \forall\,(\phi, \alpha) \in U \times U_p,
\end{equation}
where $\mathcal{C} : U \times U_p \to X$ is the first component function.

\begin{definition}\label{def:cm_alpha}
  The image $\CM(\alpha) \DEF \mathcal{C}(U, \alpha)$ is a \emph{local center manifold for \cref{eq:aie_bold_1st} at $\alpha \in U_p$}.
\end{definition}

It is a direct consequence of the above definition that for every $\alpha \in U_p$ we can parametrize $\CM(\alpha)$ by coordinates on the real center eigenspace $X_0$ that depend $C^k$-smoothly on $\alpha$. This will be important for the discussion of the normalization method following \cref{eq:ODEonCM} in \cref{sec:normal-forms}.

\begin{proposition}
  If $\alpha \in U_p$ is sufficiently small then $\CM(\alpha)$ is locally positively invariant for the semiflow generated by \cref{eq:aie_bold_1st}.
\end{proposition}
\begin{proof}
  Let $\BM{\Sigma}$ and $\Sigma$ be the semiflows generated by \cref{eq:aie_bold,eq:aie_bold_1st}, respectively. By \cref{prop:aie_bold},
  \begin{equation}
    \label{eq:sigmasigma}
    \BM{\Sigma}(s, (\psi, \alpha)) = (\Sigma(s, \psi), \alpha), \qquad \forall\,\psi \in X,\alpha \in \RR^p,
  \end{equation}
  and for all $s$ in a common interval of existence $I_{\psi,\alpha}$. By \cite[Theorem IX.5.3(i)]{diekmann1995delay} there exists $\delta > 0$ such that if $(\psi, \alpha) \in \CMBOLD$ and if
  \[
    \|\BM{\Sigma}(s, (\psi, \alpha))\| = \|\Sigma(s, \psi)\| + |\alpha| \le \delta, \qquad \forall\,s \in [0,t],
  \]
  then $\BM{\Sigma}(t, (\psi, \alpha)) \in \CMBOLD$ which by \cref{eq:sigmasigma} implies that $\Sigma(t, \psi) \in \CM(\alpha)$.

  We note that if $\psi \in \CM(\alpha)$ then by \cref{eq:cmmap_bold} it follows that $(\psi, \alpha) \in \CMBOLD$. Therefore, if $|\alpha| \le \frac{\delta}{2}$ and $\psi \in \CM(\alpha)$ then
  \[
    \|\Sigma(s, \psi)\| \le \frac{\delta}{2}, \qquad \forall\,s \in [0,t]
  \]
  implies that $\Sigma(t, \psi) \in \CM(\alpha)$. This is precisely local positive invariance of $\CM(\alpha)$ for $\Sigma$.
\end{proof}

Next we consider a solution $u : I \to X$ of \cref{eq:aie_bold_1st} that lies in $\CM(\alpha)$. By \cref{prop:aie_bold} the function $\BM{u} = (u, \alpha) : I \to \BM{X}$ is a solution of \cref{eq:aie_bold}. Also, since $u$ lies in $\CM(\alpha)$ we see from \cref{eq:cmmap_bold} that $\BM{u}$ lies in $\CMBOLD$ and therefore satisfies the differential equation
\[
  \BM{j}\dot{\BM{u}}(t) = \SUNSTAR{\BM{A}}\BM{j}\BM{u}(t) + \BM{R}(\BM{u}(t)), \qquad \forall\,t \in I,
\]
cf. \cref{eq:aode}. By \cref{eq:j_bold,prop:A_sun_star_bold,eq:R_bold} the first component of this equation gives the differential equation
\begin{equation}
  \label{eq:ode_on_cm_abstract}
  j\dot{u}(t) = \SUNSTAR{A}ju(t) + L_p\alpha + R(u(t),\alpha), \qquad \forall\,t \in I,
\end{equation}
that is satisfied by $u$.\\
In summary,

\begin{theorem}[Parameter-dependent local center manifold]
  \label{thm:cmabstract}
  Let $T_0$ be an eventually compact $\mathcal{C}_0$-semigroup on a sun-reflexive real Banach space $X$ and let $T$ be the $\mathcal{C}_0$-semigroup on $X$ defined by \cref{eq:T_T0_AIE} where $L$ is a compact perturbation. Suppose that the generator $A$ of $T$ has $1 \le n_0 < \infty$ purely imaginary eigenvalues with corresponding $n_0$-dimensional real center eigenspace $X_0$. Furthermore, assume that $R$ in \cref{eq:aie_bold_1st} is $C^k$-smooth and \cref{eq:R_conditions} holds.

  Then there exists a $C^k$-smooth map $\mathcal{C} : U \times U_p \to X$ defined in a neighborhood of the origin in $X_0 \times \RR^p$ and such that for every sufficiently small $\alpha \in \RR^p$ the manifold $\CM(\alpha) \DEF \mathcal{C}(U,\alpha)$ is locally positively invariant for the semiflow generated by \cref{eq:aie_bold_1st} at parameter value $\alpha$. Furthermore, any solution $u : I \to X$ of \cref{eq:aie_bold_1st} that lies on $\CM(\alpha)$ satisfies \cref{eq:ode_on_cm_abstract}.  
\end{theorem}

\subsection{The special case of parameter-dependent classical DDEs}\label{sec:pd:ddes}
In this section we will formulate a corollary of \cref{thm:cmabstract} that applies specifically to the classical parameter-dependent DDE \cref{eq:pd-DDE}. As in \cref{sec:ddecase} our starting point is \cref{eq:pd-DDE12} with $F = 0$,
\begin{equation}
  \label{eq:trivialdde_bold}
  \left\{
    \begin{aligned}
      \dot{x}(t) &= 0,\\
      \dot{\mu}(t) &= 0,
    \end{aligned}
  \right.
  \qquad t \ge 0,
\end{equation}
in the unknown $(x, \mu)$ with initial condition $(\phi,\alpha)$ in the state space $\BM{X} \DEF X \times \RR^p$ where $X \DEF C([-h,0],\RR^n)$. So, we interpret the first component of \cref{eq:trivialdde_bold} as a DDE but the second component as an ODE. By comparison with \cref{eq:trivialdde} it is clear that the solution of the initial value problem for \cref{eq:trivialdde_bold} defines a $\mathcal{C}_0$-semigroup $\BM{T}_0$ on $\BM{X}$,
\[
 \BM{T}_0(t) \DEF \diag{(T_0(t), I_p)},
\]
with $T_0$ the eventually compact shift semigroup on $X$ from \cref{eq:LDDE} and $I_p$ the identity on $\RR^p$. We specify the perturbations $L$ and $L_p$ in \cref{eq:B_bold} as
\[
  L\phi = (D_1F(0,0)\phi)\rss, \qquad L_p\alpha = (D_2F(0,0)\alpha)\rss,
\]
where $D_1F(0,0)$ and $D_2F(0,0)$ are as in \cref{eq:pd:splitting}. Then $L$ is of finite rank, so certainly it is compact. Also, we choose the nonlinear perturbation $R$ in \cref{eq:R_bold} as
\begin{equation}
  \label{eq:pd:Rdde}
  R(\phi,\alpha) = G(\phi,\alpha)\rss,
\end{equation}
where $G$ is defined by the splitting in \cref{eq:pd:splitting}. Then \cref{eq:G_conditions} implies that the conditions in \cref{eq:R_conditions} hold.

\begin{corollary}[Parameter-dependent local center manifold for DDEs]
\label{cor:cmdde}
Consider the particular case of the classical DDE in \cref{eq:pd-DDE},
\[
  \dot{x}(t)= F(x_t, \alpha), \qquad t \ge 0,
\]
where $F: X \times \RR^p \to \RR^n$ is $C^k$-smooth for some $k \ge 1$ with $F(0,0) = 0$. Let $T$ be the $\mathcal{C}_0$-semigroup on $X$ corresponding to the linearization of \cref{eq:pd-DDE} at $0 \in X$ for the critical parameter value $\alpha = 0$. Suppose that the generator $A$ of $T$ has $1 \le n_0 < \infty$ purely imaginary eigenvalues with corresponding $n_0$-dimensional real center eigenspace $X_0$.

Then there exists a $C^k$-smooth map $\mathcal{C} : U \times U_p \to X$ defined in a neighborhood of the origin in $X_0 \times \RR^p$ and such that for every sufficiently small $\alpha \in \RR^p$ the manifold $\CM(\alpha) \DEF \mathcal{C}(U,\alpha)$ is locally positively invariant for the semiflow generated by \cref{eq:pd-DDE} at parameter value $\alpha$.

Furthermore, if the history $x_t$ associated with a solution of \cref{eq:pd-DDE} exists on some nondegenerate interval $I$ and $x_t \in \CM(\alpha)$ for all $t \in I$, then $u : I \to X$ defined by $u(t) \DEF x_t$ is differentiable and satisfies
\[
  j\dot{u}(t) = \SUNSTAR{A}ju(t) + (D_2F(0,0)\alpha)\rss + G(u(t),\alpha)\rss, \qquad \forall\,t \in I,
\]
where $\SUNSTAR{A}$ is the {\WSTAR} generator of $\SUNSTAR{T}$ and the operator $G : X \times \RR^p \to \RR^n$ defined by \cref{eq:pd:splitting},
\[
  G(\phi,\alpha) \DEF F(\phi, \alpha) - D_1F(0,0)\phi - D_2F(0,0)\alpha,
\]
is the nonlinear part of $F$.
\end{corollary}

\begin{remark}\label{rem:nonsmooth}
  For discrete classical DDEs \cref{Eq:FiniteDDE} one may want to use one or more of the discrete delays $\tau_1,\ldots,\tau_m$ as parameters. However, in this case the operator $R$ is typically no longer differentiable, see \textup{\cite[Remark IX.9.2]{diekmann1995delay}}. Of course, if only a single discrete delay is present, this problem can be avoided by a linear rescaling of time.
\end{remark}

\section{Normal forms on the parameter-dependent center manifold\label{sec:normal-forms}}
The normalization technique described in this section goes back to \cite{Coullet1983competinginstabilities}. In \cite{Kuznetsov1999} it was applied to obtain expressions for the critical normal form coefficients of all generic codimension one and two bifurcations of equilibria in ODEs, also see \cite[\S 8.7]{Kuznetsov2004}. In this context, these expressions are independent of the (finite) dimension of the phase space and they involve only critical eigenvectors of the Jacobian matrix and its transpose as well as higher order derivatives of the right-hand side at the critical equilibrium. These properties make them suitable for both symbolic and numerical evaluation.
\par
In \cite{Kuznetsov2008} the same technique was applied to parameter-dependent normal forms to start the continuation of nonhyperbolic cycles emanating from generalized Hopf, fold-Hopf and Hopf-Hopf bifurcation points of ODEs. The resulting predictors were implemented in the freely available software package \MATCONT, a \MATLAB toolbox for continuation and bifurcation analysis of finite dimensional dynamical systems. This makes it possible to verify transversality conditions and to initialize the continuation of the nonhyperbolic cycles mentioned above. A similar switching problem for iterated maps was solved earlier in \cite{Govaerts2007maps}.
\par
In \cite{Janssens:Thesis} the normalization technique was lifted to an infinite dimensional setting, providing explicit expressions for the critical normal form coefficients of generic codimension one and two equilibrium bifurcations in DDEs. These expressions were partially implemented \cite{Wage:Thesis:2014} in the software \DDEBIFTOOL. This package can be considered as the DDE equivalent of \MATCONT in command line mode.
\par
In this section we extend the normalization method from \cite{Janssens:Thesis} to include parameters. Suppose $0 \in X$ is a stationary state of \cref{eq:pd-DDE} at the critical parameter value $0 \in \RR^p$ and assume there are $n_0 \ge 1$ eigenvalues on the imaginary axis, counting algebraic multiplicities. Let $P_0$ be the corresponding \emph{real} spectral projector on $X$, so the range $X_0$ of $P_0$ is the \emph{real} $n_0$-dimensional center eigenspace. \cref{cor:cmdde} applies to give a parameter-dependent local center manifold $\CM(\alpha)$ for \cref{eq:pd-DDE}.
\par
We allow for the introduction of a new parameter $\beta$ defined in a neighborhood of $0 \in \RR^p$ such that $\alpha = K(\beta)$ for some locally defined $C^k$-diffeomorphism $K : \RR^p \to \RR^p$ that is to be determined below, up to a certain order. If $u : I \to X$ with $u(t) \DEF x_t \in \CM(\alpha)$ is as in \cref{cor:cmdde}, then $u$ is differentiable on $I$ and satisfies
\begin{equation}
  \label{eq:ODEonCM}
  j\dot{u}(t) = \SUNSTAR{A}ju(t) + (D_2F(0,0)K(\beta))\rss + R(u(t),K(\beta)), \qquad \forall\,t \in I,
\end{equation}
where $R$ encodes the nonlinear part of $F$ as in \cref{eq:pd:Rdde}. Choose a basis $\Phi$ of $X_0$ and let $\mathcal{H} : \RR^{n_0} \times \RR^p \to X$ be a locally defined $C^k$-smooth parametrization of $\CM(\alpha)$ with respect to $\Phi$ and in terms of the new parameter $\beta$, see the remark following \cref{def:cm_alpha}. For every $t \in I$ we define $v(t) \in \RR^{n_0}$ as the coordinate vector of $P_0 u(t)$ with respect to $\Phi$. Then $v : I \to \RR^{n_0}$ satisfies a parameter-dependent ordinary differential equation of the form
\begin{equation}
  \label{eq:ODEexpansion}
  \dot{v} = \sum_{|\nu|+|\mu| \geq 1}\frac{1}{\nu!\mu!}g_{\nu\mu}v^{\nu}\beta^{\mu},
\end{equation}
where the $C^k$-smooth vector field on the right has been expanded up to some sufficiently large - but finite - order. The multi-indices $\nu$ and $\mu$ have lengths $n_0$ and $p$, respectively. \emph{We assume that \cref{eq:ODEexpansion} is in parameter-dependent normal form}, up to a certain order. Since $\mathcal{H}$ parametrizes $\CM(\alpha)$,
\[
  u(t) = \mathcal{H}(v(t), \beta), \qquad t \in I,
\]
with both $u$ and $v$ depending on the parameter, although this is left implicit in the notation. Substituting the above relation into \cref{eq:ODEonCM} produces the \emph{homological equation}
\begin{equation}
  \label{eq:homological_equation}
  \tag{HOM}
  \SUNSTAR{A}j\mathcal{H}(v,\beta)+(D_2F(0,0)K(\beta))\rss + R(\mathcal{H}(v,\beta),K(\beta)) = jD_1\mathcal{H}(v,\beta)\dot{v},
\end{equation}
with $\dot{v}$ given by the parameter-dependent normal form \cref{eq:ODEexpansion}. The unknowns in \cref{eq:homological_equation} are $\mathcal{H}$, $K$ and the coefficients $g_{\nu\mu}$ from \cref{eq:ODEexpansion}. They are determined, up to a certain order, by the assumption that \cref{eq:ODEexpansion} is in normal form. For $r, s \ge 0$ with $r + s \ge 1$ we denote by $D_1^rD_2^sF(0,0) : X^r \times [\RR^p]^s \rightarrow \RR^n$ the mixed Fr\'echet derivative of order $r + s$, evaluated at $(0,0) \in X \times \RR^p$, with the understanding that at most one of the factor spaces $X^r$ or $[\RR^p]^s$ is absent if either $r = 0$ or $s = 0$. We expand the nonlinearity $R$ as
\begin{equation}
  \label{eq:R}
  R(\phi,\alpha)=\sum_{r + s > 1}\frac{1}{r!s!}D_1^rD_2^sF(0,0)(\phi^{(r)},\alpha^{(s)})\rss,
\end{equation}
where $\phi^{(r)} \DEF (\phi,\dots,\phi)\in X^r$ and $\alpha^{(s)} \DEF (\alpha,\dots,\alpha)\in [\RR^p]^s$. The mappings $\mathcal{H}$ and $K$ can be expanded as
\begin{equation}
  \label{eq:HKexpansion}
\mathcal{H}(v,\beta)\ =\sum_{|\nu|+|\mu| \geq 1}\dfrac{1}{\nu!\mu!}H_{\nu\mu}v^{\nu}\beta^{\mu}, \qquad K(\beta)=\sum_{|\mu| \geq 1}\dfrac{1}{\mu!}K_{\mu}\beta^{\mu}.
\end{equation}
Substituting \cref{eq:ODEexpansion,eq:R,eq:HKexpansion} into \cref{eq:homological_equation}, collecting coefficients of terms $v^{\nu}\beta^{\mu}$ from lower to higher order and solving the resulting linear systems, one can solve recursively for the unknown coefficients $g_{\nu\mu}$, $H_{\nu\mu}$ and $K_{\mu}$ by applying the Fredholm alternative and taking ordinary or bordered inverses, as explained in \cref{sec:solvability}.

\section{Predictors for normal forms on center manifolds\label{sec:predictors}}
Starting from this point, we will focus exclusively on \emph{two-parameter} DDEs, i.e. we will have $p=2$. This enables the initialization of codimension one equilibrium and cycle bifurcations emanating from the generalized Hopf, fold-Hopf, Hopf-Hopf and transcritical-Hopf codimension two bifurcation points. The corresponding steady state and cycle bifurcation curves can then be continued in two parameters using the continuation capabilities of \DDEBIFTOOL, providing an unfolding of the singularity.

In this section we list known asymptotics for codimension one nonhyperbolic cycles emanating from generalized Hopf, fold-Hopf and Hopf-Hopf bifurcations obtained in \cite{Kuznetsov2008}. The asymptotics are derived from the Poincar\'e normal forms which are obtained by considering near-identity changes of coordinates generated by homogeneous polynomial functions without using time reparametrization. These normal forms are therefore ready to be used in conjunction with the homological equation \cref{eq:homological_equation} where all the time derivatives are assumed to be with the same unit of time. Following the same method as in \cite{Kuznetsov2008}, we also derive asymptotics for codimension one nonhyperbolic cycles emanating from the transcritical-Hopf bifurcation. We also provide asymptotics for the codimension one equilibrium bifurcations emanating from the degenerate Hopf points.

\subsection{Generalized-Hopf bifurcation\label{sec:GH_predictors}}
Suppose that system \cref{eq:pd-DDE} has an equilibrium $\phi=0$ at the critical parameter value $\alpha_{0}=(0,0)\in \RR^{2}$ with purely imaginary eigenvalues
\begin{equation}
\lambda_{1,2}=\pm i\omega_{0},\qquad\omega_{0}>0.\label{eq:GH_eigenvalues}
\end{equation}
Furthermore, suppose that the first Lyapunov coefficient $\ell_{1}(0)=0$ while the second Lyapunov coefficient $\ell_{2}(0)\neq0$. Then the restriction of \cref{eq:pd-DDE} to the two-dimensional center manifold $\CM(\alpha)$ can be transformed into the smooth local normal form
\begin{equation}
\dot{z}=\lambda(\alpha)z+c_{1}(\alpha)z|z|^{2}+c_{2}(\alpha)z|z|^{4}+\mathcal{O}(|z|^{6}),\label{eq:GH_nf-1}
\end{equation}
where $\lambda(\alpha),c_{1}(\alpha),c_{2}(\alpha)$ are complex-valued functions with $\ell_{1}(0)=\dfrac{1}{\omega_0}\Re\,c_{1}(0)=0$, $\lambda(0)=i\omega_{0}$ and $\ell_{2}(0)=\dfrac{1}{\omega_0}\Re\,c_{2}(0)\neq0$. Let
\[
\left\{ \begin{array}{rl}
\begin{aligned}
\lambda(\alpha) & =\mu(\alpha)+i\omega(\alpha),\\
c_{1}(\alpha) & =\Re c_{1}(\alpha)+i\Im c_{1}(\alpha),
\end{aligned}
\end{array}\right.
\]
where $\mu(\alpha)$ and $\omega(\alpha)$ are real-valued functions. Suppose that the map $\alpha\mapsto(\mu(\alpha),\Re c_{1}(\alpha))$ is regular at $\alpha=0$. Then we can introduce new parameters $(\beta_{1}(\alpha),\beta_{2}(\alpha))$ to obtain the normal form

\begin{equation}
\dot{z}=(\beta_{1}+i\omega(\beta))z+(\beta_{2}+i \Im c_{1}(\beta))z|z|^{2}+c_{2}(\beta)z|z|^{4}+\mathcal{O}(|z|^{6}),\label{eq:GH_nf-1-2}
\end{equation}
where $\omega(0)=\omega_{0}$. For convenience, we abuse notations and write $\omega(\beta)$ and $c_j(\beta)$ instead of $\omega(\alpha(\beta))$ and $c_j(\alpha(\beta))$, respectively. Similar conventions are adopted in other cases ahead.

It is well known that a curve of fold bifurcation of limit cycles (LPC) emanates from this codimension two point. To approximate this curve we substitute $z=\rho e^{i\psi}$ into \cref{eq:GH_nf-1-2} and expand $\omega(\beta)=\omega(\beta)=\omega_0+\omega_{10}\beta_1+\omega_{01}\beta_2 + \mathcal O(\|(\beta_1,\beta_2)\|^2)$.
%
%
Truncating higher order terms and separating the real and imaginary parts yields the system
\begin{equation}
\begin{cases}
\begin{aligned}
\dot{\rho} & =\rho\left(\beta_{1}+\beta_{2}\rho^{2}+\Re(c_{2}(0))\rho^{4}\right),\\
\dot{\psi} & =\omega_0+\omega_{10}\beta_1+\omega_{01}\beta_2 +\Im(c_{1}(0))\rho^{2}+\Im(c_{2}(0))\rho^{4}.
\end{aligned}
\end{cases}\label{eq:GH_amplitude_system}
\end{equation}
The curve LPC corresponds to a double zero in the amplitude equation in \cref{eq:GH_amplitude_system}. Therefore, this curve in \cref{eq:GH_nf-1-2} can be approximated by
\begin{equation}
\label{eq:GH_approximation_LPC}
\rho=\epsilon,\qquad\beta_{1}=\Re(c_{2}(0)) \epsilon^{4},\qquad\beta_{2}=-2 \Re(c_{2}(0)) \epsilon^{2}, \qquad  \epsilon>0.
\end{equation}
From the second equation in the amplitude system \cref{eq:GH_amplitude_system} we obtain using equation \cref{eq:GH_approximation_LPC} the following approximation for the period:
\begin{equation}
\label{eq:GH_approximation_LPC_T}
T =2\pi\left/\left(\omega_0+\left(\Im(c_{1}(0))-2\Re c_2(0) \omega_{01} \right)\epsilon^{2} \right) \right. .
\end{equation}
Since $\Re(c_{1}(\beta))=\beta_{2}$, it is easy to see that the Hopf curve in the original system is related to the truncated normal form by
\begin{equation}
\label{eq:GH_approximation_Hopf}
(\beta_{1},\beta_{2},z)=(0,\epsilon,0)
\end{equation}
for $\epsilon\neq0$ small.

\subsection{Fold-Hopf bifurcation\label{sec:Predictor-FH}}
Suppose that system \cref{eq:pd-DDE} has an equilibrium $x=0$ at the critical parameter value $\alpha_{0}=(0,0)\in \RR^{2}$ with eigenvalues
\begin{equation}
\lambda_{1}=0,\qquad\lambda_{2,3}=\pm i\omega_{0},\label{eq:FH_eigenvalues}
\end{equation}
where $\omega_{0}>0$. The restriction of \cref{eq:pd-DDE} to the three-dimensional center manifold $\CM(\alpha)$ can generically be transformed to the smooth local normal form
\begin{equation}
\begin{cases}
\begin{aligned}
\dot{z}_{0} & =\gamma(\alpha)+g_{200}(\alpha)z_{0}^{2}+g_{011}(\alpha)|z_{1}|^{2}+g_{300}(\alpha)z_{0}^{3}+g_{111}(\alpha)z_{0}|z_{1}|^{2}\\
 & \qquad+\mathcal{O}\left(\|\left(z_{0},z_{1},\overline{z}_{1}\right)\|^{4}\right),\\
\dot{z}_{1} & =\lambda(\alpha)z_{1}+g_{110}(\alpha)z_{0}z_{1}+g_{210}(\alpha)z_0^2z_{1}+g_{021}(\alpha)z_{1}|z_{1}|^{2}
	+ \mathcal{O}\left(\|\left(z_{0},z_{1},\overline{z}_{1}\right)\|^{4}\right),
\end{aligned}
\end{cases}\label{eq:fold-Hopf_poincare_normal_form}
\end{equation}
where $z_{0}\in \RR$, $z_{1}\in \CC$, $\gamma(0)=0,\lambda(0)=i\omega_{0}$ and the functions $g_{jkl}(\alpha)$ are real in the first equation and complex in the second. Let $\lambda(\alpha) =\mu(\alpha)+i\omega(\alpha)$ and suppose that the map $\alpha\mapsto(\gamma(\alpha),\mu(\alpha))$ is regular at $\alpha=0$. Introducing new parameters $(\beta_{1}(\alpha),\beta_{2}(\alpha))$, we obtain the truncated normal form
\begin{equation}
\begin{cases}
\begin{aligned}
\dot{z}_{0} & =\beta_{1}+g_{200}(\beta)z_{0}^{2}+g_{011}(\beta)|z_{1}|^{2}+g_{111}(\beta)z_{0}|z_{1}|^{2}+g_{300}(\beta)z_{0}^{3},\\
\dot{z}_{1} & =(\beta_{2}+i\omega_{0}+ib_{1}(\beta))z_{1}+g_{110}(\beta)z_{0}z_{1}+g_{210}(\beta)z_{0}^{2}z_{1}+g_{021}(\beta)z_{1}|z_{1}|^{2},
\end{aligned}
\end{cases}
\label{eq:FH-nf}
\end{equation}
with $b_{1}(0)=0$. Letting $z_{1}=\rho e^{i\psi}$ and separating the real and imaginary parts yields the system
\begin{equation}
\begin{cases}
\begin{aligned}
\dot{z}_{0} & =\beta_{1}+g_{200}(\beta)z_{0}^{2}+g_{011}(\beta)\rho^{2}+g_{111}(\beta)z_{0}\rho^{2}+g_{300}(\beta)z_{0}^{3},\\
\dot{\rho} & =\rho\left(\beta_{2}+\Re(g_{110}(\beta))z_{0}+\Re(g_{210}(\beta))z_{0}^{2}+\Re(g_{021}(\beta))\rho^{2}\right),\\
\dot{\psi} & =\omega_{0}+b_{1}(\beta)+\Im(g_{110}(\beta))z_{0}+\Im(g_{210}(\beta))z_{0}^{2}+\Im(g_{021}(\beta))\rho^{2}.
\end{aligned}
\end{cases}\label{eq:FH_nf_phi_psi}
\end{equation}
\subsubsection{Neimark-Sacker curve}
\label{Sec:ZH_NS_predictors}
A Hopf bifurcation in the amplitude system \cref{eq:FH_nf_phi_psi}, i.e. when the trace of the Jacobian matrix of the amplitude system vanishes but the determinant of this matrix is positive,  corresponds to a Neimark-Sacker bifurcation in the original system. Let $\rho=\epsilon$ and suppose that $g_{011}(0)\Re(g_{110}(0))<0$. Then by the implicit function theorem, we obtain the second-order predictor
\begin{equation}
\begin{cases}
\begin{aligned}
\beta_{1} & =-g_{011}(0)\epsilon^{2},\\
\beta_{2} & =\dfrac{\Re(g_{110}(0))\left(2\Re(g_{021}(0))+g_{111}(0)\right)-2\Re(g_{021}(0))g_{200}(0)}{2g_{200}(0)}\epsilon^{2}, \\
z_0&=-\dfrac{2\Re\left(g_{021}(0)\right)+g_{111}(0)}{2g_{200}(0)}\epsilon^{2},
\end{aligned}
\end{cases}\label{eq:FH_NS_predictor}
\end{equation}
for the Neimark-Sacker curve.
It follows from \cref{eq:FH_nf_phi_psi} that the period of the cycle is approximated to the same order by
\[
T =2\pi\left/\left(\omega_0+\omega_{1}\beta_{1}+\omega_{2}\beta_{2}+\Im(g_{110}(0))z_{0}+\Im(g_{021}(0))\epsilon^{2}\right)\right..
\]
Here $(z_{0},\beta_{1},\beta_{2})$ are as in \cref{eq:FH_NS_predictor} and $\omega_{i}=\dfrac{\partial b_{1}}{\partial\beta_{i}}(0)$, for $i=1,2$.

\begin{remark}
Notice that the cubic terms $g_{210}(0)$ and $g_{300}(0)$ do not appear in the predictor. However, these coefficients are needed to determine the stability of the two-dimensional torus, see \textup{\cite[Section 8.5]{Kuznetsov2004}} and \cref{ex:RH:e}.
\end{remark}

\subsubsection{Fold and Hopf curves}
\label{Sec:ZH_FoldHopf_predictors}
For the approximation to the fold and Hopf curves it is sufficient to consider the second order terms in $(z_{0},\rho)$ in the amplitude system in \cref{eq:FH_nf_phi_psi}. We obtain three equilibrium solutions given by
\[
\pm\left(\sqrt{-\frac{\beta_{1}}{g_{200}(0)}},0\right) \qquad \text{and} \qquad \left(-\frac{\beta_{2}}{\Re(g_{110}(0))},\frac{\sqrt{-\Re(g_{110}(0))^{2}\beta_{1}-g_{200}(0)\beta_{2}^{2}}}{\sqrt{g_{011}(0)}\Re(g_{110}(0))}\right).
\]
It follows that the fold curve is approximated by
\[
\left(\beta_{1},\beta_{2}\right)=\left(0,\epsilon\right)
\]
and the Hopf curve by
\[
\left(\beta_{1},\beta_{2}\right)=\left(-\frac{g_{200}(0)}{\Re(g_{110}(0))^{2}}\epsilon^{2},\epsilon\right),
\]
for $|\epsilon|$ small.

\subsection{Hopf-Hopf bifurcation\label{sec:HH_pedictors}}
Suppose that system \cref{eq:pd-DDE} at the critical parameter value $\alpha_{0}=(0,0)\in \RR^{2}$ undergoes two Hopf bifurcations simultaneously. Then the generator $A$ has two pairs of purely imaginary eigenvalues
\begin{equation}
\lambda_{1,4}=\pm i\omega_{1}(0),\qquad\lambda_{2,3}=\pm i\omega_{2}(0),\label{eq:HH_eigenvalues}
\end{equation}
where $\omega_{i}: \RR^{2}\rightarrow \CC$ for $i=1,2$ such that $\omega_{1}(0)>\omega_{2}(0)>0$. When no other eigenvalues on the imaginary axis exist, this phenomenon is called the Hopf-Hopf or double-Hopf bifurcation. Assume, furthermore that the nonresonance conditions $k\omega_{1}\neq l\omega_{2},$ with $0<k+l\text{\ensuremath{\leq}}5$ are satisfied. Then the restriction of \cref{eq:pd-DDE} to the four-dimensional center manifold $\CM(\alpha)$ can be transformed to the smooth local normal form
\begin{equation*}
\begin{cases}
\begin{aligned}
\dot{z}_{1} & =\lambda_{1}(\alpha)z_{1}+g_{2100}(\alpha)z_{1}|z_{1}|^{2}+g_{1011}(\alpha)z_{1}|z_{2}|^{2}+g_{3200}(\alpha)z_{1}|z_{1}|^{4}\\
 & \qquad+g_{2111}(\alpha)z_{1}|z_{1}|^{2}|z_{2}|^{2}+g_{1022}(\alpha)z_{1}|z_{2}|^{4}+\mathcal{O}\left(\|z_{1},\overline{z_{1}},z_{2},\overline{z_{2}}\|^{6}\right),\\
\dot{z}_{2} & =\lambda_{2}(\alpha)z_{2}+g_{1110}(\alpha)z_{2}|z_{1}|^{2}+g_{0021}(\alpha)z_{2}|z_{2}|^{2}+g_{2210}(\alpha)z_{2}|z_{1}|^{4}\\
 & \qquad+g_{1121}(\alpha)z_{2}|z_{1}|^{2}|z_{2}|^{2}+g_{0032}(\alpha)z_{2}|z_{2}|^{4}+\mathcal{O}\left(\|z_{1},\overline{z_{1}},z_{2},\overline{z_{2}}\|^{6}\right),
\end{aligned}
\end{cases}\label{eq:HH_nf}
\end{equation*}
where $z_{1,},z_{2}\in \CC^{2}$ and $g_{jklm}\in \CC.$
Let
\[
\begin{cases}
\begin{aligned}
\lambda_{1}(\alpha) & =\mu_{1}(\alpha)+i\omega_{1}(\alpha),\\
\lambda_{2}(\alpha) & =\mu_{2}(\alpha)+i\omega_{2}(\alpha),
\end{aligned}
\end{cases}
\]
and suppose that the map $\alpha\mapsto(\mu_{1}(\alpha),\mu_{2}(\alpha))$ is regular at $\alpha=0$. Then we can introduce new parameters $(\beta_{1}(\alpha),\beta_{2}(\alpha))=(\mu_{1}(\alpha),\mu_{2}(\alpha))$ to obtain the normal form
\[
\begin{cases}
\begin{aligned}
\dot{z}_{1} & =(\beta_{1}+i\omega_{1}(\beta)z_{1}+g_{2100}(\beta)z_{1}|z_{1}|^{2}+g_{1011}(\beta)z_{1}|z_{2}|^{2}+g_{3200}(\beta)z_{1}|z_{1}|^{4}\\
 & \qquad+g_{2111}(\beta)z_{1}|z_{1}|^{2}|z_{2}|^{2}+g_{1022}(\beta)z_{1}|z_{2}|^{4}+\mathcal{O}\left(\|z_{1},\overline{z_{1}},z_{2},\overline{z_{2}}\|^{6}\right),\\
\dot{z}_{2} & =(\beta_{2}+i\omega_{2}(\beta)z_{2}+g_{1110}(\beta)z_{2}|z_{1}|^{2}+g_{0021}(\beta)z_{2}|z_{2}|^{2}+g_{2210}(\beta)z_{2}|z_{1}|^{4}\\
 & \qquad+g_{1121}(\beta)z_{2}|z_{1}|^{2}|z_{2}|^{2}+g_{0032}(\beta)z_{2}|z_{2}|^{4}+\mathcal{O}\left(\|z_{1},\overline{z_{1}},z_{2},\overline{z_{2}}\|^{6}\right).
\end{aligned}
\end{cases}
\]
%
We truncate the normal form to third order
\begin{equation}
\begin{cases}
\begin{aligned}
\dot{z}_{1} & =\left(\beta_{1}+i\omega_{1}+ib_{1}(\beta)\right)z_{1}+g_{2100}(0)z_{1}|z_{1}|^{2}+g_{1011}(0)z_{1}|z_{2}|^{2},\\
\dot{z}_{2} & =\left(\beta_{2}+i\omega_{2}+ib_{2}(\beta)\right)z_{2}+g_{1110}(0)z_{2}|z_{1}|^{2}+g_{0021}(0)z_{2}|z_{2}|^{2}.
\end{aligned}
\end{cases}\label{eq:HH_nf-1}
\end{equation}
Letting $\left(z_{1},z_{2}\right)=\left(\rho_{1}e^{i\psi_{1}},\rho_{2}e^{i\psi_{2}}\right)$
and separating the real and imaginary parts yields
\begin{equation}
\begin{cases}
\begin{aligned}
\dot{\rho_{1}}& = \rho_{1}\left(\beta_{1}+\Re(g_{2100}(0))\rho_{1}^{2}+\Re(g_{1011}(0))\rho_{2}^{2}\right),\\
\dot{\rho_{2}}& = \rho_{2}\left(\beta_{2}+\Re(g_{1110}(0))\rho_{1}^{2}+\Re(g_{0021}(0))\rho_{2}^{2}\right),\\
\dot{\psi}_{1}& = \omega_{1}+b_{1}(\beta)+\Im(g_{2100}(0))\rho_{1}^{2}+\Im(g_{1011}(0))\rho_{2}^{2},\\
\dot{\psi}_{2}& = \omega_{2}+b_{2}(\beta)+\Im(g_{1110}(0))\rho_{1}^{2}+\Im(g_{0021}(0))\rho_{2}^{2}.
\end{aligned}
\end{cases}\label{eq:HH-polor_coordinates}
\end{equation}

\subsubsection{Hopf curves}
\label{Sec:FH_Hopfs}

There are two semi-trivial equilibria
\[
(\rho_{1},\rho_{2})=\left(\sqrt{-\frac{\beta_{1}}{\Re(g_{2100}(0))}},0\right),\qquad(\rho_{1},\rho_{2})=\left(0,\sqrt{-\frac{\beta_{2}}{\Re(g_{0021}(0))}}\right)
\]
of the amplitude system of \cref{eq:HH-polor_coordinates}. Translating to the original system provides the Hopf bifurcation curves
\[
H_{1}=\left\{ \left(\beta_{1},\beta_{2}\right):\beta_{1}=0\right\} ,\qquad\text{and}\qquad H_{2}=\left\{ \left(\beta_{1},\beta_{2}\right):\beta_{2}=0\right\} .
\]

\subsubsection{Neimark-Sacker curves}
\label{Sec:FH_NSs}

A nontrivial equilibrium to the amplitude system
\begin{multline*}
\left(\rho_{1},\rho_{2}\right)
=\left(\sqrt{\frac{\beta_{2}\Re g_{1011}(0)-\beta_{1}\Re g_{0021}(0)}{\Re g_{0021}(0)\Re g_{2100}(0)-\Re g_{1011}(0)\Re g_{1110}(0)}},\right. \\
\left.\sqrt{\frac{\beta_{2}\Re g_{2100}(0)-\beta_{1}\Re g_{1110}(0)}{\Re g_{1011}(0)\Re g_{1110}(0)-\Re g_{0021}(0)\Re g_{2100}(0)}}\right)
\end{multline*}
corresponds to a torus of the original system. When
\begin{equation*}
\Re(g_{1110}(0))\beta_{1}=\Re(g_{2100}(0))\beta_{2}
\end{equation*}
the nontrivial equilibrium coincides with the first semi-trivial equilibrium, thus giving a predictor for a Neimark-Sacker bifurcation curve. Similarly, when
\begin{equation*}
\Re(g_{0021}(0))\beta_{1}=\Re(g_{1011}(0))\beta_{2}
\end{equation*}
the nontrivial equilibrium coincides with the second semi-trivial equilibrium, and gives a predictor for the second Neimark-Sacker bifurcation curve. Therefore, we obtain two Neimark-Sacker bifurcation curves in \cref{eq:HH_nf-1}, with approximations given by
\begin{subequations}
\begin{align}
(\rho_{1},\rho_{2},\beta_{1},\beta_{2}) &= \left(\epsilon,0,-\Re(g_{2100}(0))\epsilon^{2},-\Re(g_{1110}(0))\epsilon^{2}\right),\label{eq:HH_NS1_pm}\\
(\rho_{1},\rho_{2},\beta_{1},\beta_{2}) &= \left(0,\epsilon,-\Re(g_{1011}(0))\epsilon^{2},-\Re(g_{0021}(0))\epsilon^{2}\right),\label{eq:HH_NS2_pm}
\end{align}\label{eq:HH_NS_asymptotics}
\end{subequations}
where $\epsilon>0$, which are the predictors given in \cite{Kuznetsov2008} and \cite{Meijer2009}. An approximation for the period of the cycle for each Neimark-Sacker predictor can be obtained from the third and fourth equation in system \cref{eq:HH-polor_coordinates}, yielding
\begin{equation*}
\begin{cases}
\begin{aligned}
T_{1} & =2\pi \left/\left(\omega_{1}+b_{11}\beta_{1}+b_{12}\beta_{2}+\Im(g_{2100}(0))\epsilon^{2}\right)\right.,\\
T_{2} & =2\pi\left/\left(\omega_{2}+b_{21}\beta_{1}+b_{22}\beta_{2}+\Im(g_{0021}(0))\epsilon^{2}\right)\right. ,
\end{aligned}
\end{cases}
\label{eq:HH_period_predictors}
\end{equation*}
where
\begin{equation*}
\label{eq:HH_bij}
b_{jk}=\dfrac{\partial b_j}{\partial \beta_{k}}(0),~~j,k=1,2.
\end{equation*}
Here we should use $(\beta_{1},\beta_{2})$ as in \cref{eq:HH_NS1_pm,eq:HH_NS2_pm}.

\subsection{Transcritical-Hopf bifurcation\label{sec:HT_predictors}}
A majority of papers in which fold-Hopf bifurcations in DDEs are studied, deals with models where the steady state remains fixed under variation of parameters. In this case the unfolding is not given by \cref{eq:fold-Hopf_poincare_normal_form} anymore and we have to consider the smooth local normal form
\begin{equation*}
\begin{cases}
\begin{aligned}
\dot{z}_{0} & =\gamma(\alpha)z_{0}+g_{200}(\alpha)z_{0}^{2}+g_{011}(\alpha)|z_{1}|^{2}+g_{300}(\alpha)z_{0}^{3}+g_{111}(\alpha)z_{0}|z_{1}|^{2}\\
 & \qquad+\mathcal{O}\left(\|\left(z_{0},z_{1},\overline{z}_{1}\right)\|^{4}\right),\\
\dot{z}_{1} & =\lambda(\alpha)z_{1}+g_{110}(\alpha)z_{0}z_{1}+g_{210}(\alpha)w^{2}z_{1}+g_{021}(\alpha)z_{1}|z_{1}|^{2}+\mathcal{O}\left(\|\left(z_{0},z_{1},\overline{z}_{1}\right)\|^{4}\right).
\end{aligned}
\end{cases}\label{eq:fold-Hopf_poincare_normal_form-2}
\end{equation*}
The bifurcation analysis can be carried out similar to the fold-Hopf case, see \cite{Guo2008} and \cite{Wang2010Hopftranscritical}. An alternative approach is presented in \cite{Wu2012}. In contrast with the fold-Hopf bifurcation,
there are in general two Neimark-Sacker bifurcation curves. Furthermore, the fold bifurcation curve becomes a transcritical bifurcation curve, and meets the Hopf bifurcation curve transversally.

Under the assumption of the same transversality condition as in the fold-Hopf bifurcation we introduce new parameters $(\beta_{1}(\alpha),\beta_{2}(\alpha))$ to obtain the truncated normal form
\begin{equation}
\begin{cases}
\begin{aligned}
\dot{z}_{0} & =\beta_{1}z_{0}+g_{200}(\beta)z_{0}^{2}+g_{011}(\beta)|z_{1}|^{2}+g_{111}(\beta)z_{0}|z_{1}|^{2}+g_{300}(\beta)z_{0}^{3},\\
\dot{z}_{1} & =(\beta_{2}+i\omega_{0}+ib_{1}(\beta))z_{1}+g_{110}(\beta)z_{0}z_{1}+g_{210}(\beta)z_{0}^{2}z+g_{021}(\beta)z_{1}|z_{1}|^{2},
\end{aligned}
\end{cases}\label{eq:Ht-nf}
\end{equation}
with $b_{1}(0)=0$. Letting $z_{1}=\rho e^{i\psi}$ and separating the real and imaginary parts yields the three dimensional system
\begin{equation}
\begin{cases}
\begin{aligned}
\dot{z}_{0} & =\beta_{1}z_{0}+g_{200}(\beta)z_{0}^{2}+g_{011}(\beta)\rho^{2}+g_{111}(\beta)z_{0}\rho^{2}+g_{300}(\beta)z_{0}^{3},\\
\dot{\rho} & =\rho\left(\beta_{2}+\Re(g_{110}(\beta))z_{0}+\Re(g_{210}(\beta))z_{0}^{2}+\Re(g_{021}(\beta))\rho^{2}\right),\\
\dot{\psi} & =\omega_{0}+b_{1}(\beta)+\Im(g_{110}(\beta))z_{0}+\Im(g_{210}(\beta))z_{0}^{2}+\Im(g_{021}(\beta))\rho^{2}.
\end{aligned}
\end{cases}\label{eq:FH_nf_phi_psi-2}
\end{equation}

\subsubsection{Neimark-Sacker bifurcation curves}
\label{Sec:HT_NS_predictors}

Following the same procedure as in \cref{sec:Predictor-FH}, we obtain that for $g_{011}(0)\Re(g_{110}(0))<0$ there are two Neimark-Sacker bifurcation curves approximated by
\begin{equation}
\begin{cases}
\begin{aligned}
\beta_{1} & =\mp 2\sqrt{g_{011}(0)g_{200}(0)}\epsilon,\\
\beta_{2} & =\mp\Re\left(g_{110}(0)\right)\sqrt{\frac{g_{011}(0)}{g_{200}(0)}}\epsilon,\\
z_0       & =\pm\sqrt{\frac{g_{011}(0)}{g_{200}(0)}}\epsilon,
\end{aligned}
\end{cases}
\label{eq:HT_NS_predictor}
\end{equation}
while the period of the corresponding cycle is approximated by
\begin{align*}
T = 2\pi \left/\left(\omega_{0}+\omega_{1}\beta_{1}+\omega_{2}\beta_{2}+\Im(g_{110}(0))z_{0}\right)\right. .
\end{align*}
Here $(z_{0},\beta_{1},\beta_{2})$ are as in \cref{eq:HT_NS_predictor} and $\omega_{i}=\dfrac{\partial b_{1}}{\partial\beta_{i}}(0)$, for $i=1,2$.

\subsubsection{Hopf and transcritical bifurcations curves}
\label{Sec:HT_HopfTrans_predictors}

The transcritical bifurcation curve in the normal form is obtained by substituting $\rho=0$ in the amplitude system of \cref{eq:FH_nf_phi_psi-2}. Then $\beta_{2}$ is unrestricted and  $z_{0}=-\beta_{1}/g_{200}(0)$. The transcritical bifurcation curve is therefore given by
\[
\left(\beta_{1},\beta_{2}\right)=\left(0,\beta_{2}\right).
\]
To obtain a predictor for the Hopf bifurcation curve we truncate \cref{eq:FH_nf_phi_psi-2} to the second order. We obtain a trivial equilibrium $(z_{0},\rho)=(0,0)$, a semi-trivial equilibrium $(z_{0},\rho)=(-\frac{\beta_{1}}{g_{200}(0)},0)$ and a nontrivial equilibrium
\[
(z_{0},\rho)=\left(-\frac{\beta_{2}}{\Re\left(g_{110}(0)\right)},\dfrac{\sqrt{\beta_{2}\left(\Re\left(g_{110}(0)\right)\beta_{1}-g_{200}(0)\beta_{2}\right)}}{\Re\left(g_{110}(0)\right)\sqrt{g_{011}(0)}}\right).
\]
It follows that the Hopf bifurcation curves are approximated by
\[
\beta_{2}=\dfrac{\Re\left(g_{110}(0)\right)}{g_{200}(0)}\beta_{1}, \qquad \beta_{2}=0.
\]

\section{Coefficients of parameter-dependent normal forms\label{sec:Coefficients-of-parameter}}
Using the method outlined in \cref{sec:normal-forms}, we derive here the coefficients needed for the predictors of the nonhyperbolic equilibria and cycles emanating from generalized Hopf, fold-Hopf, Hopf-Hopf and transcritical-Hopf bifurcations, see \cref{sec:predictors}. While doing so, we also obtain the critical normal form coefficients, which were first derived in \cite{Janssens:Thesis}. Recall, that we conside only two-parameter DDEs, so that $p=2$.

For the derivation of the coefficients in this section it is sufficient to expand and truncate the nonlinearity $R$ and the parameter-mapping $K$ in \cref{eq:R,eq:HKexpansion}, respectively, as follows
\begin{equation}
\label{eq:R_truncated}
\begin{array}{rl}
R(u,\alpha)\ =& {\displaystyle \left(\frac{1}{2} B(u,u) + A_1(u,\alpha)+\frac{1}{2} J_2(\alpha,\alpha)+\frac{1}{6} C(u,u,u) + \frac{1}{2} B_1(u,u,\alpha) \right.} \\
&\quad {\displaystyle \left. +  \frac{1}{24} D(u,u,u,u) + \frac{1}{6} C_1(u,u,u,\alpha) + \frac{1}{120} E(u,u,u,u,u) \right) \rss, }
\end{array}
\end{equation}
\begin{equation}
\label{eq:K_truncated}
\alpha=K(\beta)=K_{10} \beta_1  + K_{01} \beta_2. 
\end{equation}
Here $u \in X$, while $\alpha,\beta \in \RR^2$, and $B$, $A_1$, $J_2$, $C$, $B_1$, $D$, $C_1$ and $E$ are the standard multilinear forms arising from the expansion of $F(u,\alpha)$ (or $G(u,\alpha)$). For example,
\[
B(u,u)=D^2_1F(0,0)(u,u),~A_1(\alpha,\alpha)=D_2^2F(0,0)(\alpha,\alpha),~B_1(u,u,\alpha)=D^1_2D^2_1F(0,0)(u,u,\alpha),
\]
etc. These forms are $ \RR^n$-valued on real arguments, while they are linearly extended (`complexified') to $ \CC^n$-valued ones on complex-valued arguments. Finally, we introduce
\begin{equation}
\label{eq:J_1} 
J_1\  =\ D_2F(0,0). 
\end{equation}

We assume in all situations that $\phi_{0}=0$ is a steady state of \cref{eq:pd-DDE} at the critical parameter value $\alpha_{0}=0\in \RR^{2}$. Explicit formulas to compute the multilinear forms for the simplest DDE \cref{Eq:FiniteDDE} are given in \cref{sec:Implement}.

\subsection{Generalized Hopf bifurcation}\label{sec:GH_coef}
Since the eigenvalues \cref{eq:GH_eigenvalues} are simple, there exist eigenfunctions $\phi$ and $\SUN{\phi}$ such that
\[
A\phi=i\omega_{0}\phi, \qquad \STAR{A}\SUN{\phi}=i\omega_{0}\SUN{\phi}, \qquad \PAIR{\SUN{\phi}}{\phi} = 1.
\]
The eigenfunctions $\phi$ and $\SUN{\phi}$ are explicitly given by \cref{eq:eigenfunction,eq:eigenfunction1} with $q\in \CC^n$ and $p \in \CCC{n}$ satisfying
\[
\Delta(i\omega_{0})q=0, \qquad p\Delta(i\omega_{0})=0, \qquad p\Delta'(i\omega_0)q = 1.
\]
Any point $y\in X_{0}$ in the real critical eigenspace can be represented as
\[
y=z\phi+\bar{z}\bar{\phi},\qquad z\in \CC,
\]
where $z=\PAIR{\SUN\phi}{y} $. Therefore, the homological equation \cref{eq:homological_equation} can be written as
\[
\SUNSTAR{A} j\mathcal{H}(z,\bar{z},\beta)+J_{1}K(\beta)\rss+R(\mathcal{H}(z,\bar{z},\beta),K(\beta))=j(D_{z}\mathcal{H}(z,\bar{z},\beta)\dot{z}+D_{\bar{z}}\mathcal{H}(z,\bar{z},\beta)\dot{\bar{z}}),
\]
where $\beta=(\beta_1,\beta_2)$, $\dot{z}$ is given by the normal form \cref{eq:GH_nf-1-2} and $\mathcal{H}$ admits the expansion
\begin{equation}
\label{eq:GH-H}
\mathcal{H}(z,\bar{z},\beta) = z\phi+\bar{z}\bar{\phi}
 + H_{0010}\beta_1 + H_{0001} \beta_2
 +\sum_{j+k+|\mu| \geq 2}\dfrac{1}{j!k!\mu!}H_{jk\mu}z^{j}\bar{z}^{k}\beta^{\mu}.
\end{equation}
For the predictors derived in \cref{sec:GH_predictors} we need the truncated parameter-dependent normal form
\[
\dot{z}=\left(i\omega_0+\beta_{1}+i\omega_{10}\beta_{1}+i\omega_{01}\beta_{2}\right)z+\left(\beta_{2}+c_{1}(0)\right)z|z|^{2}+c_{2}(0)z|z|^{4},
\]
where $\Re(c_1(0))=0$.

\subsubsection{Critical normal form coefficients}
We start by calculating the critical normal form coefficients following \cite{Janssens:Thesis}. Collecting the coefficients of the quadratic terms $z^{2}$ and $z\bar{z}$ in the homological equation yields two nonsingular linear systems:
\begin{align*}
\LHS{2i\omega_{0}} jH_{2000}  &= B(\phi,\phi) \rss , \\
-\SUNSTAR{A}jH_{1100} &= B(\phi,\bar{\phi}) \rss .
\end{align*}
They are solved using \cref{lem:regular_solution} to give
\begin{align*}
H_{2000}(\theta) & =e^{2i\omega_{0}\theta}\Delta^{-1}(2i\omega_0)B(\phi,\phi),\\
H_{1100}(\theta) & =\Delta^{-1}(0)B(\phi,\bar{\phi}).
\end{align*}
For the cubic terms, the system corresponding to $z^3$ is also nonsingular,
\[
\LHS{3i\omega_{0}} jH_{3000} = \rsswp{3B(\phi,H_{2000})+C(\phi,\phi,\phi)},
\]
with solution
\[
  H_{3000}(\theta) =e^{3i\omega_{0}\theta}\Delta^{-1}(3i\omega_0)\left(3B\left(\phi,H_{2000}\right)+C\left(\phi,\phi,\phi\right)\right).
\]
On the other hand, the system corresponding to $z^2\bar{z}$ is singular,
\[
  \LHS{i\omega_{0}} jH_{2100} = \rsswp{ B\left(\bar{\phi},H_{2000}\right)+2B\left(\phi,H_{1100}\right)+C\left(\phi,\phi,\bar{\phi}\right)}-2c_{1}(0) j\phi.
\]
The Fredholm solvability condition \cref{eq:FSC} requires that
\[
c_1(0)=\frac{1}{2} p \cdot \left(B\left(\bar{\phi},H_{2000}\right)+2B\left(\phi,H_{1100}\right)+C\left(\phi,\phi,\bar{\phi}\right)\right),
\]
and from \cref{lem:bordered} we then obtain the unique solution satisfying $\PAIR{\SUN{\phi}}{H_{2100}} = 0$ as
\[
H_{2100}(\theta) = \BINV{i\omega_{0}}\left( B\left(\bar{\phi},H_{2000}\right)+2B\left(\phi,H_{1100}\right)+C\left(\phi,\phi,\bar{\phi}\right),-2c_1(0)\right)(\theta).
\]
We continue by collecting the coefficients corresponding to the fourth-order terms $z^{3}\bar{z}$ and $z^{2}\bar{z}^{2}$ in the homological equation. The corresponding nonsingular systems may be solved using \cref{lem:regular_solution} and the fact that $\Re(c_1(0))=0$. For $H_{2200}$ this easily gives
\begin{align*}
  H_{2200}(\theta) &= \Delta^{-1}(0)[2B(\bar{\phi},H_{2100})  + 2B(\phi,\BAR{H}_{2100}) + B(\BAR{H}_{2000},H_{2000})\\
                   & \qquad + 2B(H_{1100},H_{1100}) + C(\phi,\phi,\BAR{H}_{2000}) + 4C(\phi,\bar{\phi},H_{1100}) \\
                   & \qquad + C(\bar{\phi},\bar{\phi},H_{2000}) +  D(\phi,\phi,\bar{\phi},\bar{\phi})],
\end{align*}
but for $H_{3100}$ the solution is a bit more subtle. The linear system is
\begin{align*}
\LHS{2i\omega_{0}} jH_{3100} &= \left[B\left(\bar{\phi},H_{3000}\right)+3B\left(\phi,H_{2100}\right)+3B\left(H_{1100},H_{2000}\right) + 3C\left(\phi,\bar{\phi},H_{2000}\right) \right. \\
 & \qquad \left. + 3C\left(\phi,\phi,H_{1100}\right)+D\left(\phi,\phi,\phi,\bar{\phi}\right) \right] \rss -6 c_1(0) jH_{2000},
\end{align*}
so \cref{lem:regular_solution} applies with $w_0 = [\cdots] - 6c_1(0)H_{2000}(0)$ and $w = -6c_1(0)H_{2000}$ and we find
\[
  \begin{aligned}
    H_{3100}(\theta) &=e^{2i\omega_0\theta}\Delta^{-1}(2i\omega_0)[B(\bar{\phi},H_{3000}) + 3B(\phi,H_{2100}) + 3B(H_{1100},H_{2000}) \\
    & \qquad + 3C(\phi,\bar{\phi},H_{2000}) + 3C(\phi,\phi,H_{1100}) + D(\phi,\phi,\phi,\bar{\phi})]\\
    & \qquad - 6c_1(0)\Delta^{-1}(2i\omega_0)[\Delta'(2i\omega_0) - \theta\Delta(2i\omega_0)]H_{2000}(\theta).
  \end{aligned}
\]
The critical normal form coefficient $c_{2}(0)$ is calculated by applying \cref{eq:FSC} to the singular linear system corresponding to the fifth-order term $z^{3}\bar{z}^{2}$ in the homological equation. This gives
\begin{align*}
	c_{2}(0) & =\frac{1}{12}p \cdot \left[2B\left(\bar{\phi},H_{3100}\right)+3B\left(\phi,H_{2200}\right)+B\left(\overline{H}_{2000},H_{3000}\right) \right.\\
 &\qquad+6B\left(H_{1100},H_{2100}\right)+3B\left(\overline{H}_{2100},H_{2000}\right)\\
 &\qquad+6C\left(\bar{\phi},H_{1100},H_{2000}\right)+6C\left(\phi,\bar{\phi},H_{2100}\right)+C\left(\bar{\phi},\bar{\phi},H_{3000}\right)\\
 &\qquad+3C\left(\phi,\phi,\overline{H}_{2100}\right)+3C\left(\phi,\overline{H}_{2000},H_{2000}\right)+6C\left(\phi,H_{1100},H_{1100}\right)\\
 &\qquad+D\left(\phi,\phi,\phi,\overline{H}_{2000}\right)+6D\left(\phi,\phi,\bar{\phi},H_{1100}\right)+3D\left(\phi,\bar{\phi},\bar{\phi},H_{2000}\right)\\
	&\qquad \left. +E\left(\phi,\phi,\phi,\bar{\phi},\bar{\phi}\right)\right].
\end{align*}
The second Lyapunov coefficient is now given by $\ell_2(0)=\frac{1}{\omega_0} \Re(c_2(0))$.

\subsubsection{Parameter-related coefficients}
Next we derive the parameter-related coefficients that provide a linear approximation to the parameter transformation. Following \cite{Kuznetsov2008} and \cite{Beyn2002149} we first expand the eigenvalue $\lambda(\alpha)$ and $c_{1}(\alpha)$ in the normal form \cref{eq:GH_nf-1} in the original parameters $\alpha$ and truncate to fourth order,
\[
\dot{z}=\left(i\omega_{0}+\gamma_{110}\alpha_{1}+\gamma_{101}\alpha_{2}\right)z+\left(c_{1}(0)+\gamma_{210}\alpha_{1}+\gamma_{201}\alpha_{2}\right)z|z|^{2}+c_{2}(0)z|z|^{4}.
\]
The parameters $\alpha$ and $\beta$ are related via
\begin{equation}
\label{eq:GH_alpha}
\alpha=\left(\Re\left(\begin{array}{cc}
\gamma_{110} & \gamma_{101}\\
\gamma_{210} & \gamma_{201}
\end{array}\right)\right)^{-1}\beta
\end{equation}
so that
\begin{equation}
\label{eq:GH_derivative_b_1}
\dfrac{\partial b_{1}}{\partial\beta_{2}}(0)=\Im\left(\left(\begin{array}{cc}
\gamma_{110} & \gamma_{101}\end{array}\right)\left(\Re\left(\begin{array}{cc}
\gamma_{110} & \gamma_{101}\\
\gamma_{210} & \gamma_{201}
\end{array}\right)\right)^{-1}\left(\begin{array}{c}
0\\
1
\end{array}\right)\right).
\end{equation}
The homological equation \cref{eq:homological_equation} becomes
\[
\SUNSTAR{A}j\mathcal{H}(z,\bar{z},\alpha)+J_{1}\alpha \rss +R(\mathcal{H}(z,\bar{z},\alpha),\alpha)=j\left(D_{z}\mathcal{H}(z,\bar{z},\alpha)\dot{z}+D_{\bar{z}}\mathcal{H}(z,\bar{z},\alpha)\dot{\bar{z}} \right),
\]
where $\mathcal{H}$ admits the expansion
\begin{equation}
\label{eq:GH_H_alpha}
u~=~\mathcal{H}(z,\bar{z},\alpha_{1},\alpha_{2}) = z\phi+\bar{z}\bar{\phi} + H_{0010}\alpha_1 +H_{0001}\alpha_2~
 +\sum_{j+k+|\mu| \geq 2}\dfrac{1}{j!k!\mu!}H_{jk\mu}z^{j}\bar{z}^{k}\alpha^{\mu}
\end{equation}
and $R$ is given by \cref{eq:R_truncated}. Notice that the coefficients $H_{jk\mu}$ with $\mu=(00)$ in \cref{eq:GH_H_alpha} coincide with those in the expansion \cref{eq:GH-H}. Collecting the coefficients of the terms $\alpha$ and $z\alpha$ in the homological equation yields the systems
\begin{align*}
\LHSZ jH_{00\mu} & =J_{1}v_{\mu} \; \rss,\\
\LHS{i\omega_{0}} j H_{10\mu}
	& = \rsswp{A_1\left(\phi,v_{\mu}\right)+B\left(\phi,H_{00\mu}\right)}-\gamma_{1\mu}j\phi,
\end{align*}
where $\mu=(10),(01)$ and $v_{10}=(1,0)^{T},v_{01}=(0,1)^{T}$. We solve these systems using \cref{sec:solvability}. By the first part of \cref{lem:regular_solution} the first system has the (constant) solutions
\[
H_{00\mu}(\theta)=\Delta^{-1}(0)J_{1}v_{\mu}
\]
and \cref{eq:FSC} gives
\[
\gamma_{1\mu}=p\left(A_{1}\left(\phi,v_{\mu}\right)+B\left(\phi,H_{00\mu}\right)\right).
\]
Using \cref{lem:bordered} we obtain the solutions
\[
H_{10\mu}(\theta)=\BINV{i\omega_{0}}(A_{1}\left(\phi,v_{\mu}\right)+B\left(\phi,H_{00\mu}\right),-\gamma_{1\mu})(\theta)
\]
for the second equation. To determine $\gamma_{2\mu}$ we first collect the coefficients corresponding to the $z^{2}\alpha$ and $z\bar{z}\alpha$ terms in the homological equation. We obtain the equations
\begin{align*}
\LHS{2i\omega_{0}}jH_{20\mu}
	& = \left[A_{1}\left(H_{2000},v_{\mu}\right)
		+2B\left(\phi,H_{10\mu}\right) + B\left(H_{2000},H_{00\mu}\right)
		+B_1\left(\phi,\phi,v_{\mu}\right) \right. \\
 	& \left. \qquad
 	+ C\left(\phi,\phi,H_{00\mu}\right) \right] \rss -2\gamma_{1\mu}jH_{2000}, \\
\LHSZ jH_{11\mu}
	& = \left[A_{1}\left(H_{1100},v_{\mu}\right)
		+2\Re\left(B\left(\bar{\phi},H_{10\mu}\right)\right)+B\left(H_{1100},H_{00\mu}\right)
		+B_1\left(\phi,\bar{\phi},v_{\mu}\right) \right. \\
	& \left. \qquad 
	+ C\left(\phi,\bar{\phi},H_{0\mu}\right) \right] \rss - 2\Re(\gamma_{1\mu}) jH_{1100}.
\end{align*}
\cref{lem:regular_solution} implies that solutions of the first two equations are given by
\begin{align*}
    H_{20\mu}(\theta) &=e^{2i\omega_0\theta}\Delta^{-1}(2i\omega_0)\left[A_1\left(H_{2000},v_{\mu}\right)+2B\left(\phi,H_{10\mu}\right) +B\left(H_{2000},H_{00\mu}\right)+B_1\left(\phi,\phi,v_{\mu}\right) \right.\\
    & \qquad \left.+C\left(\phi,\phi,H_{00\mu}\right)\right] -2\gamma_{1\mu}\Delta(2i\omega_0)^{-1}\left(\Delta'(2i\omega)-\theta\Delta(2i\omega_0)\right) H_{2000}(\theta), \\
    H_{11\mu}(\theta) & =\Delta^{-1}(0)\left[A_{1}\left(H_{1100},v_{\mu}\right)+2\Re\left(B\left(\bar{\phi},H_{10\mu}\right)\right)+B\left(H_{1100},H_{00\mu}\right)+B_1\left(\phi,\bar{\phi},v_{\mu}\right) \right.\\
& \qquad \left.+C\left(\phi,\bar{\phi},H_{0\mu}\right)\right]-2\Re(\gamma_{1\mu})\Delta(0)^{-1}\left(\Delta'(0)-\theta\Delta(0)\right)\,H_{1100}(\theta).
\end{align*}
Applying \cref{eq:FSC} to $z^2\bar{z}\alpha$ terms in the homological equation results in
\begin{align*}
\gamma_{2\mu} & =\frac{1}{2}~p\cdot \left[A_{1}\left(H_{2100},v_{\mu}\right)+B\left(\bar{\phi},H_{20\mu}\right)+2B\left(\phi,H_{11\mu}\right) \right.\\
 & \qquad+B\left(H_{2100},H_{00\mu}\right)+B\left(H_{2000},\bar H_{10\mu}\right)+2B\left(H_{1100},H_{10\mu}\right)\\
 & \qquad+B_1\left(H_{2000},\bar{\phi},v_{\mu}\right)+2B_1\left(\phi,H_{1100},v_{\mu}\right)+2C\left(\phi,\bar{\phi},H_{10\mu}\right)\\
 & \qquad+C\left(H_{2000},\bar{\phi},H_{00\mu}\right)+C\left(\phi,\phi,H_{01\mu}\right)+2C\left(\phi,H_{1100},H_{00\mu}\right)\\
 & \qquad+\left. C_{1}\left(\phi,\phi,\bar{\phi},v_{\mu}\right)+D\left(\phi,\phi,\bar{\phi},H_{00\mu}\right)\right].
\end{align*}

\subsubsection{Hopf and LPC predictors} \label{sec:genh_predictors}
Now we are ready to specify the predictors for the original parameter-dependent DDE \cref{eq:pd-DDE}. To approximate the Hopf parameter values $\alpha$ and the corresponding equilibrium, we  merely substitute $\beta$ from \cref{eq:GH_approximation_Hopf} into \cref{eq:GH_alpha}, and then put the result together with $z=0$ into the expansion \cref{eq:GH_H_alpha}.

To approximate the LPC parameter values, we substitute $\beta$ from \cref{eq:GH_approximation_LPC} into \cref{eq:GH_alpha}. The cycle period is approximated by \cref{eq:GH_approximation_LPC_T} with \cref{eq:GH_derivative_b_1}. To obtain a predictor for the periodic orbit in the phase space, we set $z=\epsilon e^{i\psi}$ into \cref{eq:GH_H_alpha} using the obtained $\alpha$ values. Truncating to the second order in $\epsilon$ then yields
\begin{align*}
u & =2\Re(e^{i\psi}\phi)\epsilon+\left(H_{1100}-2\Re(c_{2}(0))H_{0001}+\Re(e^{2i\psi}H_{2000})\right)\epsilon^{2}, \qquad \psi\in[0,2\pi].
\end{align*}

\subsection{Fold-Hopf bifurcation\label{sec:fold-Hopf}}
Since the eigenvalues \cref{eq:FH_eigenvalues} are simple, there exist eigenfunctions $\phi_{0,1}$ and $\SUN{\phi_{0,1}}$ satisfying
\[
A\phi_0=0,\qquad A\phi_1=i\omega_0\phi_1,\qquad \STAR{A}\SUN{\phi_0}=0,\qquad \STAR{A}\SUN{\phi_1}=i\omega_0\SUN{\phi_1},
\]
as well as the mutual normalization condition
\[
\PAIR{\SUN{\phi_i}}{\phi_j} =\delta_{ij},\qquad 0\leq i,j\leq 1.
\]
The eigenfunctions $\phi_{0,1}$ and $\SUN{\phi_{0,1}}$ can be explicitly computed using \cref{eq:eigenfunction,eq:eigenfunction1} with $q_{0}\in \RR^{n}$, $q_{1} \in \CC^n$, $p_{0} \in \RRR{n}$ and $p_{1} \in \CCC{n}$ satisfying
\[
\Delta(0)q_{0}=0, \qquad \Delta(i\omega_{0})q_{1}=0,\qquad p_{0}\Delta(0)=0,\qquad p_{1}\Delta(i\omega_{0})=0,
\]
as well as
\[
  p_0\Delta'(0)q_0 = 1, \qquad p_1\Delta'(i\omega_0)q_1 = 1.
\]
Any point $y\in X_{0}$ in the real critical eigenspace can be represented as
\[
y=z_{0}\phi_{0}+z_{1}\phi_{1}+\bar{z}_{1}\bar{\phi}_{1},\qquad (z_0, z_1) \in \RR \times \CC,
\]
where $z_{0}=\PAIR{\SUN{\phi_0}}{y} $ and $z_1=\PAIR{\SUN{\phi_1}}{y}$. Therefore, the homological equation \cref{eq:homological_equation} can be written as
\begin{equation}
\begin{multlined}
\label{eq:FH_HOM}
\SUNSTAR{A}j\mathcal{H}(z,\beta)+J_{1}(\beta)\rss+R(\mathcal{H}(z,\beta),K(\beta))\\
=j\left( D_{z_{0}}\mathcal{H}(z,\beta)\dot{z}_{0}+D_{z_{1}}\mathcal{H}(z,\beta)\dot{z}_{1}+D_{\bar{z}_{1}}\mathcal{H}(z,\beta)\dot{\bar{z}}_{1}\right),
\end{multlined}
\end{equation}
where $z=(z_0,z_1,\bar{z}_1)$, $\beta=(\beta_1,\beta_2)$ and $\dot{z}$ is given by the normal form \cref{eq:FH-nf}. Here, the mapping $\mathcal{H}$ admits the expansion
\begin{equation}
\label{eq:FH_H}
\begin{array}{rcl}
\mathcal{H}(z_{0},z_{1},\bar{z}_{1},\beta) &=& z_{0}\phi_{0}+z_{1}\phi_{1}+\bar{z}_{1}\bar{\phi}_{1}
+ H_{00010}\beta_1+H_{00001}\beta_2\\
&+&{\displaystyle \sum_{j+k+l+|\mu|\geq 2}\dfrac{1}{j!k!l!\mu!}H_{jkl\mu}z_{0}^{j}z_{1}^{k}\bar{z}_{1}^{l}\beta^{\mu}},
\end{array}
\end{equation}
and the functions $K$ and $R$ are as in \cref{eq:K_truncated,eq:R_truncated}, respectively.

\subsubsection{Critical normal form coefficients}
We start by computing the critical normal form coefficients following \cite{Janssens:Thesis}. Collecting the quadratic terms $z_{0}^{2}$, $z_{1}^{2}$, $z_{0}z_{1}$ and $z_{1}\bar{z}_{1}$ we obtain one nonsingular and three singular linear systems. By  \cref{eq:FSC} the singular systems are consistent if and only if
\begin{align*}
g_{200}(0)=\frac{1}{2} & p_{0}B(\phi_{0},\phi_{0}),\qquad g_{110}(0)=p_{1}B(\phi_{0},\phi_{1}),\qquad g_{011}(0)=p_{0}B(\phi_{1},\bar{\phi}_{1}).
\end{align*}
This yields the three quadratic normal form coefficients. The corresponding solutions may then be obtained using \cref{lem:regular_solution,lem:bordered}. Namely,
\begin{align*}
H_{20000}(\theta) & = \BINV{0}(B(\phi_{0},\phi_{0}),-2g_{200}(0))(\theta),\\
H_{02000}(\theta) & = e^{2i\omega_{0}\theta}\Delta^{-1}(2i\omega_{0})B(\phi_{1},\phi_{1}),\\
H_{11000}(\theta) & = \BINV{i\omega_{0}}(B(\phi_{0},\phi_{1}),-g_{110}(0))(\theta),\\
H_{01100}(\theta) & = \BINV{0}(B(\phi_{1},\bar{\phi}_{1}),-g_{011}(0))(\theta).
\end{align*}
For the four remaining cubic normal form coefficients, we collect the coefficients of the resonant terms $z_{0}^{j}z_{1}^{k}\bar{z}_{1}^{l}$  in \cref{eq:FH_HOM} with $j+k+l=3$. This yields four singular linear systems. As before, by \cref{eq:FSC} these systems are consistent if and only if
\begin{align*}
g_{300}(0) & =\frac{1}{6}p_{0} \left(3B(\phi_{0},H_{20000})+C(\phi_{0},\phi_{0},\phi_{0})\right),\\
g_{111}(0) & =p_{0} \left(B(\phi_{0},H_{01100})+B(\phi_{1},\bar{H}_{11000})+B(\bar{\phi}_{1},H_{11000})+C(\phi_{0},\phi_{1},\bar{\phi}_{1})\right),\\
g_{210}(0) & =\frac{1}{2}p_{1} \left(2B(\phi_{0},H_{11000})+B(\phi_{1},H_{20000})+C(\phi_{0},\phi_{0},\phi_{1})\right),\\
g_{021}(0) & =\frac{1}{2}p_{1} \left(2B(\phi_{1},H_{01100})+B(\bar{\phi}_{1},H_{02000})+C(\phi_{1},\phi_{1},\bar{\phi}_{1})\right).
\end{align*}

\subsubsection{Parameter-related coefficients}
The parameter-related linear terms in \cref{eq:FH_HOM} give
\begin{align*}
\LHSZ jH_{00010} & = J_{1}K_{10}\rss-j\phi_{0},\\
\LHSZ jH_{00001} & = J_{1}K_{01}\rss.
\end{align*}
Let $\gamma=(\gamma_{1},\gamma_{2})=p_{0}^{T}J_{1}$. Then by \cref{eq:FSC} we obtain the orthogonal frame
\begin{equation}
K_{10}=s_{1}+\delta_{1}s_{2},\qquad K_{01}=\delta_{2}s_{2},\label{eq:FH_Ks}
\end{equation}
where
\[
s_{1}^{T}=\gamma/\|\gamma\|^{2},\qquad s_{2}^{T}=(-\gamma_{2},\gamma_{1})
\]
and $\delta_{1,2}\in \RR$ are constants. Using \cref{lem:bordered} from \cref{{sec:solvability}} we get
\begin{equation}
\begin{aligned}
H_{00010}(\theta) & =\INV{\Delta}(0)\left(J_{1}K_{10}-\Delta'(0)q_{0}\right)+\delta_{3}q_{0}+\theta q_{0}\\
 & =r_{1}+\delta_{1}r_{2}+\delta_{3}q_{0}-r_{3}(\theta),\\
H_{00001}(\theta) & = \delta_{2}r_{2}+\delta_{4}q_{0},
\end{aligned}\label{eq:FH_H000mu}
\end{equation}
where
\[
r_{1}=\INV{\Delta}(0)\left(J_{1}s_{1}\right),\qquad r_{2}=\INV{\Delta}(0)\left(J_{1}s_{2}\right),\qquad r_{3}(\theta)=\INV{\Delta}(0)\left(\Delta'(0)q_{0}\right)-\theta q_{0},
\]
and the real constants $\delta_{3}$ and $\delta_{4}$ are not chosen such that $\PAIR{\SUN{\phi_0}}{H_{00010}} =0$ and $\PAIR{\SUN{\phi_0}}{H_{00001}} =0$, but will be determined below. Collecting the $z_{0}\beta$ and $z_{1}\beta$ terms in the homological equation yields the systems
\begin{align} \label{eq:FH_secondorder_systems}
\begin{split}
\LHSZ jH_{10010} & = \rsswp{B(\phi_{0},H_{00010}) +A_{1}(\phi_{0},K_{10})} - jH_{20000},\\
\LHSZ jH_{10001} & = \rsswp{B(\phi_{0},H_{00001})+A_{1}(\phi_{0},K_{01})},\\
\LHS{i\omega_{0}} jH_{01010} & = \rsswp{B(\phi_{1},H_{00010})+A_{1}(\phi_{1},K_{10})}
	-j\left(i\omega_{1} \phi_{1}+H_{11000}\right),\\
\LHS{i\omega_{0}} jH_{01001}
	& = \rsswp{B(\phi_{1},H_{00001})+A_{1}(\phi_{1},K_{01})}-\left(1+i\omega_{2}\right) j\phi_{1}.
\end{split}
\end{align}
Notice that the coefficients $\omega_{1,2}$ were introduced in \cref{Sec:ZH_NS_predictors}. To determine $\delta_i(i=1,2,3,4)$ we substitute \cref{eq:FH_Ks,eq:FH_H000mu} into \cref{eq:FH_secondorder_systems}. Then by \cref{eq:FSC} we obtain the system
\begin{multline*}
\left(\begin{array}{cc}
p_{0}B(\phi_{0},r_{2})+p_{0}A_{1}(\phi_{0},s_{2}) & 2 g_{200}(0)\\
\Re(p_{1}B(\phi_{1},r_{2})+p_{1}A_{1}(\phi_{1},s_{2})) & \Re(g_{110}(0))
\end{array}\right) \left(\begin{array}{cc}
\delta_{1} & \delta_2\\
\delta_{3} & \delta_4
\end{array}\right) = \\ \left(\begin{array}{cc}
-p_{0}\left(A_{1}(\phi_{0},s_{1})+B(\phi_{0},r_{1}-r_3)\right)		& 0\\
-\Re(p_{1}\left(A_{1}(\phi_{1},s_{1})+B(\phi_{1},r_{1}-r_3))\right)	& 1
\end{array}\right).
\end{multline*}
Subsequently, the coefficients $\omega_{1}$ and $\omega_{2}$ are given by
\begin{align*}
\omega_{1} & =\Im\left(p_{1}B(\phi_{1},H_{00010})+p_{1}A_{1}(\phi_{1},K_{10})\right),\\
\omega_{2} & =\Im\left(p_{1}B(\phi_{1},H_{00001})+p_{1}A_{1}(\phi_{1},K_{01})\right).
\end{align*}

\subsubsection{Hopf, fold, and Neimark-Sacker predictors} \label{sec:fold-Hopf_predictors}
To approximate the fold and Hopf curves and their corresponding equilibria, one should substitute the expressions for $\beta$ and the equilibrium coordinates given in \cref{Sec:ZH_FoldHopf_predictors} into the expansions \cref{eq:FH_H} and \cref{eq:K_truncated}.

To approximate the periodic orbit at the Neimark-Sacker bifurcation, we substitute $z_{1}=\epsilon e^{i\psi}$ and \cref{eq:FH_NS_predictor} into \cref{eq:FH_H}. After a truncation this gives
\begin{align*}
u & =2\Re\left(e^{i\psi}\phi_{1}\right)\epsilon+\left(\frac{\Re(g_{110}(0))\left(2\Re(g_{021}(0))+g_{111}(0)\right)-2\Re(g_{021}(0))g_{200}(0)}{2g_{200}(0)}H_{00001} \right. \\
	& \phantom{=} \left. -g_{011}(0)H_{00010}+H_{01100}-\left(\frac{2\Re(g_{021}(0))+g_{111}(0)}{2g_{200}(0)}\right)\phi_{0}+\Re\left(e^{2i\psi}\bar{H}_{02000}\right)\right)\epsilon^{2},
\end{align*}
where $\psi \in [0,2\pi]$.

\subsection{Hopf-Hopf bifurcation}
\label{sec:HH_coef}
Since the eigenvalues \cref{eq:HH_eigenvalues} are simple, there exist eigenfunctions $\phi_{1,2}$ and $\SUN{\phi_{1,2}}$,
\begin{equation}
A\phi_{1}=i\omega_{1}\phi_{1},\qquad A\phi_{2}=i\omega_{2}\phi_{2},\qquad \STAR{A}\SUN{\phi_1}=i\omega_1\SUN{\phi_1},\qquad \STAR{A}\SUN{\phi_2}=i\omega_2\SUN{\phi_2},\label{eq:HH_eigenfunctions}
\end{equation}
satisfying the mutual normalization conditions
\[
\PAIR{\SUN{\phi_i}}{\phi_j} = \delta_{ij}, \qquad 1 \leq i,j \leq 2.
\]
The eigenfunctions $\phi_{1,2}$ and $\SUN{\phi_{1,2}}$ can be explicitly computed using \cref{eq:eigenfunction,eq:eigenfunction1} with $q_{1,2} \in \CC^n$ and $p_{1,2} \in \CCC{n}$ such that both
\[
\Delta(i\omega_{1})q_{1}=0,\qquad\Delta(i\omega_{2})q_{2}=0,\qquad p_{1}\Delta(i\omega_{1})=0,\qquad p_{2}\Delta(i\omega_{2})=0,
\]
as well as
\[
  p_1\Delta'(i\omega_1)q_1 = 1, \qquad p_2\Delta'(i\omega_2)q_2 = 1.
\]
Any point $y\in X_{0}$ in the real critical eigenspace can be represented as
\[
y=z_{1}\phi_{1}+\bar{z}_{1}\bar{\phi}_{1}+z_{2}\phi_{2}+\bar{z}_{2}\bar{\phi}_{2},\qquad z_{1,2}\in \CC,
\]
where $z_1=\PAIR{\SUN{\phi_1}}{y}$ and $z_2=\PAIR{\SUN{\phi_2}}{y}$. Therefore, the homological equation \cref{eq:homological_equation} can be written as
\begin{equation}
\begin{multlined}
\SUNSTAR{A}\mathcal{H}(z,\beta)+J_{1}(\beta) \rss +R(\mathcal{H}(z,\beta),K(\beta))\\
=j\left(D_{z_{1}}\mathcal{H}(z,\beta)\dot{z}_{1}+D_{\bar{z}_{1}}\mathcal{H}(z,\beta)\dot{\bar{z}}_{1}+D_{z_{2}}\mathcal{H}(z,\beta)\dot{z}_{2}+D_{\bar{z}_{2}}\mathcal{H}(z,\beta)\dot{\bar{z}}_{2}\right),
\end{multlined}
\label{eq:HH_homological_eq}
\end{equation}
where $z=(z_{1},\bar{z}_{1},z_{2},\bar{z}_{2})$, $\beta=(\beta_1,\beta_2)$ and $\dot{z}$ is given by the normal form \cref{eq:HH_nf-1}. The mapping $\mathcal{H}$ admits the expansion
\begin{align}
\mathcal{H}(z_{1},\bar{z}_{1},z_{2},\bar{z}_{2},\beta_{1},\beta_{2})& = z_{1}\phi_{1}+\bar{z}_{1}\bar{\phi}_{1}+z_{2}\phi_{2}+\bar{z}_{2}\bar{\phi}_{2}+H_{000010}\beta_1 + H_{000001}\beta_2 \nonumber \\
 & +\sum_{j+k+l+m+|\mu| \geq 2}\dfrac{1}{j!k!l!m!\mu!}H_{jklm\mu}z_{1}^{j}\bar{z}_{1}^{k}z_{2}^{l}\bar{z}_{2}^{m}\beta^{\mu} \label{eq:H_expansion-1-2}
\end{align}
and the functions $K$ and $R$ are as in \cref{eq:K_truncated,eq:R_truncated}, respectively.

\subsubsection{Critical normal form coefficients}
For initialization of the Neimark-Sacker curves \cref{eq:HH_NS_asymptotics} we need the cubic critical normal form coefficients $g_{2100}(0)$, $g_{1011}(0)$, $g_{1110}(0)$ and $g_{0021}(0)$. We compute these coefficients following \cite{Janssens:Thesis}.

Collecting the coefficients of the quadratic terms $|z_{1}|^{2},z_{1}^{2},z_{1}z_{2},|z_{2}|^{2},z_{1}\bar{z}_{2}$ and $z_{2}\bar{z}_{1}$ in the homological equation, we obtain six nonsingular linear systems. By \cref{lem:regular_solution} their solutions are
\begin{align*}
H_{110000}(\theta)& = \Delta^{-1}(0)B(\phi_{1},\bar{\phi}_{1}),\\
H_{200000}(\theta)& = e^{2i\omega_{1}\theta}\Delta^{-1}(2i\omega_{1})B(\phi_{1},\phi_{1}),\\
H_{101000}(\theta)& = e^{i\left(\omega_{1}+\omega_{2}\right)\theta}\Delta^{-1}(i\left(\omega_{1}+\omega_{2}\right))B(\phi_{1},\phi_{2}),\\
H_{001100}(\theta)& = \Delta^{-1}(0)B(\phi_{2},\bar{\phi}_{2}),\\
H_{100100}(\theta)& = e^{i\left(\omega_{1}-\omega_{2}\right)\theta}\Delta^{-1}(i\left(\omega_{1}-\omega_{2}\right))B(\phi_{1},\bar{\phi}_{2}),\\
H_{002000}(\theta)& = e^{2i\omega_{2}\theta}\Delta^{-1}(2i\omega_{2})B(\phi_{2},\phi_{2}).
\end{align*}
The desired cubic critical normal form coefficients are obtained by collecting the coefficients of the resonant cubic terms $z_{1}|z_{1}|^{2}$, $z_1|z_{2}|^{2}$, $|z_{1}|^{2}z_2$ and $|z_{2}|^{2}z_{2}$ in the homological equation. This leads to four singular linear systems. By \cref{eq:FSC} these systems are solvable if and only if
\begin{align*}
g_{2100}(0) & =\frac{1}{2}p_{1}\left(2B(\phi_{1},H_{110000})+B(\bar{\phi}_{1},H_{200000})+C(\phi_{1},\phi_{1},\bar{\phi}_{1})\right),\\
g_{1011}(0) & =p_{1}\left(B(\bar{\phi}_{2},H_{101000})+B(\phi_{1},H_{001100})+B(\phi_{2},H_{100100})+C(\phi_{1},\phi_{2},\bar{\phi}_{2})\right),\\
g_{1110}(0) & =p_{2}\left(B(\bar{\phi}_{1},H_{101000})+B(\phi_{1},\bar{H}_{100100})+B(\phi_{2},H_{110000})+C(\phi_{1},\bar{\phi}_{1},\phi_{2})\right),\\
g_{0021}(0) & =\frac{1}{2}p_{2}\left(2B(\phi_{2},H_{001100})+B(\bar{\phi}_{2},H_{002000})+C(\phi_{2},\phi_{2},\bar{\phi}_{2})\right).
\end{align*}

\subsubsection{Parameter-related coefficients}
The linear terms in \cref{eq:HH_homological_eq} give back the eigenfunctions \cref{eq:HH_eigenfunctions} and the parameter-related equations
\[
\LHSZ jH_{0000\mu}=J_{1}K_{\mu} \rss,
\]
where $\mu=(10),(01)$. Let
\begin{equation}
K_{\mu}=\gamma_{1\mu}e_{1}+\gamma_{2\mu}e_{2},\label{eq:HH_Kmu}
\end{equation}
where $e_{1}=(1,0)$, $e_{2}=(0,1)$ and $\gamma_{i\mu}(i=1,2)\in \RR$ are constants to be determined. Then \cref{lem:regular_solution} from \cref{{sec:solvability}} implies
\begin{align}
H_{0000\mu}(\theta) & =\gamma_{1\mu}\Delta^{-1}(0)J_{1}e_{1}+\gamma_{2\mu}\Delta^{-1}(0)J_{1}e_{2}.\label{eq:HH_H000mu}
\end{align}
Collecting the $z_{i}\beta_{j}$-terms with $1\leq i,j\leq 2$ yields the systems
\begin{align}
\begin{split}\label{eq:HH_secondorder_systems}
\LHS{i\omega_{1}} jH_{100010}
	& = \rsswp{A_{1}(\phi_{1},K_{10})+B(\phi_{1},H_{000010})} - (1+ib_{11})j\phi_{1},\\
\LHS{i\omega_{1}} jH_{100001}
	& =\rsswp{A_{1}(\phi_{1},K_{01})+B(\phi_{1},H_{000001})}-ib_{12} j\phi_{1},\\
\LHS{i\omega_{2}} jH_{001010}
	& =\rsswp{A_{1}(\phi_{2},K_{10})+B(\phi_{2},H_{000010})}-ib_{21} j\phi_{2},\\
\LHS{i\omega_{2}} jH_{001001}
	& =\rsswp{A_{1}(\phi_{2},K_{01})+B(\phi_{2},H_{000001})}-(1+ib_{22})j\phi_{2},
\end{split}
\end{align}
where $b_{jk}$ are defined in \cref{Sec:FH_NSs}. To determine $\gamma_{i\mu}(i=1,2)$ we substitute \cref{eq:HH_Kmu,eq:HH_H000mu} into \cref{eq:HH_secondorder_systems}. Then by \cref{eq:FSC} we obtain the system
\[
\Re\left[\left(\begin{array}{cc}
\Gamma_{11} & \Gamma_{12}\\
\Gamma_{31} & \Gamma_{32}
\end{array}\right)\right]\left(\begin{array}{cc}
\gamma_{110} & \gamma_{210}\\
\gamma_{101} & \gamma_{201}
\end{array}\right)=\left(\begin{array}{cc}
1 & 0 \\
0 & 1
\end{array}\right),
\]
where
\[
\Gamma_{ij} \DEF A_1(\phi_i,e_j) + B(\phi_i,\Delta^{-1}(0) J_1 e_j), \qquad 1\leq i,j \leq 2.
\]
Note that $\Delta^{-1}(0)J_{1}e_{i}$ is a constant function of $\theta$.

It now follows from \cref{eq:HH_secondorder_systems} that the coefficients $b_{11},b_{12},b_{21}$ and $b_{22}$, needed for the second order approximation of the periods, are given by
\begin{align*}
b_{11} & =\Im\left(p_{1}\left(A_{1}(\phi_{1},K_{10})+B(\phi_{1},H_{000010})\right)\right),\\
b_{12} & =\Im\left(p_{1}\left(A_{1}(\phi_{1},K_{01})+B(\phi_{1},H_{000001})\right)\right),\\
b_{21} & =\Im\left(p_{2}\left(A_{1}(\phi_{2},K_{10})+B(\phi_{2},H_{000010})\right)\right),\\
b_{22} & =\Im\left(p_{2}\left(A_{1}(\phi_{2},K_{01})+B(\phi_{2},H_{000001})\right)\right).
\end{align*}

\subsubsection{Hopf and Neimark-Sacker predictors}
To approximate the Hopf curves and their corresponding equilibria, one should substitute the expressions for $\beta$ and the equilibrium coordinates given in \cref{Sec:FH_Hopfs} into the expansions \cref{eq:H_expansion-1-2} and \cref{eq:K_truncated}.

To approximate the Neimark-Sacker periodic orbits, we substitute $(z_1,z_2)=(\epsilon e^{i\psi_{1}},0)$ and \cref{eq:HH_NS1_pm}, and $(z_1,z_2)=(0,\epsilon e^{i\psi_{2}})$ and \cref{eq:HH_NS2_pm} into \cref{eq:H_expansion-1-2}. After a truncation, we obtain
\begin{align*}
	u_1 & =2\Re\left(e^{i\psi_{1}}\phi_{1}\right)\epsilon+\Bigl(-\Re(g_{1110}(0))H_{000001}-\Re(g_{2100}(0))H_{000010} \Bigr.\\
	& \phantom{=} \left.+H_{110000}+\Re\left(e^{2i\psi_{1}}H_{200000}\right)\right)\epsilon^{2}, \qquad \psi_1\in[0,2\pi]
\end{align*}
and
\begin{align*}
u_2 & =2\Re\left(e^{i\psi_{2}}\phi_{2}\right)\epsilon+\Bigl(-\Re(g_{0021}(0))H_{000001}-\Re(g_{1011}(0))H_{000010}\Bigr.\\
	& \phantom{=} \left.+H_{001100}+\Re\left(e^{2i\psi_{2}}H_{002000}\right)\right)\epsilon^{2}, \qquad \psi_2\in[0,2\pi].
\end{align*}

\subsection{Transcritical-Hopf bifurcation} \label{sec:transcritical-Hopf}
Compared with the fold-Hopf bifurcation in \cref{sec:fold-Hopf}, the eigenvalues, eigenfunctions, the homological equation, and the functions $\mathcal{H}$, $K$ and $R$ remain unchanged. It is only the ODE on the center manifold $\dot{z}$ that changes to the normal form \cref{eq:Ht-nf}. Furthermore, also the critical normal form coefficients for the transcritical-Hopf bifurcation remain the same as for the fold-Hopf bifurcation. Therefore, we proceed only with the parameter-related equations.

Collecting the coefficients of the $z_{0}\beta$ and $z_{1}\beta$ terms in the homological equation we obtain the systems
\begin{align} \label{eq:TH_secondorder_systems}
\begin{split}
\LHSZ jH_{10010} & = A_{1}(\phi_{0},K_{10})\rss - j\phi_{0},\\
\LHSZ jH_{10001} & = A_{1}(\phi_{0},K_{01}) \rss, \\
\LHS{i\omega_{0}} jH_{01010} & = A_{1}(\phi_{1},K_{10})-i\omega_{1}j\phi_{1} \rss, \\
\LHS{i\omega_{0}} jH_{01001} & = A_{1}(\phi_{1},K_{01})-(1+i\omega_{2})j\phi_{1} \rss.
\end{split}
\end{align}
Let
\begin{equation} \label{eq:TH_Kmu}
K_{\mu}=\gamma_{1\mu}e_{1}+\gamma_{2\mu}e_{2}, \qquad \mu=(10),(01),
\end{equation}
where $e_{1}=(1,0)$, $e_{2}=(0,1)$ and $\gamma_{i\mu}(i=1,2)\in \RR$. To determine $\gamma_{i\mu}(i=1,2)$ we substitute \cref{eq:TH_Kmu} into \cref{eq:TH_secondorder_systems}. Then by \cref{eq:FSC} we obtain the system
\begin{align*}
\left(\begin{array}{cc}
p_{0}A_{1}(\phi_{0},e_{1}) & p_{0}A_{1}(\phi_{0},e_{2})\\
\Re \left(p_{1}A_{1}(\phi_{1},e_{1}) \right) & \Re\left(p_{1}A_{1}(\phi_{1},e_{2})\right)
\end{array}\right) & \left(\begin{array}{cc}
\gamma_{110} & \gamma_{210}\\
\gamma_{101} & \gamma_{201}
\end{array}\right)=\left(\begin{array}{cc}
1 & 0\\
0 & 1
\end{array}\right).
\end{align*}
In order to make the last two systems in \cref{eq:TH_secondorder_systems} consistent we must have that
\begin{align*}
\omega_{1}& = \Im\left(p_{1}A_{1}(\phi_{1},K_{10})\right),\\
\omega_{2}& = \Im\left(p_{1}A_{1}(\phi_{1},K_{01})\right).
\end{align*}

\subsubsection{Neimark-Sacker predictors}
The predictors for the Hopf and transcritical bifurcation curves, as well as those for the Neimark-Sacker bifurcation curves (including the cycle periods), can be easily obtained using the asymptotics from \cref{Sec:HT_HopfTrans_predictors,Sec:HT_NS_predictors}. In particular,
to approximate the periodic orbits along the Neimark-Sacker curves, we substitute $z_{1}=\epsilon e^{i\psi}$ and \cref{eq:HT_NS_predictor} into \cref{eq:FH_H}. This gives the following linear approximations:

\begin{align*}
u & =\left(\mp \sqrt{\frac{ g_{011}(0)}{g_{200}(0)}}\phi_0
+2\Re\left(e^{i\psi}\phi_{1}\right)\right)\epsilon \qquad \psi\in[0,2\pi].
\end{align*}

\section{Computation of derivatives for discrete DDEs}
\label{sec:Implement}
All predictors described in the previous sections are implemented in version 3.2a of \DDEBIFTOOL for models of the type \cref{Eq:FiniteDDE}. The discrete DDE \cref{Eq:FiniteDDE} is a particular instance of \cref{eq:pd-DDE} with $h=\tau_m$ and
$$
F(\phi, \alpha) = f(\Xi \phi, \alpha),
$$
where the linear evaluation operator $\Xi: X \to \RR^{n \times (m+1)}$ is defined by
\begin{equation}
  \label{Eq:Handle}
  \Xi \phi \DEF \left(\phi(-\tau_0),\phi(-\tau_{1}),\dots,\phi(-\tau_{m})\right).
\end{equation}
with the convention $\tau_0 \DEF 0$. In particular, by the chain rule,
\[
D_1F(0,0)\phi = D_1f(0,0)\Xi\phi = \sum_{j=0}^m{D_{1,j}f(0,0)\phi(-\tau_j)}, \qquad \phi \in X,
\]
with $M_j \DEF D_{1,j}f(0,0) \in \RR^{n \times n}$ the partial derivative of $f$ at the origin with respect to its $j$th state argument. So, if \cref{Eq:FiniteDDE} has a equilibrium at the origin for $\alpha = 0$, then the linear part of the splitting \cref{eq:DDE-RHS} at $\alpha = 0$ is precisely the right-hand side of the above equation. Therefore $\zeta : [0,h] \to \RR^{n \times n}$ must be such that
$$
\PAIR{\zeta}{\phi} = \sum_{j=0}^m{M_j\phi(-\tau_j)}, \qquad \forall\,\phi \in X.
$$
Hence $\zeta$ has jump discontinuities $M_j$ at the points $\tau_j$ for $j = 0,\ldots,m$ and is constant otherwise. So, in this case the characteristic matrix \cref{eq:CharMatrix} is given by
\[
\Delta(z) = z I  - \sum_{j=0}^m  M_j {\rm e}^{-z \tau_j}, \qquad z \in \CC.
\]
The multilinear forms appearing in \cref{eq:R_truncated} can be expressed in terms of the derivatives of the function  $f: \RR^{n \times (m+1)} \times \RR^p \to \RR^n$ from \cref{Eq:FiniteDDE}. For $r, s \ge 0$ with $r + s \ge 1$ the mixed derivative of order $r + s$ of $f$ at $(0,0)$ is an $(r + s)$-linear form on $[\RR^{n \times (m + 1)}]^r \times [\RR^p]^s$, with the understanding that at most one factor may be absent in case $r = 0$ or $s = 0$. Let $Q, Q^{1},\ldots,Q^{r}$ be matrices in $\RR^{n \times (m + 1)}$ and let $\alpha, \alpha^{1},\ldots,\alpha^{s}$ be vectors in $\RR^p$. Then this derivative acts as
\begin{multline}
  \label{eq:derivsf}
  D_1^r D_2^s f(0,0)(Q^1,\ldots,Q^r, \alpha^1,\ldots,\alpha^s) =\\
  \sum_{j,k,\ell} \left.\frac{\partial^{r + s} f(Q, \alpha)}{\partial q_{j_1k_1}\ldots \partial q_{j_rk_r}\partial \alpha_{\ell_1} \ldots \partial \alpha_{\ell_s}}\right|_{(Q, \alpha) = (0,0)} q^1_{j_1k_1}\cdots q^r_{j_rk_r}\alpha^1_{\ell_1}\cdots\alpha^s_{\ell_s},
\end{multline}
where the multidimensional sum runs over
\[
1 \le j_1, \ldots, j_r \le n, \qquad 0 \le k_1, \ldots, k_r \le m, \qquad 1 \le \ell_1, \ldots, \ell_s \le p.
\]
The multilinear forms appearing in \cref{eq:R_truncated}, as well as \cref{eq:J_1}, are computed from \cref{eq:derivsf} by composition with $\Xi$ from \cref{Eq:Handle} as
\begin{equation*}
\label{eq:multilinearforms}
D_1^rD_2^s F(0,0)(\phi_1,\ldots,\phi_r, \alpha_1, \ldots, \alpha_s) = D_1^rD_2^sf(0,0)(\Xi \phi_1,\ldots, \Xi \phi_r, \alpha_1, \ldots, \alpha_s),
\end{equation*}
for $\phi_1,\ldots,\phi_r \in X$ and $\alpha_1,\ldots, \alpha_s \in \RR^p$. For given $r$ and $s$ the multidimensional array of partial derivatives inside the sum in \cref{eq:derivsf} is of course symmetric under permutation of the state indices $j_1k_1,\ldots,j_rk_r$ and the parameter indices $\ell_1,\ldots,\ell_s$. This can be exploited for efficient storage and access.

\section{Examples\label{sec:Examples}}
In this section we will demonstrate the correctness of the normal form coefficients and the accuracy of the predictors in four different models. We do this twofold. Firstly, by comparing the predictors in parameter-space with the computed in \DDEBIFTOOL bifurcation curves, and, secondly, by performing simulations near the bifurcation point under consideration. The simulation is done either with the build-in routine \blist{dde23} of \MATLAB or with the Python package \texttt{pydelay} \cite{Flunkert2009Flunkert}. The latter gives significant speed performance when considering simulation over longer time intervals. This usually is the case when one wants to demonstrate the existence of stable invariant manifolds. Since in this section only the main results are given, we provide details (including simulation results) in the \hyperlink{mysupplement}{Supplement}. Furthermore, the source code of the examples has been included into the \DDEBIFTOOL software package. This will hopefully provide a good starting point when considering other models.

\subsection{Generalized Hopf bifurcation in a coupled FHN neural system with delay}
\label{sec:example_FHN}

In \cite{Xu2010} the following system is considered
\begin{equation}
\begin{cases}
\begin{aligned}
\dot{u}_{1}(t) & =-\dfrac{u^3_{1}(t)}{3}+(c+\alpha)u^2_{1}(t)+du_{1}(t)-u_{2}(t)+2\beta f(u_{1}(t-\tau)),\\
\dot{u}_{2}(t) & =\varepsilon(u_{1}(t)-bu_{2}(t)).
\end{aligned}
\end{cases}\label{eq:DDE_FHN}
\end{equation}
Here $(u_{1},u_{2})$ is the completely synchronous solution of the three coupled FitzHugh\textendash Nagumo (FHN) neuron system
\begin{equation}
\begin{cases}
\begin{aligned}
\dot{u}_{1}(t) & =-\dfrac{u^{3}_{1}(t)}{3}+(c+\alpha)u_{1}^{2}(t)+du_{1}(t)-u_{2}(t)+\beta\left[f(u_{3}(t-\tau))+f(u_{5}(t-\tau))\right],\\
\dot{u}_{2}(t) & =\varepsilon(u_{1}(t)-bu_{2}(t)),\\
\dot{u}_{3}(t) & =-\dfrac{u^{3}_{3}(t)}{3}+(c+\alpha)u_{3}^{2}(t)+du_{3}(t)-u_{4}(t)+\beta\left[f(u_{3}(t-\tau))+f(u_{5}(t-\tau))\right],\\
\dot{u}_{4}(t) & =\varepsilon(u_{3}(t)-bu_{4}(t)),\\
\dot{u}_{5}(t) & =-\dfrac{u^{3}_{5}(t)}{3}+(c+\alpha)u_{5}^{2}(t)+du_{5}(t)-u_{6}(t)+\beta\left[f(u_{1}(t-\tau))+f(u_{3}(t-\tau))\right],\\
\dot{u}_{6}(t) & =\varepsilon(u_{5}(t)-bu_{6}(t)),
\end{aligned}
\end{cases}\label{eq:thee_coupled_neurons}
\end{equation}
where $\alpha,\beta$ measure the synaptic strength in self-connection and neighborhood-interaction, respectively. The parameters $b$ and $\epsilon$ are assumed to be positive such that $0<b<1$ and $0<\epsilon\ll1$. The function $f$ is a sufficiently smooth sigmoid amplification function and $\tau>0$ represents the time delay in signal transmission. For the derivation of \cref{eq:DDE_FHN} from the system \cref{eq:thee_coupled_neurons}, as well as for stability conditions of the completely synchronous solution, we refer to  \cite{Xu2010}. In that article a generalized Hopf point was analyzed using the traditional formal adjoint method and the two-step center manifold reduction, see \cite{Hale@1977}. Numerical simulations where made to confirm their results. For this $(\beta,\alpha)$ are taken as the unfolding parameters and the parameters
\[
b=0.9,\qquad\varepsilon=0.08,\qquad c=2.0528,\qquad d=-3.2135,\qquad\tau=1.7722
\]
are fixed. The sigmoid amplification function $f(u)=\tanh(u)$ is used.

\begin{figure}[htbp]
\centering
\ifcompileimages
  \tikzsetnextfilename{FHN_bifdia}%
  \input{tikz/FHN_bifdia}%

\else
\includegraphics{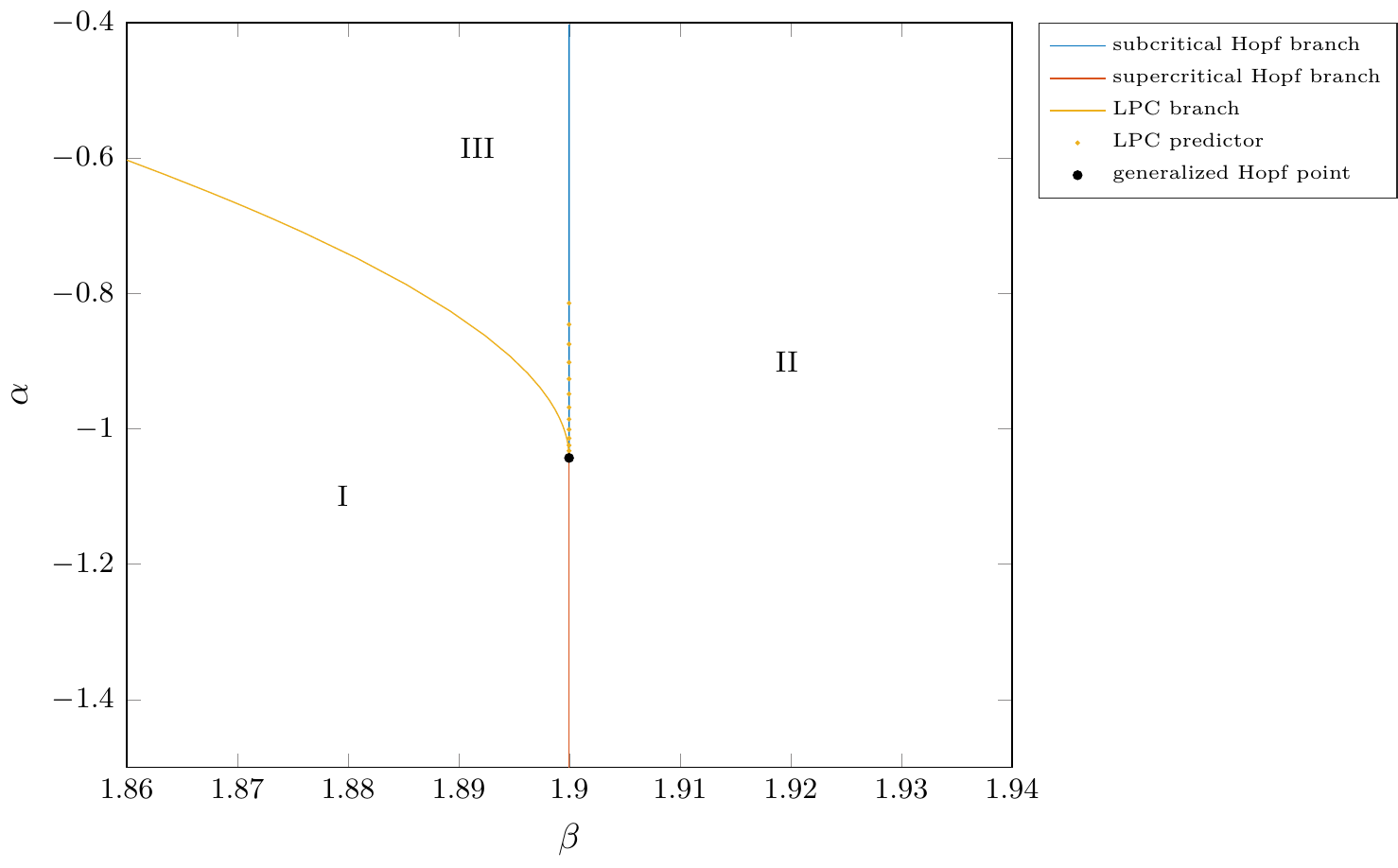}
\fi
\caption{Bifurcation diagram near the generalized Hopf point in the system \cref{eq:DDE_FHN} with unfolding parameters $(\beta,\alpha)$. The bifurcation curves are nearly identical to those in the bifurcation diagram of the topological normal form as presented in \textup{\cite[page 314]{Kuznetsov2004}}.}
\label{fig:FHN-bifurcation-diagram}
\end{figure}

According to \cite{Xu2010}, a generalized Hopf point is present at the origin with the parameter values $(\beta,\alpha)=(1.9,-0.9710)$. We took this point and calculated its stability and the corresponding normal form coefficients. Although we do confirm that the point under consideration is a Hopf point, the first Lyapunov coefficient does not vanish and we conclude that the point cannot be a generalized Hopf point. However, the simulation in \cite{Xu2010} do suggest a generalized Hopf point for nearby parameter values. Therefore we continued the Hopf point in $(\beta,\alpha)$. Then a generalized Hopf point is located at $(\beta,\alpha)=(1.9,-1.0429)$ with \emph{negative} second Lyapunov coefficient $\ell_2(0)=-15.6733$, indicating the existence of a stable steady state inside a unstable cycle, which in turn is located inside a stable cycle. We remark that the second Lyapunov coefficient found in \cite{Xu2010} is \emph{positive}. This contradicts the simulation of the dynamics made in the same article. Indeed when the second Lyapunov coefficient is positive a time-reversal must be taking into account when considering the bifurcation diagram in the case the second Lyapunov coefficient is negative, see \cite{Kuznetsov2004}. Then the situation of a stable steady state inside a stable cycle (separated by an unstable cycle) does not occur.

Using the predictors from \cref{sec:GH_predictors} the fold and LPC bifurcation curves emanating from the generalized Hopf point. In \cref{fig:FHN-bifurcation-diagram} the resulting bifurcation diagram is shown.

\subsection{Fold-Hopf bifurcation of the Rose\textendash Hindmarsh model with time delay}
\label{sec:ex_Rose_Hindmarsh}
In \cite{Ma2011} a Rose-Hindmarsh model \cite{Hindmarsh1982,Hindmarsh1984} with time delay in the self-feedback process, which takes the form
\begin{equation}
\begin{cases}
\begin{aligned}
\dot{x}(t)& = y(t)-ax^{3}(t)+bx^{2}(t-\tau)-cz(t)+I_{app},\\
\dot{y}(t)& = c-dx^{2}(t)-y(t),\\
\dot{z}(t)& = r(S(x(t)-\chi)-z(t)),
\end{aligned}
\end{cases}\label{eq:Rose-Hindmarsh}
\end{equation}
is considered. Here $x$ represents membrane potential, $y$ represents a recovery variable, $z$ denotes the adaption current, and $a,b,c,d>0,S$ and $\chi$ are real constants. The external current $I_{app}$ and $r$ are control parameters, and $\tau$ denotes the synaptic transmission delay. The constants $a,b,c,d,\chi$ and $r$ are fixed. Let $(x_{\star},y_{\star},z_{\star})$ be a steady state of \cref{eq:Rose-Hindmarsh}, then
\begin{equation} \label{eq:ystar_zstar}
y_{\star}=c-dx_{\star}^{2},\qquad z_{\star}=S(x_{\star}-\chi).
\end{equation}
The conditions for a fold-Hopf bifurcation have been derived in \cite{Ma2011} analytically. Indeed, let $S$ be arbitrary and set
\begin{align}
x_{\star} & =\frac{1}{3a}\left(b-d\pm\sqrt{\left(b-d\right)^{2}-3acS}\right), \nonumber \\
I_{app} & =x_{\star}^{2}(ax_{\star}-b+d)+c(S(x_{\star}-\chi)-1), \label{eq:I_app}  \\
A & =x_{\star}^{2}\left((3ax_{\star}+2d)^{2}-4b^{2}\right)-2rx_{\star}(2dx_{\star}-1)(3ax_{\star}-2b+2d) \nonumber \\
 & \qquad+r^{2}(4dx_{\star}(-2bx_{\star}+dx_{\star}-1)+1),\nonumber \\
B & =9a^{2}x_{\star}^{4}+2rx_{\star}(3ax_{\star}-2b+2d)-4b^{2}x_{\star}^{2}-4dx_{\star}+r^{2}+1,\nonumber  \\
\omega_{1,2} & =\sqrt{-B\pm\sqrt{B^{2}-4A}}. \nonumber
\end{align}
Then a fold-Hopf bifurcation occurs when
\[
\tau=\begin{cases}
\frac{1}{\omega_{1,2}}\left(\arcsin Y+2k\pi\right), & Z\geq0,\\
\frac{1}{\omega_{1,2}}\left(\pi-\arcsin Y+2k\pi\right), & Z\leq0,
\end{cases}
\]
where $k=0,1,2,\dots$ and
\begin{align*}
Y &= \frac{\omega_{1,2}}{2b}\left(\frac{r  (2b-2d-3ax_\star)}{r^{2}+\omega_{1,2}^{2}}+\frac{2d}{\omega_{1,2}^{2}+1}-\frac{1}{x_{\star}}\right), \\
Z &= \frac{\omega_{1,2}}{2b}\left(\frac{r^2(2b-2d-3ax_\star)}{r^2+\omega_{1,2}^2}+\frac{2 d}{\omega_{1,2}^2+1} +3ax_\star \right).
\end{align*}
In \cite{Ma2011} the parameters values
\begin{equation}
a=1.0,\qquad b=3.0,\qquad c=1.0,\qquad d=5.0,\qquad\chi=-1.6,\qquad r=0.001\label{eq:rose_hindmarsh_pm1}
\end{equation}
are fixed. It follows that a fold-Hopf bifurcation is located at
\[
x_{\star}=0.1308,\qquad S=-0.57452592,\qquad\tau=5.768830916,
\]
and $I_{app}$, $(y_\star,z_\star)$ given by \cref{eq:I_app} and \cref{eq:ystar_zstar}, respectively. To unfold the singularity the `parameters' $(x_{\star},S)$ are used, see \cite{Ma2011}. Here we will take the more natural unfolding parameter $(I_{app},S)$. Calculating the stability with \DDEBIFTOOL gives the eigenvalues
\[
0.001\pm1.0081i,\qquad-0.000+0.000i.
\]
All other eigenvalues lie in the open left half of the complex plane. Calculating the normal form coefficients reveals that
\[
s = \mbox{sgn}(g_{200}(0)g_{011}(0)) = \mbox{sgn}(1.8487\mathrm{e}{-05}),
\qquad  \theta(0)=\frac{\Re(g_{110}(0))}{g_{200}(0)}=-139.0315
\]
and
\begin{equation} \label{ex:RH:e}
\begin{array}{l}
e(0)=\Re\left[g_{210}(0)+g_{110}(0)\left(\frac{\Re g_{021}(0)}{g_{011}(0)}-\frac{3}{2}\frac{g_{300}(0)}{g_{200}(0)}+\frac{g_{111}(0)}{2g_{011}(0)}\right) -\frac{g_{021}(0) g_{200}(0) }{g_{011}}\right]=15.6941.
\end{array}
\end{equation}
Since $s=1$ and $\theta(0)<0$ , a global bifurcation curve or invariant tori are present for parameters sufficiently close to the bifurcation, see \cite[page 342]{Kuznetsov2004}. However, since the sign of $e(0)$ is \emph{positive} the tori are unstable. Thus according to our analysis the simulated torus in \cite{Ma2011} cannot be attributed to the fold-Hopf bifurcation.
\begin{figure}[ht]
\centering
\ifcompileimages
  \tikzsetnextfilename{RH_bifurcation_diagram_II}%
  \input{tikz/RH_bifurcation_diagram_II}%

\else
\includegraphics{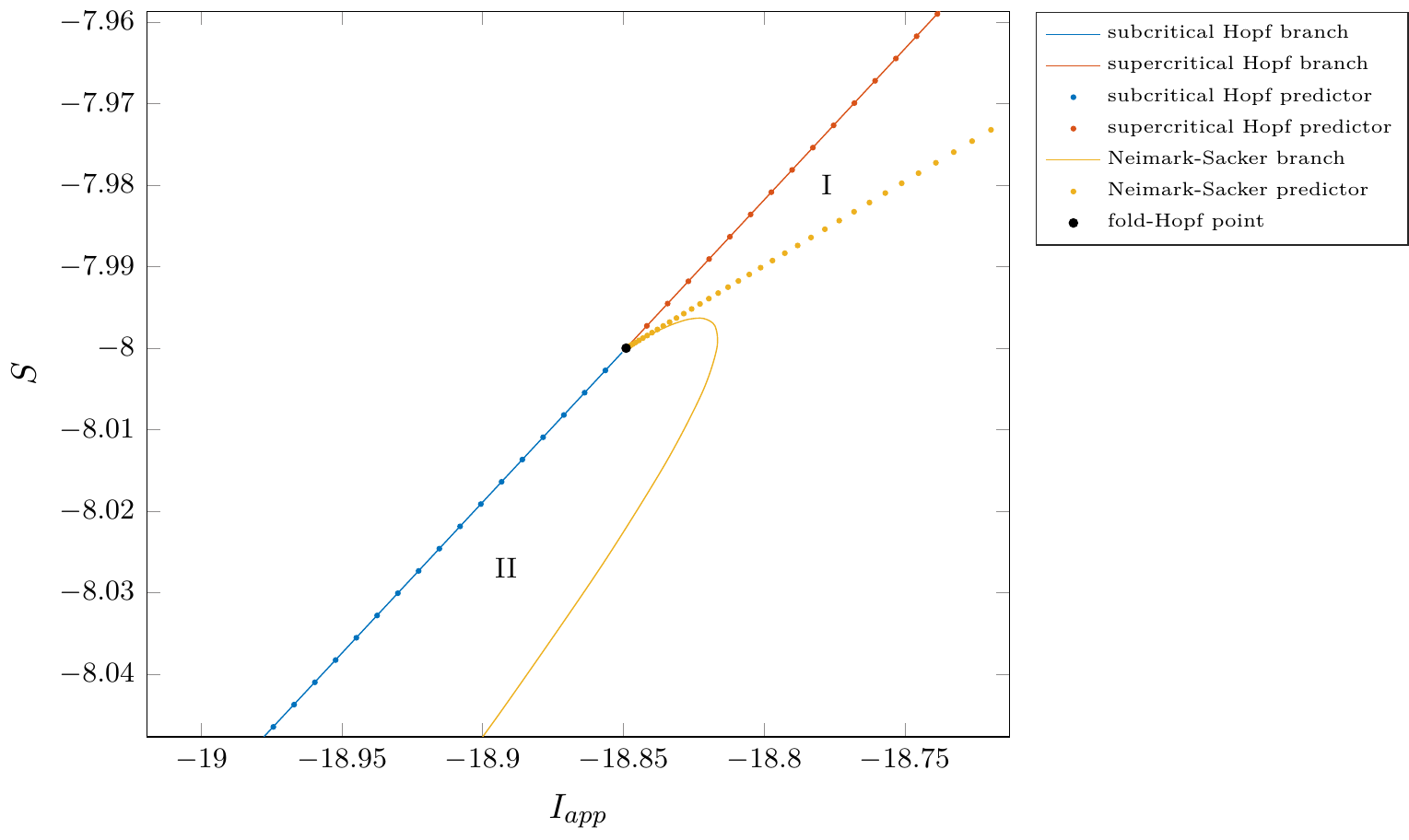}
\fi
\caption{\label{fig:Rose-Hindmarsh_bifurcation_diagram}Bifurcation diagram near the fold-Hopf point in \cref{eq:Rose-Hindmarsh} with $(r,S)=(1.4,-8)$. The fold branch is not included here since it is indistinguishable from the Hopf curve at this scale.
}
\end{figure}

For demonstration purposes, we take the parameters $r=1.4$ and $S=-8$, while keeping the other parameters as in \cref{eq:rose_hindmarsh_pm1}. Then a fold-Hopf bifurcation is located at $x_{\star}=1.0972$, $\tau=0.9402$, $I_{app}$ as in \cref{eq:I_app}, and $(y_\star,z_\star)$ given by \cref{eq:ystar_zstar}. The leading eigenvalues become
\[
0.000\pm5.6042i, \qquad 0.000+0.000i,
\]
while the normal form coefficients are given by
\[
s=\mbox{sgn}(1.7700), \qquad \theta(0)=-0.1569 \qquad \mbox{and} \qquad e(0)=-0.0378.
\]
Thus the sign of $s$ and $\theta(0)$ remain unchanged. However, since the sign of $e(0)$ is  \emph{negative}, there is a time reversal to take into account. Therefore, we expect a stable torus to be present for nearby parameter values. Using the predictors from \cref{sec:Predictor-FH}, we successfully continued the fold, Hopf, and Neimark-Sacker bifurcation curves emanating from the point, see \cref{fig:Rose-Hindmarsh_bifurcation_diagram}. 

\subsection{Hopf-Hopf and generalized Hopf bifurcations in Active control system}
\label{sec:acs_example}

Active control system is used to control the response of structures to internal or external excitation. The mathematical model with time delay can be described as follows \cite{Peng2013}
\begin{equation}
m\ddot{x}(t)+c\dot{x}(t)+kx(t)+ux(t-\tau)+v\dot{x}(t-\tau)=\tilde{f}(t).\label{eq:acs1}
\end{equation}
Here $x(t)$ is the displacement of the controlled system, $m>0$ is the mass, $c$ and $k$ are the damping and the stiffness, respectively, $\tau$ is the time delay represented in the relative displacement feedback loop and in the relative velocity feedback loop, $u$ and $v$ are feedback strengths, respectively, and $\tilde{f}$ represents the external excitation. Let $t^{*}=\sqrt{k/m}t$, $\zeta=c/2m\sqrt{m/k}$, $g_{u}=u/k,g_{v}=v/m\sqrt{m/k}$ and $f(t)=\tilde{f}(t)/k$ . Then equation \cref{eq:acs1} becomes
\[
\ddot{x}(t)+2\zeta\dot{x}(t)+x(t)+g_{u}x(t-\tau)+g_{v}\dot{x}(t-\tau)=f(t),
\]
where the asterisks are omitted for simplicity. Following \cite{Ding@2016} and \cite{Peng2013} we consider the case when $f$ is replaced by a nonlinear position time delay feedback given  by $\beta x^{3}(t-\tau)$, see also \cite{xu2003vanderPolDuffing}. As in \cite{Ding@2016} we fix the parameters
\[
g_{u}=0.1,\quad g_{v}=0.52,\quad\beta=0.1
\]
and take $\zeta$ and $\tau$ as control parameters.
Let $\dot{x}(t)=y(t)$, then we obtain
\begin{equation}
\begin{cases}
\begin{aligned}
\dot{x}(t)&=\tau y(t),\\
\dot{y}(t)&=\tau\left(-x(t)-g_{u}x(t-1)-2\zeta y(t)-g_{v}y(t-1)+\beta x^{3}(t-1) \right).
\end{aligned}
\end{cases}\label{eq:acs3-1}
\end{equation}
Here the delay is scaled by using the transformation of time $t\rightarrow t/\tau$. In this way the delay can treated as an ordinary parameter.

\begin{figure}[htbp]
\centering
\ifcompileimages
  \tikzsetnextfilename{acs_hoho_predictors}%
  \input{tikz/acs_hoho_predictors}%

\else
\includegraphics{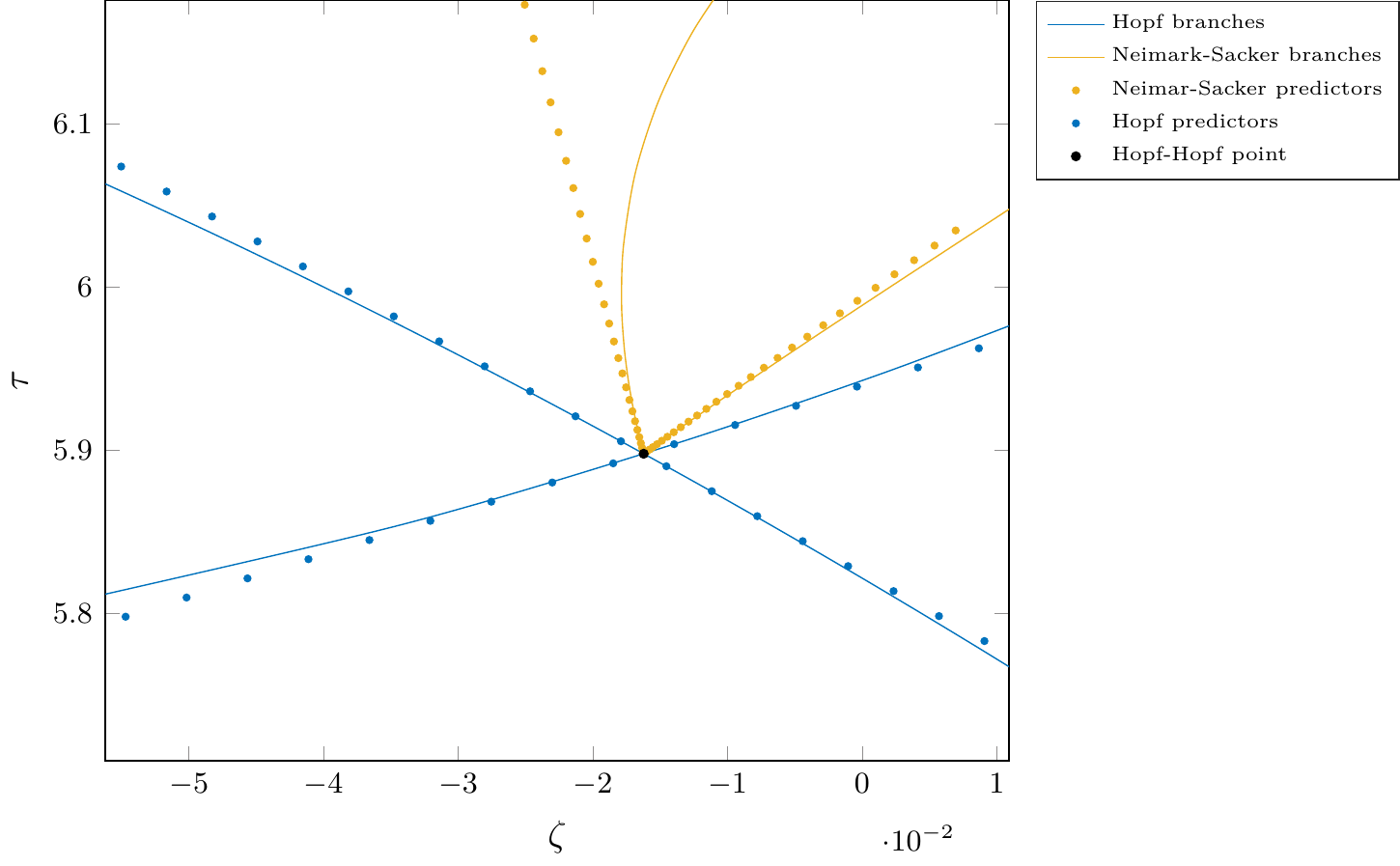}
\fi
\caption{Bifurcation diagram near the Hopf-Hopf point at parameter values \cref{eq:acs-HH-pm} in an active control system with time delay given by \cref{eq:acs3-1}. There are two supercritical Hopf curves (blue) and two Neimark-Sacker curves (yellow).  We see that the predictors (dotted) give good approximations near the codimension two point.}
\label{fig:acs_hoho_predictors}
\end{figure}

The trivial equilibrium undergoes a Hopf-Hopf bifurcation at the parameter values
\begin{equation}
(\zeta_{c},\tau_{c})=(-0.016225,5.89802),
\label{eq:acs-HH-pm}
\end{equation}
see \cite{Ding@2016} for the derivation. Using \DDEBIFTOOL we manually construct the Hopf-Hopf point and compute its stability and normal form coefficients. We obtain the eigenvalues $0.0000\pm4.5275i$ and $-0.0000\pm7.6449i$. The quadratic critical normal form coefficients are
\begin{align*}
g_{2100}(0) & =-0.0915+0.1214i, & g_{1110}(0) & =0.2151+0.3876i,\\
g_{1011}(0) & =-0.3084+0.4096i, & g_{0021}(0) & =0.1813+0.3268i.
\end{align*}
From
\[
(\text{Re }g_{2100}(0))(\text{Re }g_{0021}(0))=-0.0166<0,
\]
we conclude that this Hopf-Hopf bifurcation is of `difficult' type, see \cite{Kuznetsov2004}. Furthermore, since the quantities
\[
\theta=\theta(0)=\frac{\text{Re }g_{1011}(0)}{\text{Re }g_{0021}(0)}=-1.7009,\qquad\delta=\delta(0)=\frac{\text{Re }g_{1101}(0)}{\text{Re }g_{2100}(0)}=-2.3517
\]
are such that $\theta<0,\,\delta<0,\,\theta\delta>1$ it follows that we are in case VI. %
We continue the Neimark-Sacker and Hopf bifurcation curses emanating from the Hopf-Hopf point using the predictors from \cref{sec:HH_pedictors}. In  \cref{fig:acs_hoho_predictors} a close-up is given near the Hopf-Hopf point comparing the computed curves with the predictors in parameter space.

\begin{figure}[ht]
\noindent \begin{centering}
\ifcompileimages
  \tikzsetnextfilename{acs_hoho_genh_lpc_ns}%
  \input{tikz/acs_hoho_genh_lpc_ns}%

\else
\includegraphics{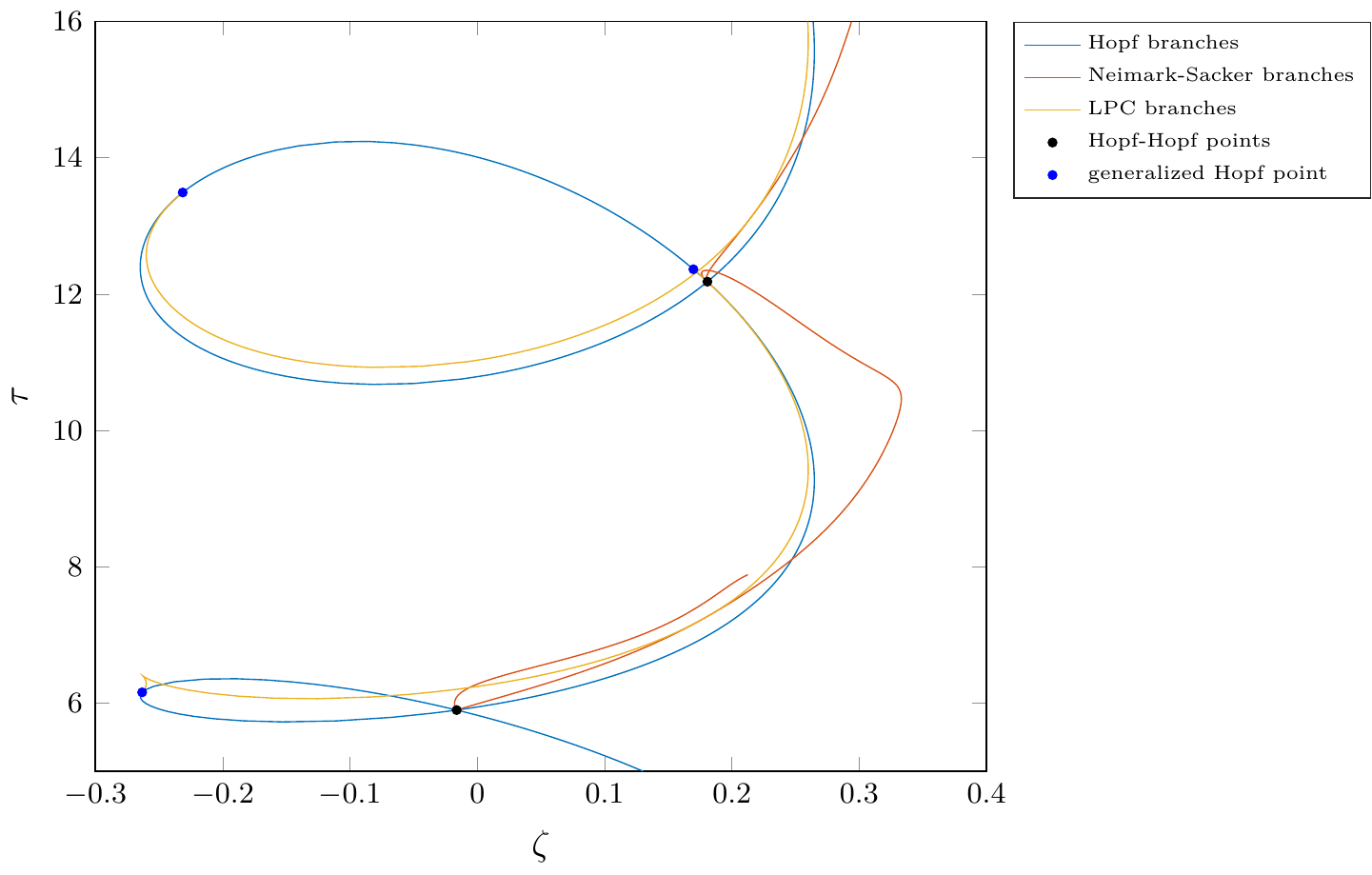}
\fi
\par\end{centering}
\caption{\label{fig:acs_unfolding} Bifurcation diagram obtained by continuing Hopf-Hopf, Neimark-Sacker and LPC bifurcation curves in the active control system \cref{eq:acs3-1} using the predictors from \cref{sec:GH_predictors,sec:HH_pedictors} combined with the continuation capabilities from \DDEBIFTOOL. Two Hopf-Hopf points are connected by a Neimark-Sacker bifurcation curve. Also two of the three generalized Hopf points are connected by a single LPC curve.}
\end{figure}

Using the detection capabilities of \DDEBIFTOOL one additional Hopf-Hopf point and three generalized Hopf points are located on the continued Hopf branches. The normal form coefficients of the second Hopf-Hopf point are such that
\[
(\text{Re }g_{2100}(0))(\text{Re }g_{0021}(0))=1.7331\mathrm{e}{-04}>0
\]
and
\[
\theta\geq\delta>0,\qquad \theta\delta>1.
\]
We conclude that we are in case I of the `simple' type, see \cite[page 360]{Kuznetsov2004}. Therefore, no stable invariant two-dimensional torus is predicted for nearby parameter values, only two stable period orbits expected. Using the predictors from \cref{sec:GH_predictors,sec:HH_pedictors} we can easily continue the codimension one cycle bifurcations from the located degenerate Hopf points, showing complicated bifurcation diagram in  \cref{fig:acs_unfolding}.

\subsection{Transcritical-Hopf bifurcation in Van der Pol's oscillator with delayed position and velocity feedback}
\label{sec:HT_example}
In \cite{Bramburger2014} a generalization of Van der Pol's oscillator with delayed feedback
\begin{equation}\label{eq:dde_vanderPol}
\ddot{x}(t)+\varepsilon(x^{2}(t)-1)\dot{x}(t)+x(t)=g(\dot{x}(t-\tau),x(t-\tau)),
\qquad 0<\tau<\infty,
\end{equation}
is considered. Here $g\in C^{3}$ satisfies the conditions $g(0,0)=0,g_{\dot{x}}(0,0)=a$ and $g_{x}(0,0)=b$. The linearization of equation \cref{eq:dde_vanderPol} around the trivial solution $x=0$ gives
\[
\ddot{x}(t)-\varepsilon\dot{x}(t)+x(t)=a\dot{x}(t-\tau)+bx(t-\tau).
\]
From which we obtain the characteristic equation
\[
\Delta(\lambda,\tau)=\lambda^{2}-\varepsilon\lambda+1-(a\lambda+b)e^{-\lambda\tau}=0.
\]
Let
\begin{equation}
b=1,\qquad\tau=\tau_{0}\neq\varepsilon+a,\qquad\varepsilon^{2}-a^{2}<2,\label{eq:fold-Hopf_conditions}
\end{equation}
then the characteristic equation has a simple zero and a pair of purely imaginary roots $\lambda=\pm i\omega_{0}$. Here $\omega_{0}$ and $\tau_{0}$ are defined by
\[
\omega_{0}=\sqrt{2-\varepsilon^{2}+a^{2}},\qquad\tau_{0}=\frac{1}{\omega_{0}}\arccos\left(\dfrac{1-(1+\varepsilon a)\omega_{0}^{2}}{a^{2}\omega_{0}^{2}+1}\right),
\]
see \cite[Proposition 2.1]{Bramburger2014}. We set the function $g$ to
\begin{align*}
g(\dot{x}(t-\tau),x(t-\tau))=&(1+\mu_{1}) x(t-\tau)-0.2 \dot{x}(t-\tau)-0.2 x(t-\tau)^{2}\\
						&-0.2 x(t-\tau) \dot{x}(t-\tau)-0.2 x(t-\tau)^2+0.5 x(t-\tau)^3
\end{align*}
and $\varepsilon=0.3$. Then the conditions \cref{eq:fold-Hopf_conditions} are satisfied and
\begin{equation}
\omega_{0}\approx1.396424004376894,\qquad\tau_{0}\approx1.757290761249588.\label{eq:vdpo_omega0_tau0}
\end{equation}
\begin{figure}[htbp]
\centering
\ifcompileimages
  \tikzsetnextfilename{VDPO_bifurcation_diagram}%
  \input{tikz/VDPO_bifurcation_diagram}%

\else
\includegraphics{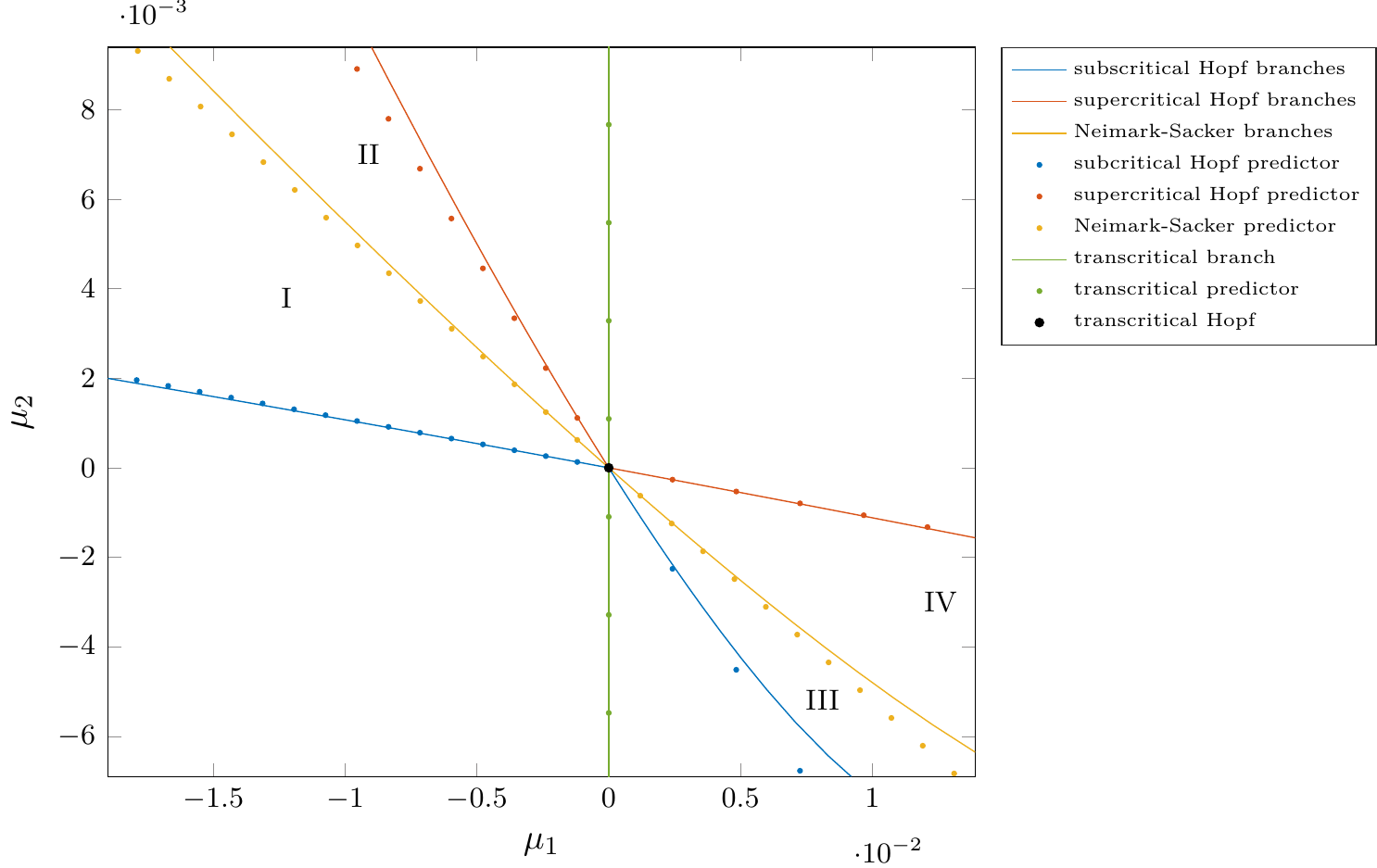}
\fi
\caption{\label{fig:HT_bifurcation_diagram} Bifurcation diagram near the transcritical-Hopf bifurcation in the delayed Van der Pol's oscillator given by \cref{eq:vdp}. There are two supercritical Hopf curves (blue), two subcritical Hopf curves (red), two Neimark-Sacker curves (yellow) and one transcritical curve (green).  We see that the predictors (dotted) give good approximation for nearby values.}
\end{figure}

To analyze the system with \DDEBIFTOOL we set $y(t)=\dot{x}(t)$ and transform the time with $t\rightarrow t/\tau$ to obtain the two-component system
\begin{equation}
\begin{cases}
\begin{aligned}
\dot{x}(t)&=\left(\tau_0+\mu_2\right)y(t),\\[0.5em]
\dot{y}(t)&=\left(\tau_0+\mu_{2}\right)\left[-x(t)-\varepsilon(x^2(t)-1)y(t)+(1+\mu_1)x(t-1)-0.2y(t-1) \right. \\[0.5em]
&\qquad \left. -0.2x^2(t-1)-0.2x(t-1)y(t-1)-0.2y^2(t-1)+0.5x^3(t-1)\right].
\end{aligned}
\end{cases}\label{eq:vdp}
\end{equation}
Here we introduced the unfolding parameters $(\mu_1,\mu_{2}):=(b-1,\tau-\tau_{0})$ to translate the singularity to the origin. One immediately sees that the trivial equilibrium $(\dot{x},x)=(0,0)$ is an equilibrium for all parameter values $(\mu_{1},\mu_{2})$. Therefore, the parameter-dependent normal form for the generic fold-Hopf cannot be used here. Instead the normal form for the transcritical-Hopf bifurcation must be used. Using \DDEBIFTOOL we compute the stability and the normal form coefficients. The leading eigenvalues are $0.000+0.000i$ and $-0.000+2.4539i$, where $2.4539\approx\omega_{0}\tau_{0}$, see \cref{eq:vdpo_omega0_tau0}. Furthermore, the normal form coefficients are such that
\[
g_{011}(0)\times\text{Re\,}\left(g_{110}(0)\right)=0.4241\times \Re\left(-0.1337+0.2672i\right)<0.
\]
Therefore, there are two Neimark-Sacker bifurcation curves predicted, see \cref{sec:HT_predictors}. Using the predictors from \cref{sec:HT_predictors} we continue the transcritical, Hopf and Neimark-Sacker bifurcation curves emanating from the transcritical-Hopf bifurcation point. In \cref{fig:HT_bifurcation_diagram} the bifurcation diagram is shown.

\section{Concluding remarks}
We have provided explicit formulas for the normal form coefficients needed to initialize codimension one equilibrium and nonhyperbolic cycle bifurcations emanating from generalized Hopf, fold-Hopf, Hopf-Hopf and transcritical-Hopf points in DDEs. Applications to four different models were given, confirming the correctness of the derivation of the normal form coefficients and the asymptotics. A paper providing a second-order predictor for the homoclinic orbits emanating from the generic and transcritical codimension two Bogdanov-Takens bifurcations in DDEs, along the lines of \cite{kouznetsov2014improved}, is in preparation.

Our proof of the existence of a smooth parameter-dependent center manifold is given in the general context of perturbation theory for dual semigroups (sun-star calculus). Consequently the applicability of this result extends beyond classical DDEs, although here we did restrict to the case of an eventually compact $\mathcal{C}_0$-semigroup on a sun-reflexive state space. It follows that the results from \cref{sec:pd:duality,sec:spectral_bold,sec:nonlinear_bold,sec:cm_bold} are valid as well for other classes of delay equations such as renewal equations (also known as Volterra functional equations) and systems of mixed type \cite{diekmann2007stability}.

Furthermore, in \cite{VanGils2013,Dijkstra2015} the technique was used to calculate the critical normal form coefficients for Hopf and Hopf-Hopf bifurcations occurring in neural field models with propagation delays. For these models sun-reflexivity is lost, which is typical for delay equations in abstract spaces or with infinite delay. However, it is often possible to overcome this functional analytic complication, so dual perturbation theory can still be employed successfully \cite{Diekmann2008,Diekmann2012blending,VanGils2013,Janssens2019}. It has also been used in the context of semilinear hyperbolic systems \cite{Lichtner2009hyperbolicsystems}.

It is demonstrated - at a formal level - in \cite{Sieber@2017} that the normalization technique described in \cref{sec:normal-forms} still works for DDEs with state-dependent delays. However, as already mentioned in \cref{rem:nonsmooth}, for DDEs the nonlinearity generally does \emph{not} depend differentiably on the delay parameters. Therefore, in the case of state-dependent DDEs it is generally not possible to justify differentiation of the nonlinearity with respect to the state, let alone to rely on higher order smoothness. So, as far as we know there is still no proof of the validity of the normalization technique for state-dependent DDEs.

Returning to the setting of classical DDEs, the most obvious next challenge is to derive normal forms for bifurcations of periodic orbits by generalizing \cite{Kuznetsov2005,DeWitte2013,DeWitte2014} to DDEs. The resulting formulas can then be implemented in \DDEBIFTOOL to facilitate numerical bifurcation analysis of periodic orbits in classical DDEs.

\section*{Acknowledgments}
The authors would like to thank Prof. Odo Diekmann (Utrecht University) for very useful discussions on parameter-dependent perturbation of linear semigroups. We also thank Prof. Peter De Maesschalck (Hasselt University) for supporting this research project.

\ifarxiv
\clearpage
\hypertarget{mysupplement}{}

%
%

\newif\ifsiam


\siamtrue
\siamfalse

\newif\ifarxiv
\ifsiam \arxivfalse \else \arxivtrue \fi


\ifsiam
\documentclass[supplement]{siamonline0516}

\ifpdf
\hypersetup{
  pdftitle={Supplementary Materials: \TheTitle},
  pdfauthor={\TheAuthors}
}
\fi

\externaldocument{main}

\begin{document}
\maketitle
\pagestyle{myheadings}
\thispagestyle{plain}
\markboth{\TheAuthors}{SWITCHING TO NONHYPERBOLIC CYCLES IN DDES}
\fi


\ifarxiv
\begin{mytitle}
  \title{\textsc{supplementary materials for:}\\ \TheTitle}
  \maketitle
\end{mytitle}
\thispagestyle{plain}

\ResetCounters
\pdfbookmark[0]{Supplement}{supplement}
\fi


\newif\ifcompileimages
\compileimagesfalse


In this supplement we provide walkthroughs of the examples given in \cref{sec:Examples} with \DDEBIFTOOL\footnote{\url{http://ddebiftool.sourceforge.net/}} \cite{DDEBIFTOOL}. These walkthroughs enable other researchers to reproduce the results obtained in the main text.

Additionally, we will show the code used for simulation near the bifurcations points under consideration. Either using the build-in routine \lstinline|dde23| from \MATLAB \cite{Shampine01solvingdelay} or the Python package \PYDELAY\footnote{\url{http://pydelay.sourceforge.net/}} \cite{Flunkert2009Flunkert}. Other DDE models, undergoing one of the degenerate Hopf bifurcations treated in this paper, can easily be studied by making minor modifications to the given code.

The focus will be on the initialization and continuation of the various codimension one equilibrium and cycles bifurcation curves emanating from the degenerate Hopf points and on simulation near the bifurcation points. For a complete overview of the capabilities and functionality for \DDEBIFTOOL, we refer to the online tutorials files and also the manual and the references therein.

All code has been included into the \DDEBIFTOOL package version 3.2a on the SourceForge repository and can be executed without the need to copy and paste. Note that the code is tested on \MATLAB 2018b and \OCTAVE 4.2.2. Different results may occur with other versions of \MATLAB and \OCTAVE.

\section{Generalized Hopf bifurcation in a coupled FHN neural system with delay}
A completely synchronous solution of the three coupled FitzHugh-Nagumo (FHN) neuron system is given by the system
\begin{equation}
\begin{cases}
\begin{aligned}
\dot{u}_{1}(t) & =-\dfrac{u_{1}^{3}(t)}{3}+(c+\alpha)u_{1}^{2}(t)+du_{1}(t)-u_{2}(t)+2\beta f(u_{1}(t-\tau)),\\
\dot{u}_{2}(t) & =\varepsilon(u_{1}(t)-bu_{2}(t)).
\end{aligned}
\end{cases}\label{sm:eq:DDE_FHN}
\end{equation}
see \cref{sec:example_FHN} and \cite{Ma2011}. As before, we fix the parameters
\[
b=0.9,\qquad\varepsilon=0.08,\qquad c=2.0528,\qquad d=-3.2135,\qquad\tau=1.7722,
\]
and take for $f:\mathbb{R}\to \mathbb{R}$ the sigmoid amplification function $f(u)=\tanh(u)$. The parameters $(\beta,\alpha)$ are used to unfold the singularity. 

\begin{remark}
This demonstration can be found in the directory \lstinline|demos/tutorial/VII/FHN| relative to the main directory of the \DDEBIFTOOL package.
\end{remark}

\subsection{Generate system files}
Before we start to analyze the system with \DDEBIFTOOL, we first create a \emph{system file}. This file contains the definition of the system \cref{sm:eq:DDE_FHN}, the standard derivatives needed for calculation of the eigenvalues and eigenvectors, the continuation of bifurcation points and cycles, and also the multilinear forms, see \cref{eq:multilinearforms}, used for the calculation of the coefficients of the critical and parameter-dependent normal forms. Alternatively, one can only supply the system itself, see \cref{sm:lst:wo_system_file}. Then finite difference is used to approximate the derivatives. However, this is less efficient and accurate, and therefore not recommended. A separate script \lstinline|gen_sym_FHN.m| is used to create a system file. The most important parts of this script are listed and discussed below.
\begin{lstlisting}[style=customMatlab,escapechar=!]
%% Add paths and load sym package if GNU Octave is used
clear
ddebiftoolpath='../../../../';
addpath(strcat(ddebiftoolpath,'ddebiftool'),...
  strcat(ddebiftoolpath,'ddebiftool_extra_symbolic'));
if dde_isoctave()
  pkg load symbolic
end
!\matlabrule!
%% Create parameter names as strings and define fixed parameters
% The demo has the parameters |beta|, |alpha| and |tau|
parnames={'beta','alpha','tau'};
b=sym(0.9,'r');
epsilon=sym(0.08,'r');
c=sym(2.0528,'r');
d=sym(-3.2135,'r');
!\matlabrule!
%% Create symbols for parameters, states and delays states
% |par| is the array of symbols in the same order as parnames.
% Due to the following two lines we may, for example,
% use either beta or par(1) to refer to the delay.
syms(parnames{:});       % create symbols for beta, alpha and tua
par=cell2sym(parnames);  % now beta is par(1) etc
!\matlabrule!
%% Define system using symbolic algebra
% create symbols for u1(t) u1(t-tau), u2(t), u2(t-tau)
syms u1 u1t u2 u2t
du1_dt=-u1^3/3+(c+alpha)*u1^2+d*u1-u2+2*beta*tanh(u1t);
du2_dt=epsilon*(u1-b*u2);
!\matlabrule!
%% Differentiate and generate code (multi-linear forms)
[fstr,derivs]=dde_sym2funcs(...
[du1_dt;du2_dt],... % n x 1 array of derivative symbolic expressions
[u1,u1t;u2,u2t],... % n x (ntau+1) array of symbols for states (current & delayed)
par,... % 1 x np (or np x 1) array of symbols used for parameters
'filename','sym_FHN_mf',... % optional argument specifying output file
'directional_derivative',false); 
!\matlabrule!
%% Differentiate and generate code (directional derivatives)
[fstr,derivs]=dde_sym2funcs(...
[du1_dt;du2_dt],... % n x 1 array of derivative symbolic expressions
[u1,u1t;u2,u2t],... % n x (ntau+1) array of symbols for states (current & delayed)
par,... % 1 x np (or np x 1) array of symbols used for parameters
'filename','sym_FHN',...  % optional argument specifying output file
'directional_derivative',true);
\end{lstlisting}
The variable \lstinline|ddebiftoolpath| is directed to the \DDEBIFTOOL main folder, which should have been extracted somewhere on the computer. Here a path relative to the current working directory is used. Note that although we only use the parameters $(\beta,\alpha)$ as unfolding parameters, in the current version of \DDEBIFTOOL, we also need to include the delay(s) in the list of parameters. After running the script, the function \lstinline|dde_sym2funcs| creates two system files \lstinline|sym_FHN_mf.m| and \lstinline|sym_FHN.m|. The first file \lstinline|sym_FHN_mf.m| implements the higher order derivatives as multilinear forms, as explained in \cref{sec:Implement}, and therefore the file we will solely be using. The second file \lstinline|sym_FHN.m| uses directional derivatives to implement the higher order derivatives. The directional derivatives approach \emph{formally} allows the use of state-dependent delays, see \cite{Sieber@2017}. Although both approaches yields (up to rounding errors) identical normal form coefficients, multilinear forms are much faster.

\subsection{Loading the \DDEBIFTOOL package}\label{sm:sec:loading_DDE-BIFTool}
Now that a system file is created we continue with \DDEBIFTOOL to analyze \cref{sm:eq:DDE_FHN}. The code in the following sections highlights the import parts of the file \lstinline|FHN.m|.
\DDEBIFTOOL consists of a set of \MATLAB routines. Thus, to start using \DDEBIFTOOL, we only need to add \DDEBIFTOOL directories to the search path.
\begin{lstlisting}[style=customMatlab,caption=Add \DDEBIFTOOL scripts to the search path, label={sm:lst:searchpath}]
%% Clean workspace and add DDE-BifTool scripts to 
% the MATLAB search path
clear;      % clear variables
close all;  % close figures
ddebiftoolpath='../../../../';
addpath(strcat(ddebiftoolpath,'ddebiftool'),...
        strcat(ddebiftoolpath,'ddebiftool_extra_psol'),...
        strcat(ddebiftoolpath,'ddebiftool_extra_nmfm'),...
        strcat(ddebiftoolpath,'ddebiftool_utilities'));
\end{lstlisting}
There are four subdirectories added to the search path:
\par
\medskip
\begin{description}
\item[ddebiftool] Containing the core files of \DDEBIFTOOL.
\item[ddebiftool\_extra\_psol] An extension for enabling continuation of periodic orbit bifurcations for delay-differential equations with constant or state-dependent delay.
\item[ddebiftool\_extra\_nmfm] An extension for normal form computation.
\item[ddebiftool\_utilities] Containing various utilities.
\end{description}

\subsection{Set parameter names}
The following code allows us to use \lstinline[keywordstyle=\color{darkblue}]{ind.beta} instead of remembering the index of the parameter $\beta$ in the parameter array, and similarly for the other parameters.
\begin{lstlisting}[style=customMatlab]
%% Set parameter names
parnames={'beta','alpha','tau'};
cind=[parnames;num2cell(1:length(parnames))];
ind=struct(cind{:});
\end{lstlisting}
In this way, fewer mistakes are likely to be made and the code is easier to read.

\subsection{Initialization}
Next, we set up the \lstinline|funcs| structure, containing information about where the system and its derivatives are stored, a function pointing to which parameters are delays, and various other settings.
\begin{lstlisting}[style=customMatlab]
%% Set the funcs structure
% We load the precalculated multilinear forms. These have been
% generated with the file gen_sym_FHN.m.
funcs=set_symfuncs(@sym_FHN_mf,'sys_tau',@()ind.tau);
\end{lstlisting}
Alternatively, when no system files have been generated, one could initialize the system \cref{sm:eq:DDE_FHN} as follows.
\begin{lstlisting}[style=customMatlab,caption=Define system without a system file, label={sm:lst:wo_system_file}]
%% Define the system
b=0.9; epsilon=0.08; c=2.0528; d=-3.2135; % fixed parameters
FHN_sys = @(xx,par) [...
-xx(1,1,:).^3/3+(c+par(1,ind.alpha,:)).*xx(1,1,:).^2+d*xx(1,1,:)...
-xx(2,1,:)+2*par(1,ind.beta,:).*tanh(xx(1,2,:));
epsilon*(xx(1,1,:)-b*xx(2,1,:))];
%% Set funcs structure
funcs=set_funcs('sys_rhs',FHN_sys,'sys_tau',@()ind.tau,...
    'x_vectorized',true,'p_vectorized',true);			
\end{lstlisting}
Inspecting the output of the \lstinline|funcs| handle gives.
\begin{lstlisting}[style=customBash,keepspaces=true] 
>> funcs

funcs = 

  ?struct? with fields:

                     sys_rhs: @(x,p)wrap_rhs(x,p,funcs.sys_rhs,funcs.x_vectorized,funcs.p_vectorized)
                    sys_ntau: @()0
                     sys_tau: @()ind.tau
                    sys_cond: @dummy_cond
                    sys_deri: @(x,p,nx,np,v)dde_gen_deriv(funcs.sys_dirderi,x,p,nx,np,v,1)
                    sys_dtau: []
                  sys_mfderi: {}
                 sys_dirderi: {[function_handle]  [function_handle]}
                 sys_dirdtau: []
                x_vectorized: 1
                p_vectorized: 1
                        hjac: @(ord)eps^(1/(2+ord))
    sys_unfolding_parameters: []
                      tp_del: 0
           sys_deri_provided: 0
        sys_dirderi_provided: 0
\end{lstlisting}
The output shows that no derivative file is supplied. In this case, the derivatives are calculated using finite-difference approximations with the function \lstinline|dde_dirderiv|. Again, we do not recommend using the latter approach. However, it can be useful for debugging purposes.

\subsection{Stability and normal form coefficients of the generalized-Hopf point}
We manually specify a steady-state at the generalized-Hopf point found in \cite{Ma2011} and calculate its stability.
\begin{lstlisting}[style=customMatlab]
% construct steady-state point
beta0=1.9; alpha0=-0.9710; tau0=1.7722;
stst=dde_stst_create('x',[0;0]);
stst.parameter(ind.beta)  = beta0;
stst.parameter(ind.alpha) = alpha0;
stst.parameter(ind.tau)   = tau0;
% Calculate stability
method=df_mthod(funcs,'stst');
stst.stability=p_stabil(funcs,stst,method.stability);
\end{lstlisting}
Inspecting the \lstinline|stst.stability| structure yields
\begin{lstlisting}[style=matlabConsole]
>> stst.stability.l1(1:6)

ans =

   0.0000 + 0.0720i
   0.0000 - 0.0720i
  -0.0818 + 3.1068i
  -0.0818 - 3.1068i
  -0.3478 + 6.4411i
  -0.3478 - 6.4411i

>> 
\end{lstlisting}
The eigenvalues confirm that the point under consideration is indeed a Hopf point. Next, we convert the steady-state point to a Hopf point and calculate the normal form coefficients with the function \lstinline|nmfm_genh|, which implements the coefficients derived in \cref{sec:GH_coef}.
\begin{lstlisting}[style=customMatlab]
%% Calculate critical normal form coefficients
hopf=p_tohopf(funcs,stst);
method=df_mthod(funcs,'hopf');
hopf.stability=p_stabil(funcs,hopf,method.stability);
genh=p_togenh(hopf);
genh=nmfm_genh(funcs,genh);
\end{lstlisting}
The normal form coefficients are stored in the \lstinline|genh.nmfm| structure.
\begin{lstlisting}[style=matlabConsole,keepspaces=true]
>> genh.nmfm
ans = 

  ?struct? with fields:

    L2: -18.1302
    L1: 0.3980

>> 
\end{lstlisting}
Clearly, the first Lyapunov coefficient (\lstinline|L1|) is nonzero. It follows that the Hopf point is not degenerate. 

\subsection{Continue Hopf point}
\begin{figure}
	\centering
	\subfloat[Hopf curve]{
	    \ifcompileimages
  \tikzsetnextfilename{FHN_Hopf_curve}%
  \input{tikz/FHN_Hopf_curve}%

		\else
		    \includegraphics{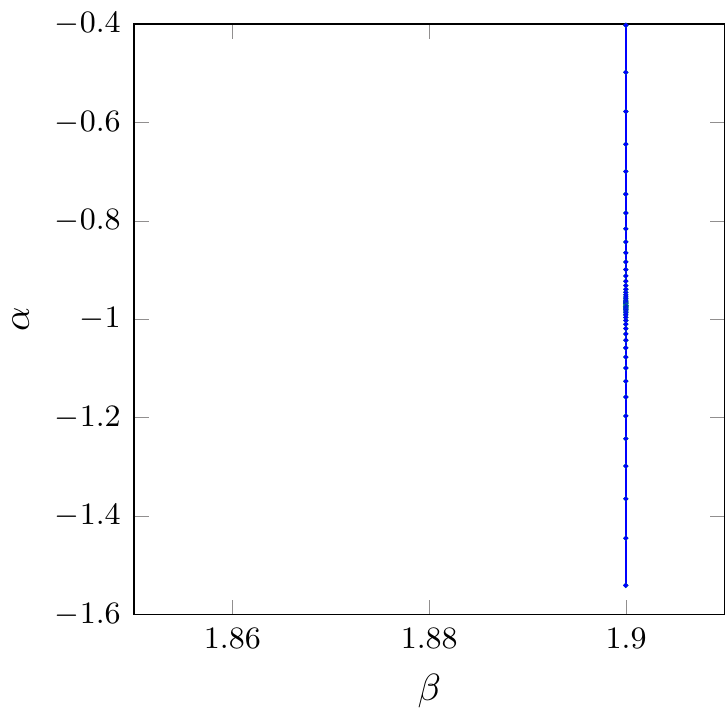}
		\fi
		\label{fig:FHN_Hopf_curve}}
	\hfill
	\subfloat[Limit point of cycles curve]{
		\ifcompileimages
  \tikzsetnextfilename{FHN_LPC_curve}%
  \input{tikz/FHN_LPC_curve}%

		\else
		    \includegraphics{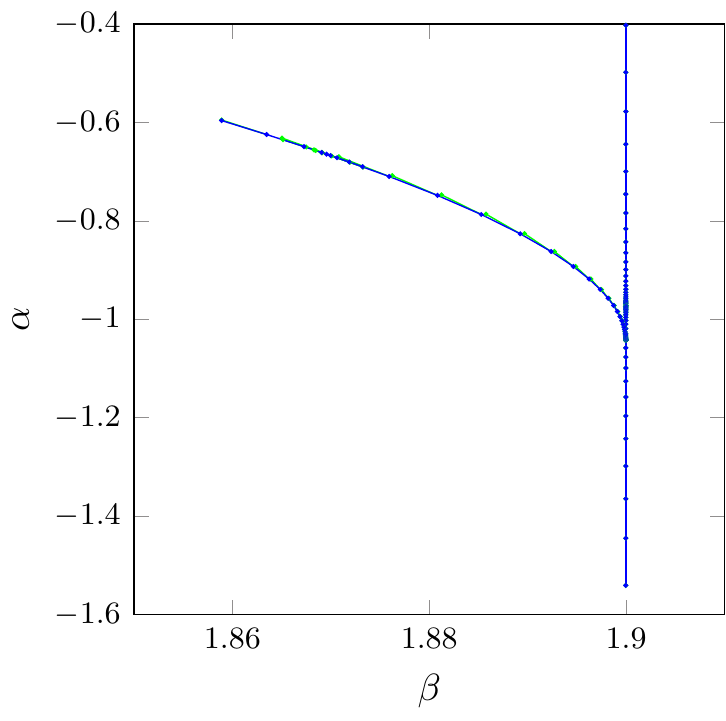}
		\fi
		\label{fig:FHN_LPC_curve}}
	\caption{\textup{(a)} Hopf curve continued from the manually constructed point. \textup{(b)} LPC curve continued from the detected generalized Hopf point using our predictors.}
\end{figure}

Since the simulations in \cite{Ma2011} do indicate a generalized-Hopf point nearby, we continue the Hopf point.
\begin{lstlisting}[style=customMatlab,escapechar=!]
%% Initialize Hopf branch
unfolding_pars=[ind.beta, ind.alpha];
hbr=df_brnch(funcs,unfolding_pars,'hopf');
hbr.point=hopf;
hbr.point(2)=hopf;
hbr.point(2).parameter(ind.alpha)=...
hbr.point(2).parameter(ind.alpha)+0.001;
method=df_mthod(funcs,'hopf');
method.point.print_residual_info=1;
hbr.point(2)=p_correc(funcs,hbr.point(2),ind.beta,[],method.point);
!\matlabrule!
%% Continue Hopf branch
figure(1); clf;
hbr=br_contn(funcs,hbr,30);
hbr=br_rvers(hbr);
hbr=br_contn(funcs,hbr,30);
title('Continued Hopf branch');
xlabel('$\beta$','Interpreter','LaTex')
ylabel('$\alpha$','Interpreter','LaTex')
box on
\end{lstlisting}
The continued branch \lstinline{hbr} is shown in \cref{fig:FHN_Hopf_curve}.

\subsection{Detect bifurcation points}
\begin{figure}
\centering
\ifcompileimages
  \tikzsetnextfilename{FHN_testfunction_genh}%
  \input{tikz/FHN_testfunction_genh}%

\else
\includegraphics{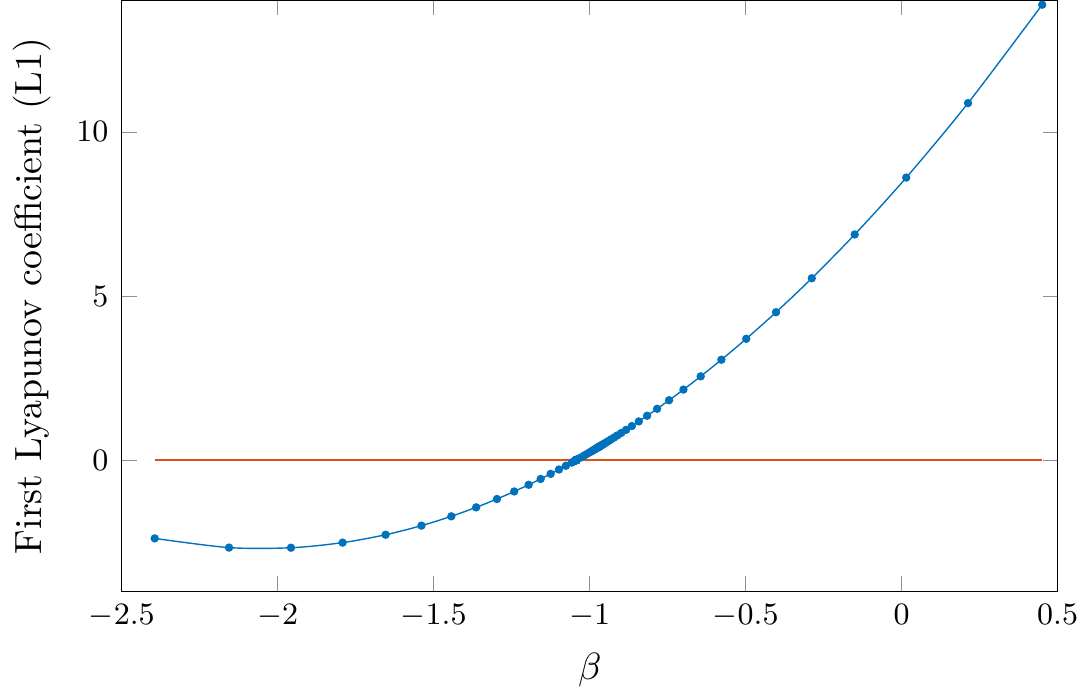}
\fi
\label{fig:FHN_testfunction}
\caption{Plot of the test function for a generalized Hopf point.}
\end{figure}

To detect bifurcation points on the Hopf branch, we use the function \lstinline|LocateSpecialPoints|.
\begin{lstlisting}[style=customMatlab]
[hbr_wbifs,hopftests,hc2_indices,hc2_types]=...
    LocateSpecialPoints(funcs,hbr);
\end{lstlisting}
The \MATLAB console shows the following output.
\begin{lstlisting}[style=matlabConsole]
HopfCodimension2: calculate stability if not yet present
HopfCodimension2: calculate L1 coefficients
HopfCodimension2: (provisional) 1 gen. Hopf  detected.
br_insert: detected 1 of 1: genh. Normalform:
    L2: -15.6733
    L1: -1.6801e-12
\end{lstlisting}
Thus a generalized Hopf point is indeed present on the Hopf branch \lstinline|hbr|. The returned branch \lstinline|hbr_wbifs| contains this point. The array \lstinline{hc2_indices} is used to subtract the generalized Hopf point below. If there would be more bifurcation points detected, \lstinline|hc2_types| can be used to inspect their types. Lastly, \lstinline|hopftests| stores the test functions to detect a bifurcation point. A change in sign in one of these functions indicates a bifurcation. The code below plots the test function for the generalized Hopf point, i.e. the first Lyapunov coefficient (L1), see \cref{fig:FHN_testfunction}.
\begin{lstlisting}[style=customMatlab]
al=arrayfun(@(x)x.parameter(ind.alpha),hbr_wbifs.point);
figure(2); clf;
plot(al,hopftests.genh(1,:),'.-',al,zeros(size(al)));
xlabel('$\beta$','Interpreter','LaTex');
ylabel('First Lyapunov coefficient (L1)')
title('Criticality along Hopf bifurcation curve')
\end{lstlisting}

%

%
\subsection{Continue limit point of cycle curve}
First, we subtract the detected generalized-Hopf point from the branch.
\begin{lstlisting}[style=customMatlab]
%% Calculate parameter-dependent normal form coefficients
% Select the located generalized Hopf point on the hopf_br_wbifs.
% Then convert the Hopf point to a generalized Hopf point.
% As before, we use the function nmfm_genh to calculate the normal 
% form coefficients. By adding the option free_pars and providing  
% the unfolding parameters the parameter-dependent normal  
% form coefficients are calculated.
hopf=hbr_wbifs.point(hc2_indices); % select generelized Hopf point
genh=p_togenh(hopf);
genh=nmfm_genh(funcs,genh,'free_pars',unfolding_pars);
genh.nmfm.L1
\end{lstlisting}
By inspecting the \lstinline|genh| structure, we obtain the correct parameter values of the generalized Hopf bifurcation.
\begin{lstlisting}[style=matlabConsole]
>> genh.parameter

ans =

    1.9000   -1.0429    1.7722

>> 
\end{lstlisting}
To continue the limit point of cycles curve emanating from the generalized Hopf point, we use the function \lstinline|C1branch_from_C2point|.
\begin{lstlisting}[style=customMatlab]
%% Continue LPC curve emanating from generalized-Hopf point
figure(1); [lpcfuncs,lpcbr,~]=C1branch_from_C2point(funcs,genh,...
    unfolding_pars,'codim2',genh.kind,'codim1','POfold');
nop=50; [lpcbr,suc]=br_contn(lpcfuncs,lpcbr,nop); assert(suc>0)
\end{lstlisting}
This function uses the file \lstinline|nmfm_POfold_from_genh_init.m| in which the predictor, see \cref{sec:GH_predictors,sec:genh_predictors}, is implemented. Using this function, a branch with three initial corrected cycles is created which is continued in the standard way, see \cref{fig:FHN_LPC_curve}.

\subsection{Calculate predicted periodic orbits}
To compare the computed parameter values and periodic orbits on the branch \lstinline|lpcbr| with the predictor, we again use the function \lstinline|C1branch_from_C2point|, but with the additional argument \lstinline|predictor| set to \lstinline|1| and \lstinline|step| to an interval of $\varepsilon$-values. Now the cycles are left uncorrected.
\begin{lstlisting}[style=customMatlab]
%% Predictor LPC curve emanating from generalized-Hopf point
[~,lpcbr_pred,~]=C1branch_from_C2point(funcs,genh,...
    unfolding_pars,'codim2',genh.kind,'codim1','POfold',...
    'step',linspace(0,1,45),'predictor',1);
\end{lstlisting}

\subsection{Bifurcation diagram}
\begin{figure}[ht]
\centering
\ifcompileimages
  \tikzsetnextfilename{FHN_bifdia}%
  \input{tikz/FHN_bifdia}%

\else
\includegraphics{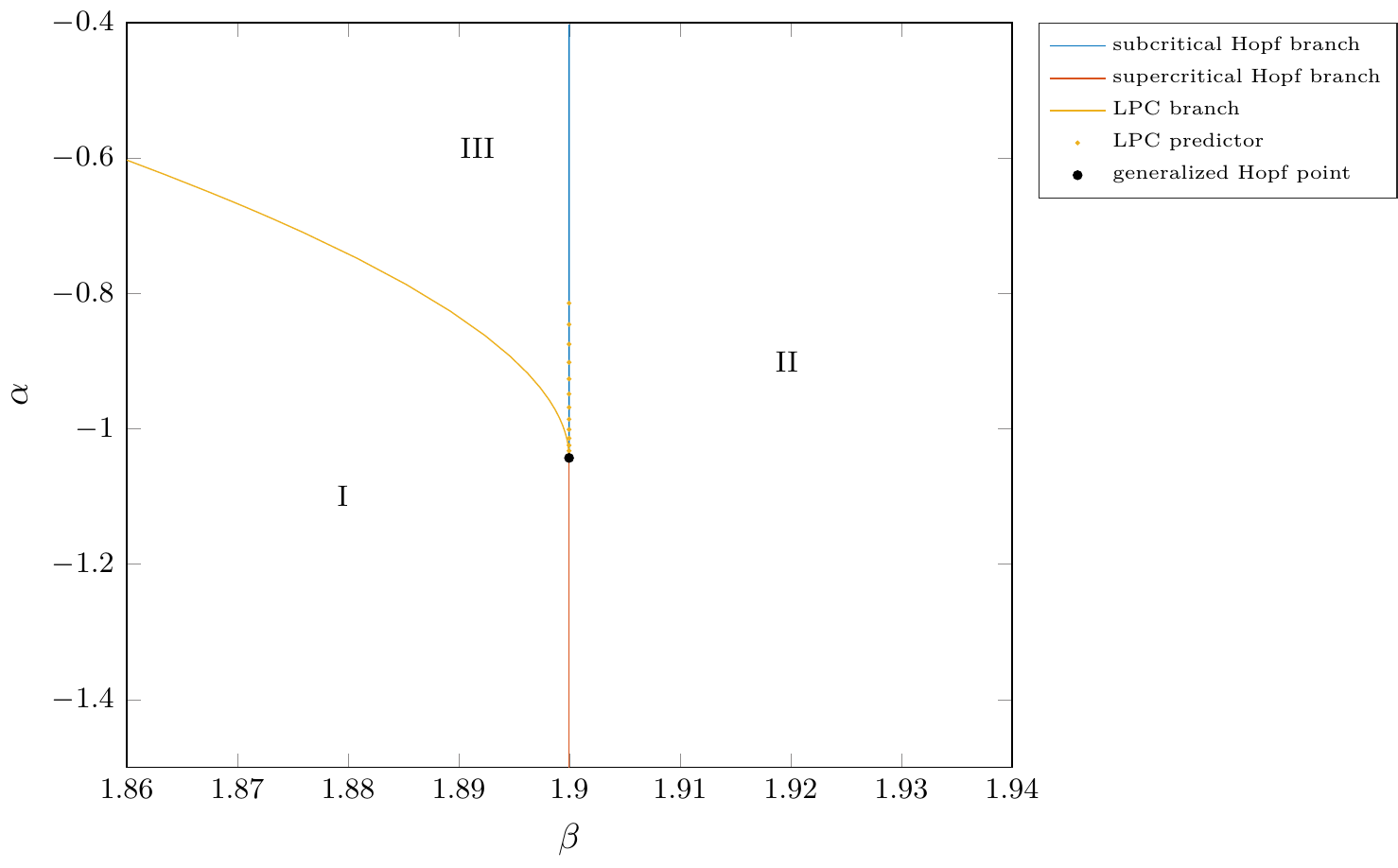}
\fi
\caption{Bifurcation diagram near the generalized
Hopf point in the system \cref{sm:eq:DDE_FHN} with unfolding parameters
$(\beta,\alpha)$. The bifurcation curves are nearly identical to
those in the bifurcation diagram of the topological normal form as
presented in \textup{\cite[page 314]{Kuznetsov2004}}.}
\label{sm:fig:FHN-bifurcation-diagram}
\end{figure}

The following code produces the bifurcation diagram presented in the main text, see \cref{fig:FHN-bifurcation-diagram}, and has been reproduced here in \cref{sm:fig:FHN-bifurcation-diagram}. The figure was exported with the \MATLAB and \OCTAVE compatible package {\tt matlab2tikz}, see \cite{matlab2tikz}. 
\begin{lstlisting}[style=customMatlab,caption=\MATLAB code for bifurcation diagram, label={sm:lst:bifurcation_diagram}]
%% Bifurcation diagram
figure(3); clf; hold on;
% Inline function to subtract parameters
getpars=@(points,ind) arrayfun(@(p)p.parameter(ind),points);
cm=colormap('lines');
L1s=hopftests.genh(1,:); % L1 along the hopf branch
% Plot sub- and supercritical Hopf branches
plot(getpars(hbr_wbifs.point(L1s>0),ind.beta),...
  getpars(hbr_wbifs.point(L1s>0),ind.alpha),'Color',cm(1,:),...
  'DisplayName','subcritical Hopf branch');
plot(getpars(hbr_wbifs.point(L1s<0),ind.beta),...
  getpars(hbr_wbifs.point(L1s<0),ind.alpha),'Color',cm(2,:),...
  'DisplayName','supercritical Hopf branch');
plot(getpars(lpcbr.point,ind.beta),...
  getpars(lpcbr.point,ind.alpha),'Color',cm(3,:),...
  'DisplayName','LPC branch');
plot(getpars(lpcbr_pred.point,ind.beta),...
  getpars(lpcbr_pred.point,ind.alpha),'.','Color',cm(3,:),...
  'DisplayName','LPC predictor');
plot(getpars(genh,ind.beta),getpars(genh,ind.alpha),'k.',...
  'MarkerSize',8,'DisplayName','generalized Hopf point');
title('Bifurcation diagram near generalized Hopf point')
xlabel('$\beta$','Interpreter','LaTex');
ylabel('$\alpha$','Interpreter','LaTex');
text(1.8779,-1.1001,'I');
text(1.9130,-0.9008,'II');
text(1.8890,-0.5854,'III');
axis([1.86 1.945 -1.5000 -0.4000])
legend('Location','NorthEast'); box on 
\end{lstlisting}

\subsection{Plot comparing computed and predicted periodic orbits}
\begin{figure}
	\centering
	\ifcompileimages
  \tikzsetnextfilename{FHN_compare_orbits}%
  \input{tikz/FHN_compare_orbits}%

	\else
	\includegraphics{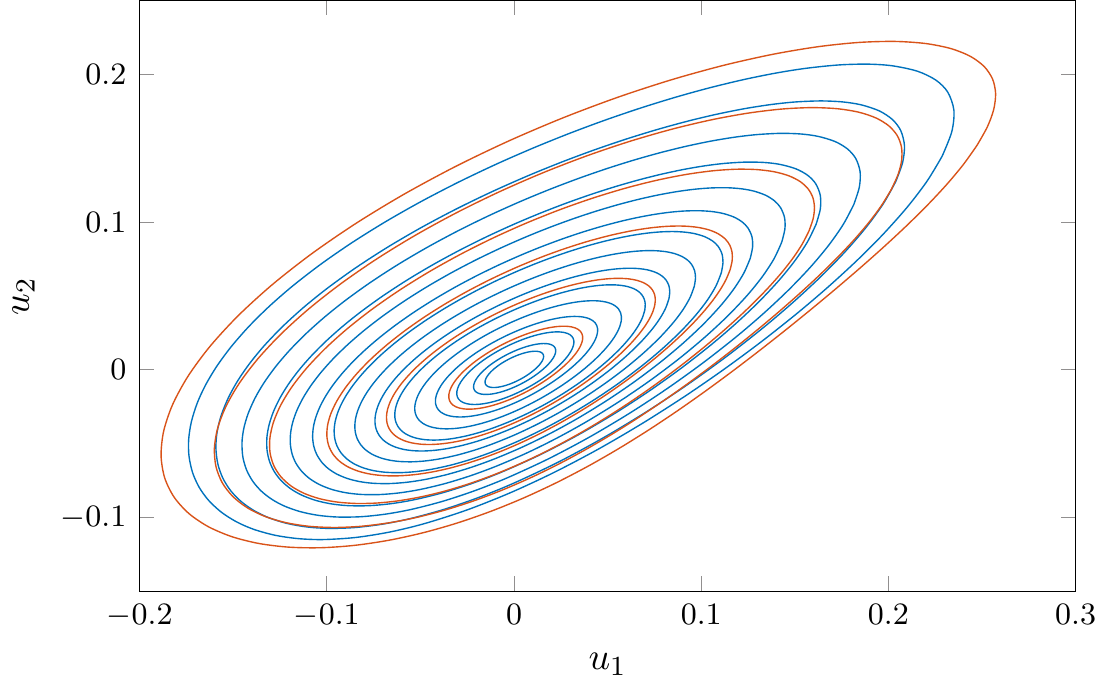}
	\fi
	\caption{Comparison between computed periodic orbits (blue) and predicted periodic orbits (red) emanating from the generalized Hopf bifurcation.}
	\label{sm:fig:FHN:compare_orbits}
\end{figure}

Lastly, we create a plot to compare the computed and predicted periodic orbits.
\begin{lstlisting}[style=customMatlab]
%% Plot comparing computed and predicted periodic orbits
figure(4); clf; hold on;
for i=1:15
    plot(lpcbr.point(i).profile(1,:),...
        lpcbr.point(i).profile(2,:),'Color',cm(1,:));
end
% Plot predicted periodic orbits
for i=1:7
    plot(lpcbr_pred.point(i).profile(1,:),...
        lpcbr_pred.point(i).profile(2,:),'Color',cm(2,:));
end
xlabel('$u_1$','Interpreter','LaTex');
ylabel('$u_2$','Interpreter','LaTex');
title('Compare computed and predicted periodic orbits')
box on
\end{lstlisting}
The resulting plot is shown in \cref{sm:fig:FHN:compare_orbits}. Note that the cycles shown have \emph{different} underlying parameter values. Nonetheless, we see that the cycles are in good agreement.

\subsection{Simulation with \MATLAB}
Next, we simulate the dynamics near the generalized Hopf point. For this, we take a point in each of the three regions as shown in \cref{sm:fig:FHN-bifurcation-diagram}. The following code, from the file \lstinline|FHN_simulation.m|, uses the \MATLAB function \lstinline|dde23|.
\begin{lstlisting}[style=customMatlab,escapechar=!]
%% Clean workspace and add DDE-BifTool scripts
%  to the MATLAB search path.
clear;      % clear variables
close all;  % close figures
ddebiftoolpath='../../../../';
addpath(strcat(ddebiftoolpath,'ddebiftool'),...
    strcat(ddebiftoolpath,'ddebiftool_extra_psol'),...
    strcat(ddebiftoolpath,'ddebiftool_extra_nmfm'),...
    strcat(ddebiftoolpath,'ddebiftool_utilities'));
load('FHN_results.mat')
!\matlabrule!
%% Point in region I
beta0=1.8779;
alpha0=-1.1001;
% Point near the steady-state
x1=0;
x2=0.01;
% Integrate
tfinal=1000;
sol = dde23(@(t,y,Z) funcs.sys_rhs([y,Z],...
	[beta0 alpha0 tau0]),tau0,[x1 x2],[0 tfinal]);
t=linspace(0,tfinal,1000);
y=deval(sol,t);
% Plot
title('Point in region I')
figure(1);clf;
xlabel('$u_1$','Interpreter','LaTex')
ylabel('$u_2$','Interpreter','LaTex')
plot(y(1,:),y(2,:))
!\matlabrule!
%% Point in region II
beta0=1.9130;
alpha0=-0.9008;
% Point near the steady-state
x2=0.1;
% Integrate
sol = dde23(@(t,y,Z) funcs.sys_rhs([y,Z],[beta0 alpha0 tau0]),...
tau0,[x1 x2],[0 tfinal]);
t=linspace(0,tfinal,1000);
y=deval(sol,t);
% Plot
figure(2);clf;
title('Point in region II')
xlabel('$u_1$','Interpreter','LaTex')
ylabel('$u_2$','Interpreter','LaTex')
plot(y(1,:),y(2,:))
!\matlabrule!
%% Point in region III
beta0=1.8890;
alpha0=-0.6081;
% Orbit converging to periodic orbit
x2=0.1;
% Integrate
tfinal=3000; % use lager time interval
sol = dde23(@(t,y,Z) funcs.sys_rhs([y,Z],[beta0 alpha0 tau0]),...
tau0,[x1 x2],[0 tfinal]);
t=linspace(0,tfinal,4000);
y1=deval(sol,t);
% Plot
figure(3);clf;
title('Point in region III')
xlabel('$u_1$','Interpreter','LaTex')
ylabel('$u_2$','Interpreter','LaTex')
plot(y1(1,:),y1(2,:))
% Orbit converging to the stable steady-state
% Integrate
x2=0.093;
sol = dde23(@(t,y,Z) funcs.sys_rhs([y,Z],[beta0 alpha0 tau0]),...
tau0,[x1 x2],[0 tfinal]);
y2=deval(sol,t);
% Add to plot
hold on; plot(y2(1,:),y2(2,:),'Color',cm(2,:));
!\matlabrule!
%% Time series of the previous solutions in region III
figure(4);clf;
plot(t,y1(1,:),t,y2(1,:))
title('Time series of solutions in region III')
xlabel('$t$','Interpreter','LaTex')
ylabel('$u_1$','Interpreter','LaTex')
\end{lstlisting}
In \cref{fig:FHN_sim1,fig:FHN_sim2,fig:FHN_sim3,fig:FHN_sim3_time_series} the resulting plots are shown, confirming the dynamics near the generalized Hopf point as predicted in \cite{Kuznetsov2004}.
\begin{figure}[ht!]
\centering
\ifcompileimages
\subfloat[\label{fig:FHN_sim1}]{
  \tikzsetnextfilename{FHN_sim1}%
  \input{tikz/FHN_sim1}%

}\hspace*{\fill}\subfloat[\label{fig:FHN_sim2}]{
  \tikzsetnextfilename{FHN_sim2}%
  \input{tikz/FHN_sim2}%

}\\
\subfloat[\label{fig:FHN_sim3}]{
  \tikzsetnextfilename{FHN_sim3}%
  \input{tikz/FHN_sim3}%

}\hspace*{\fill}\subfloat[\label{fig:FHN_sim3_time_series}]{
  \tikzsetnextfilename{FHN_sim3_time_series}%
  \input{tikz/FHN_sim3_time_series}%

}
\else
\subfloat[\label{fig:FHN_sim1}]{
\includegraphics{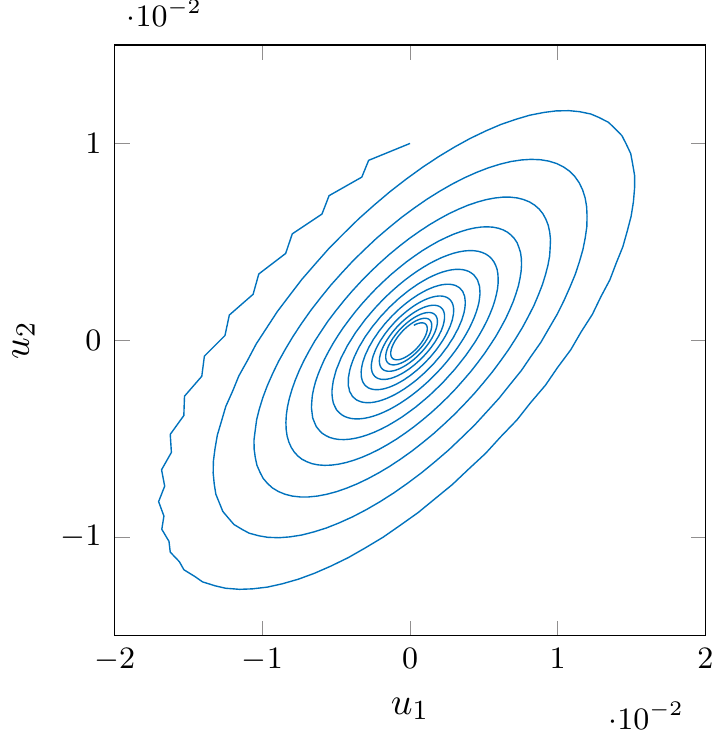}
}\hspace*{\fill}\subfloat[\label{fig:FHN_sim2}]{
\includegraphics{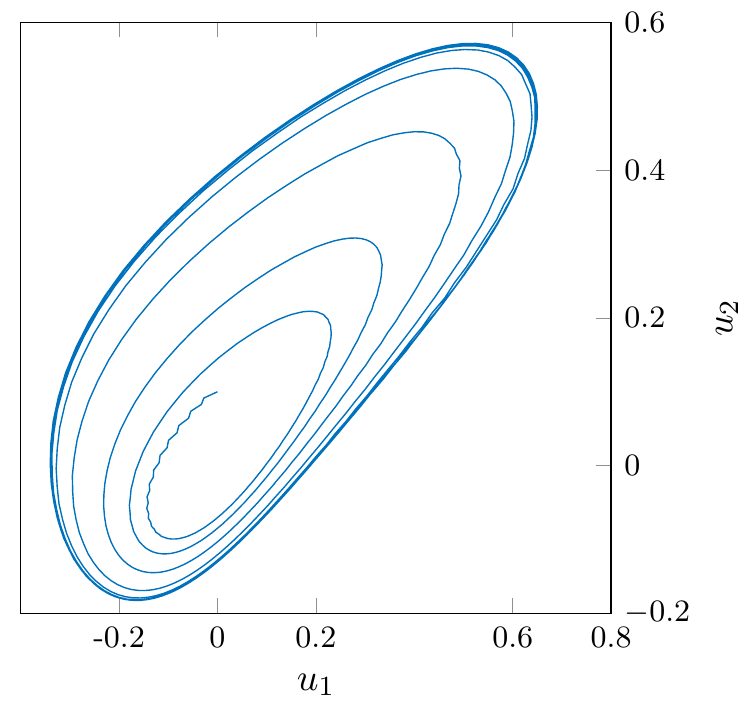}
}\\
\subfloat[\label{fig:FHN_sim3}]{
\includegraphics{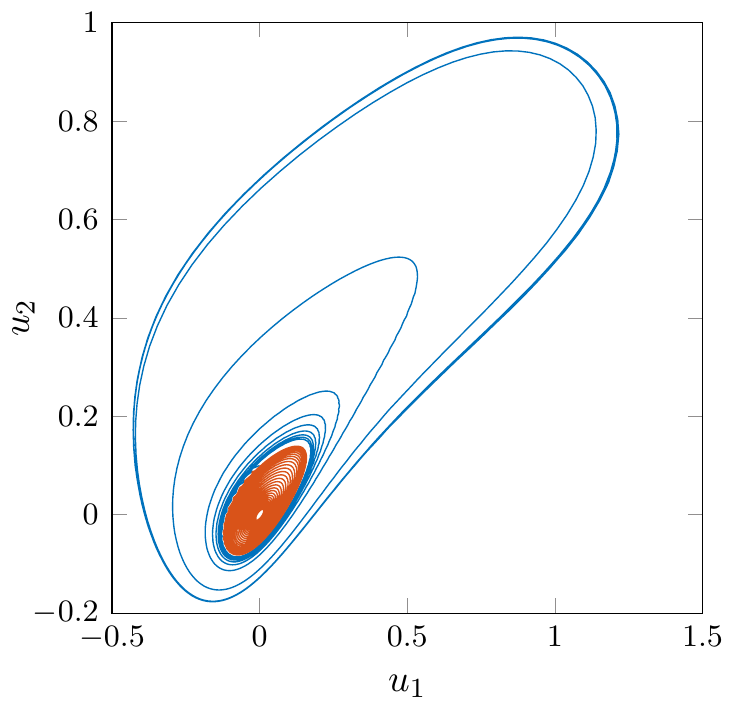}
}\hspace*{\fill}\subfloat[\label{fig:FHN_sim3_time_series}]{
\includegraphics{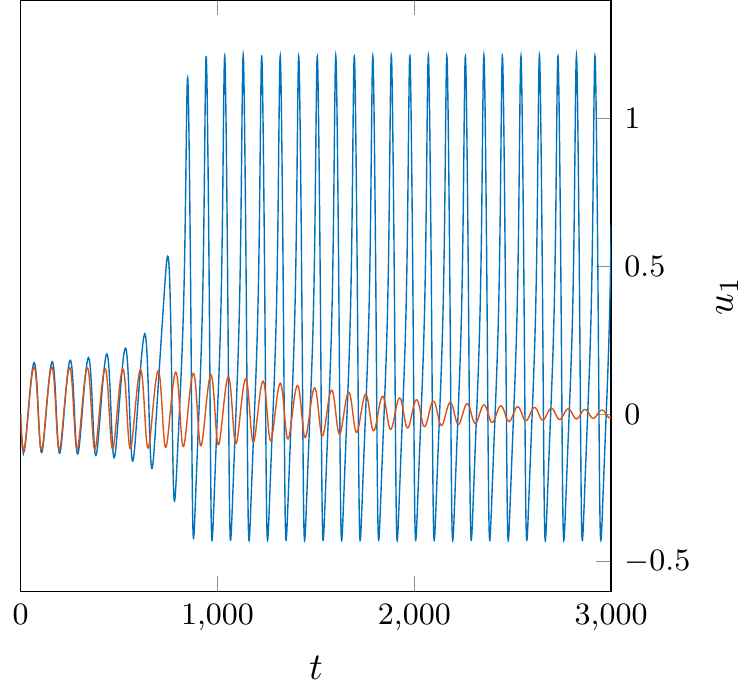}
}
\fi
\caption{Simulation near the generalized Hopf point in the system \cref{sm:eq:DDE_FHN}. In \textup{(a)} we see a stable steady-state corresponding to a point in region I. When we enter region II, the stability of the steady-state is lost and a stable cycle appears, as seen in  \textup{(b)}. In region III, there is a stable steady-state inside a stable cycle. This is confirmed in \textup{(c)} and \textup{(d)}. In \textup{(c)}, the initial point of the orbit in blue is just outside the unstable cycle and converges to the stable cycles. The initial point of the orbit in red is just inside the unstable cycle and converges to the stable steady-state. In \textup{(d)}, the time series of these orbits are shown in the $(t,u_1(t))$ plane.}
\end{figure}

\section{Fold-Hopf bifurcation in the Rose\textendash Hindmarsh model with time delay}
In \cite{Ma2011} a Rose-Hindmarsh model \cite{Hindmarsh1982,Hindmarsh1984} with time delay in the self-feedback process, 
\begin{equation}
\begin{cases}
\begin{aligned}
\dot{x}(t)& = y(t)-ax^3(t)+bx^2(t-\tau)-cz(t)+I_{app},\\
\dot{y}(t)& = c-dx^2(t)-y(t),\\
\dot{z}(t)& = r(S(x(t)-\chi)-z(t)),
\end{aligned}
\end{cases}\label{sm:eq:Rose-Hindmarsh}
\end{equation}
is considered, see \cref{sec:ex_Rose_Hindmarsh}. The parameters values
\begin{equation}
a=1.0,\qquad b=3.0,\qquad c=1.0,\qquad d=5.0,\qquad\chi=-1.6,\qquad r=0.001\label{sm:eq:rose_hindmarsh_pm1}
\end{equation}
are fixed and $(I_{app},S)$ are the unfolding parameters.

\begin{remark}
	This demonstration can be found in the directory \lstinline|demos/tutorial/VII/RH| relative to the main directory of the \DDEBIFTOOL package. Here, we omit the code to generate a system file. The system file \lstinline|sym_RH_mf.m| has been generated with the script \lstinline|gen_sym_RS.m|. Also, we assume that the \DDEBIFTOOL package has been loaded as in \cref{sm:lst:searchpath}. The code in \crefrange{sm:sec:RH:pars_and_funcs}{sm:sec:RH:comparing_period_orbits} highlights the important parts of the file \lstinline|RH.m|.
\end{remark}

\subsection{Set parameter names and funcs structure} \label{sm:sec:RH:pars_and_funcs}
As in the previous example, we set the parameter names and define the \lstinline|funcs| structure.
\begin{lstlisting}[style=customMatlab]
%% Set parameter names
parnames={'Iapp','S','r','tau'};
cind=[parnames;num2cell(1:length(parnames))];
ind=struct(cind{:});
%% Set funcs structure
% We load the precalculated multilinear forms. These have been
% generated with the file gen_sym_RH.m.
funcs=set_symfuncs(@sym_RH_mf,'sys_tau',@()ind.tau);
\end{lstlisting}

\subsection{Stability and normal form coefficients of the fold-Hopf point} We construct a steady-state at the fold-Hopf point and calculate its stability.
\begin{lstlisting}[style=customMatlab]
%% Construct fold-Hopf point
a=1.0; b=3.0; c=1.0; d=5.0; chi=-1.6;
r=1.0e-03;
S=-0.57452592;
[xstar,Iapp,tau]=bifurcationvalues(a,b,c,d,chi,r,S);
% Construct steady-state point
stst=dde_stst_create('x',[xstar; c-d*xstar^2; S*(xstar-chi)]);
stst.parameter([ind.Iapp ind.S ind.r ind.tau])=[Iapp S r tau];
% Calculate stability
method=df_mthod(funcs,'stst');
stst.stability=p_stabil(funcs,stst,method.stability);
stst.stability.l1(1:5)
\end{lstlisting}
The function \lstinline|bifurcationvalues| calculates $(x_\star,I_{app},\tau)$ according to the formulas as given in \cref{sec:ex_Rose_Hindmarsh}. The \MATLAB console shows the following output.
\begin{lstlisting}[style=matlabConsole]
ans =

    0.0000 + 0.0000i
   -0.0000 + 1.0079i
   -0.0000 - 1.0079i
   -0.0994 + 1.9324i
   -0.0994 - 1.9324i
\end{lstlisting}
We have a zero eigenvalue and a pair of purely imaginary eigenvalues. Furthermore, the remaining eigenvalues have negative real parts. Next, we calculate the normal form coefficients and the transformation to the center manifold with the function \lstinline|nmfm_zeho|, which implements the coefficients as derived in \cref{sec:fold-Hopf}. For this we need to set the argument \lstinline|free_pars| to the unfolding parameter $(I_{app},S)$. These coefficients will be used to start the continuation of the various branches emanating from the fold-Hopf point.
\begin{lstlisting}[style=customMatlab]
% Calculate normal coefficients
hopf=p_tohopf(funcs,stst);
zeho=p_tozeho(hopf);
unfolding_pars=[ind.Iapp, ind.S];
zeho=nmfm_zeho(funcs,zeho,'free_pars',unfolding_pars);
zeho.nmfm
\end{lstlisting}
The \MATLAB console shows the following output.
\begin{lstlisting}[style=matlabConsole,keepspaces=true]
ans = 

  ?struct? with fields:

             g200: -0.0024
             g110: 0.3296 + 0.7006i
             g011: -0.0078
             g300: 0.0106
             g111: -0.0146
             g210: -2.7764 - 3.1806i
             g021: -0.3745 - 2.4754i
                b: -0.0024
                c: -0.0078
                d: 0.3296 + 5.2274i
                e: 15.6941
                s: 1.8487e-05
            theta: -139.0315
    transcritical: 0
             h200: [1x1 struct]
             h011: [1x1 struct]
             h020: [1x1 struct]
             h110: [1x1 struct]
                K: [2x2 double]
           h000mu: [1x2 struct]
           omega1: 7.4540
           omega2: 2.1259
\end{lstlisting}
Since $s > 0$ and $\theta(0) < 0$ global bifurcations or invariant tori are present. However, since the sign of $e$ is 
positive the stability of the invariant tori will be unstable for nearby parameter values. It follows that the torus observed in \cite{Ma2011} by simulations does not originate from the fold-Hopf point under consideration.

\subsection{Set bifurcation parameter range and step size bounds}
Before continuing the various branches emanating from the fold-Hopf point, we create the variable \lstinline|brpars| containing parameter bounds and maximal stepsizes.
\begin{lstlisting}[style=customMatlab]
%% Set bifurcation parameter range and step size bounds
brpars={'min_bound',[ind.Iapp -20; ind.S -12],...
        'max_bound',[ind.Iapp  10; ind.S  5],...
        'max_step', [ind.Iapp 4e-02; ind.S 4e-02]};
\end{lstlisting}
\subsection{Continue NS, Hopf and fold branch} \label{sm:sec:RH:continuation}
\begin{figure}
\centering
\ifcompileimages
  \tikzsetnextfilename{RH_bifurcation_diagram_I}%
  \input{tikz/RH_bifurcation_diagram_I}%

\else
\includegraphics{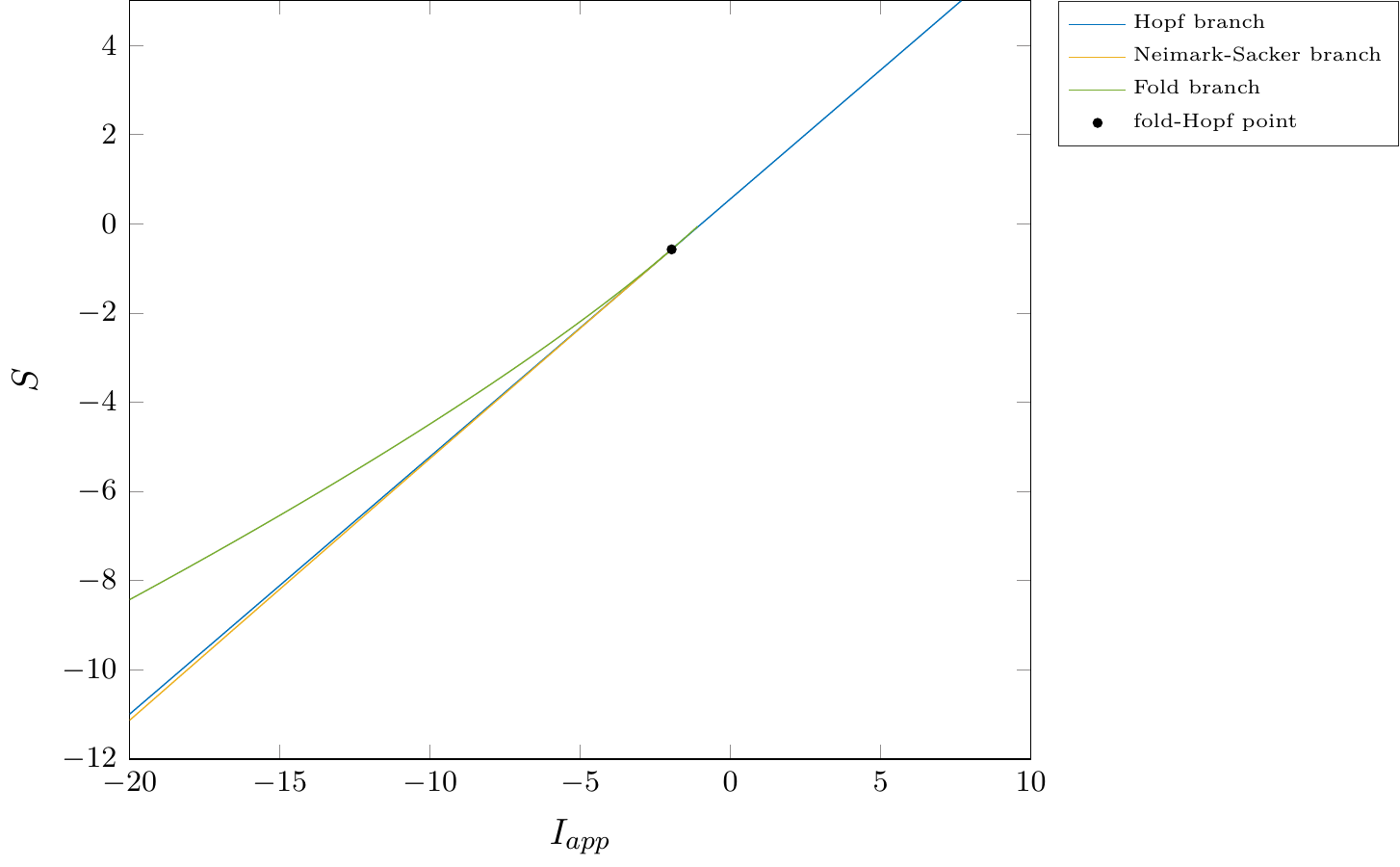}
\fi
\caption{Bifurcation diagram near the fold-Hopf point in \cref{sm:eq:Rose-Hindmarsh} with $(r,S)=(0.001,-0.57452592)$.}
\label{sm:fig:RH_bifurcation_diagram_I}
\end{figure}

As in the previous example, we use the function \lstinline|C1branch_from_C2point| to start to continue the branches emanating from the fold-Hopf point. \cref{sm:fig:RH_bifurcation_diagram_I} is created using similar code as in Listing \ref{sm:lst:bifurcation_diagram}. We remark that even when there would be stable tori present for nearby parameter values, the window in which these tori would exist is quite small. Indeed, the parameter values would have to be below the Hopf curve and above the Neimark-Sacker curve to the left of the fold-Hopf point.
\begin{lstlisting}[style=customMatlab,escapechar=!]
%% Continue Neimark-Sacker curve emanating from fold-Hopf point
[trfuncs,nsbr,~]=C1branch_from_C2point(funcs,zeho,unfolding_pars,...
    'codim2',zeho.kind,'codim1','TorusBifurcation',...
    brpars{:},'step',1e-03,'plot',0);
ntrsteps=1000; [nsbr,suc]=br_contn(trfuncs,nsbr,ntrsteps);
!\matlabrule!
%% Continue Hopf curve emanating from fold-Hopf point
[~,hbr,suc]=C1branch_from_C2point(funcs,zeho,unfolding_pars,...
    'codim2',zeho.kind,'codim1','hopf',...
    brpars{:},'step',1e-03,'plot',0);
nop=1000; [hbr,suc]=br_contn(funcs,hbr,nop); assert(suc>0)
hbr=br_rvers(hbr);
[hbr,suc]=br_contn(funcs,hbr,nop);
!\matlabrule!
%% Continue fold curve emanating from fold-Hopf point
[~,fbr,suc]=C1branch_from_C2point(funcs,zeho,unfolding_pars,...
    'codim2',zeho.kind,'codim1','fold','step',1e-03,'plot',0);
nop=1000; [fbr,suc]=br_contn(funcs,fbr,nop);
fbr=br_rvers(fbr);
[fbr,suc]=br_contn(funcs,fbr,nop);
\end{lstlisting}

\subsection{Stable invariant tori}
We change the parameters $r = 1.4$ and $S = -8$, while keeping the other fixed parameters as in \cref{sm:eq:rose_hindmarsh_pm1}. Using the formulas given in \cref{sec:ex_Rose_Hindmarsh}, we calculate $x_\star,I_{app}$ and $\tau$.
\begin{lstlisting}[style=customMatlab]
%% Different parameters with stable torus
r=1.4; S=-8;
[xstar,Iapp,tau]=bifurcationvalues(a,b,c,d,chi,r,S);
% Construct steady-state point
stst=dde_stst_create('x',[xstar; c-d*xstar^2; S*(xstar-chi)]);
stst.parameter([ind.Iapp ind.S ind.r ind.tau])=[Iapp S r tau];
% Calculate stability
method=df_mthod(funcs,'stst');
stst.stability=p_stabil(funcs,stst,method.stability);
stst.stability.l1(1:5)
\end{lstlisting}
The \MATLAB console outputs
\begin{lstlisting}[style=matlabConsole]
ans =

    0.00000 +  0.00000i
    0.00000 +  5.60424i
    0.00000 -  5.60424i
   -0.27870 +  0.00000i
   -0.66607 + 11.94839i
   -0.66607 - 11.94839i
\end{lstlisting}
which are indeed the eigenvalues that should be present at a fold-Hopf bifurcation. Next, we calculate the normal form coefficients.
\begin{lstlisting}[style=customMatlab]
% Calculate normal coefficients
hopf=p_tohopf(funcs,stst);
zeho=p_tozeho(hopf);
unfolding_pars=[ind.Iapp, ind.S];
zeho=nmfm_zeho(funcs,zeho,'free_pars',unfolding_pars);
\end{lstlisting}
Inspecting the normal form coefficients yields
\begin{lstlisting}[style=matlabConsole]
>> fprintf('s=%f, theta=%f, e=%f\n',...
       zeho.nmfm.s,zeho.nmfm.theta,zeho.nmfm.e)

s=1.770013, theta=-0.156886, e=-0.037794
>>
\end{lstlisting}
The coefficients $s$ and $\theta(0)$ reveal that we are in case III of the fold-Hopf bifurcation, see \cite[page 342]{Kuznetsov2004}. Since the sign of $e(0)$ is negative, there is a time reversal to take into account. Therefore, we expect stable tori to be present for nearby parameter values.

\subsection{Adjusting bifurcation parameter range}
We adjust the variable \lstinline|brpars| to reflect the current situation.
\begin{lstlisting}[style=customMatlab]
%% Set bifurcation parameter range and step size bounds
brpars={'min_bound',[ind.Iapp -20; ind.S -9],...
        'max_bound',[ind.Iapp -18; ind.S -7],...
        'max_step', [ind.Iapp 0.04; ind.S 0.04]};
\end{lstlisting}

\subsection{Detect special points on the Hopf branch}
Since the code to continue the Neimark-Sacker, Hopf and fold curves is identical to the code in \cref{sm:sec:RH:continuation}, we continue with detecting bifurcations on the Hopf branch. The Hopf points on the branch \lstinline|hbr_wbifs| will contain the normal form coefficients \lstinline|L1| and \lstinline|L2|. These will be used to visualize the criticality of the Hopf points (sub or super) in the bifurcation diagram.
\begin{lstlisting}[style=customMatlab]
[hbr_wbifs,hopftests,hc2_indices,hc2_types]=...
    LocateSpecialPoints(funcs,hbr);
\end{lstlisting}

\subsection{Predictors}
As in the previous example, we obtain predictors for the various branches simply by setting the argument \lstinline|predictor| to 1 and the argument \lstinline|step| to a range of $\varepsilon$-values when calling the function \lstinline|C1branch_from_C2point|.
\begin{lstlisting}[style=customMatlab]
%% Predictors for Neimark-Sacker and Hopf curves
[~,nsbr_pred]=C1branch_from_C2point(funcs,zeho,unfolding_pars,...
    'codim2','zeho','codim1','TorusBifurcation',...
    'step',linspace(1e-03,2.2e-01,40),'predictor',1);
[~,hbrsub_pred]=C1branch_from_C2point(funcs,zeho,unfolding_pars,...
    'codim2','zeho','codim1','hopf',...
    'step',linspace(0,1e-03,20),'predictor',1);
[~,hbrsup_pred]=C1branch_from_C2point(funcs,zeho,unfolding_pars,...
    'codim2','zeho','codim1','hopf',...
    'step',linspace(-1e-03,0,20),'predictor',1);
\end{lstlisting}

\subsection{Bifurcation diagram}
\begin{figure}
\centering
\ifcompileimages
  \tikzsetnextfilename{RH_bifurcation_diagram_II}%
  \input{tikz/RH_bifurcation_diagram_II}%

\else
\includegraphics{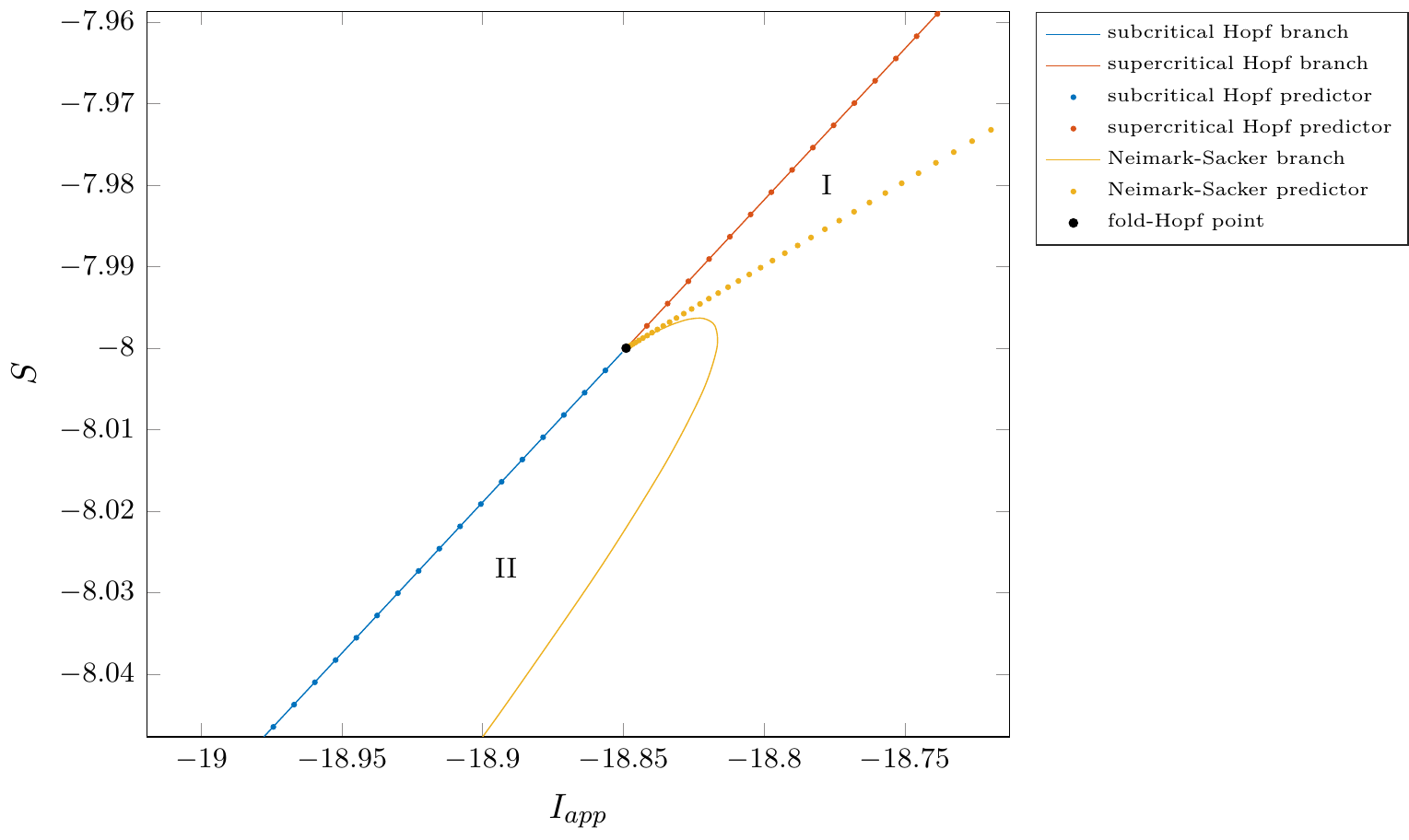}
\fi
\caption{Bifurcation diagram near the fold-Hopf point in \cref{sm:eq:Rose-Hindmarsh} with $(r,S)=(1.4,-8)$. The fold branch is not included here since it is indistinguishable from the Hopf curve at this scale.}
\label{sm:fig:RH_bifurcation_diagram_II}
\end{figure}

We plot the obtained curves and the predictor for the Neimark-Sacker and Hopf curve with the following code.
\begin{lstlisting}[style=customMatlab]
%% Plot comparing computed and predicted Neimark-Sacker curve
figure(8); clf; hold on;
nsbr2_pm_pred = [getpars(nsbr_pred,ind.Iapp); ...
    getpars(nsbr_pred,ind.S)];
hbrsub_pm_pred = [getpars(hbrsub_pred,ind.Iapp); ....
    getpars(hbrsub_pred,ind.S)];
hbrsup_pm_pred = [getpars(hbrsuper_pred,ind.Iapp); ...
    getpars(hbrsuper_pred,ind.S)];
plot(hbrsub_pm(1,:),hbrsub_pm(2,:),'Color',cm(1,:),...
    'DisplayName','subcritical Hopf branch');
plot(hbrsuper_pm(1,:),hbrsuper_pm(2,:),'Color',cm(2,:),...
    'DisplayName','supercritical Hopf branch');
plot(hbrsub_pm_pred(1,:),hbrsub_pm_pred(2,:),'.','Color',cm(1,:),...
    'DisplayName','subcritical Hopf predictor');
plot(hbrsup_pm_pred(1,:),hbrsup_pm_pred(2,:),'.','Color',cm(2,:),...
    'DisplayName','supercritical Hopf predictor');
plot(nsbr1_pm(1,:),nsbr1_pm(2,:),'Color',cm(3,:),...
    'DisplayName','Neimark-Sacker branch');
plot(nsbr2_pm_pred(1,:),nsbr2_pm_pred(2,:),'.','Color',cm(3,:),...
    'DisplayName','Neimark-Sacker predictor');
plot(zeho.parameter(ind.Iapp),zeho.parameter(ind.S),'k.',...
    'MarkerSize',12,'DisplayName','fold-Hopf point')
title('Neimark-Sacker curve emanating from the fold-Hopf point')
axis([-19.0193  -18.7128   -8.0477   -7.9587])
xlabel('$I_{app}$','Interpreter','LaTex');
ylabel('$S$','Interpreter','LaTex');
text(-18.784,-7.981,  'I','FontSize',14); % stable period orbit
text(-18.905,-8.032, 'II','FontSize',14); % stable 2d torus
legend('Location','NorthWest'); box on
\end{lstlisting}
\cref{sm:fig:RH_bifurcation_diagram_II} shows the resulting bifurcation diagram.

\subsection{Plots comparing computed and predicted periodic orbits} \label{sm:sec:RH:comparing_period_orbits}
We create a plot to compare the computed and predicted periodic orbits.
\begin{lstlisting}[style=customMatlab]
%% Plot comparing computed and predicted periodic orbits
figure(9); clf; hold on;
for i=1:14
  plot3(nsbr.point(i).profile(1,:),nsbr.point(i).profile(2,:),...
        nsbr.point(i).profile(3,:),'Color',cm(1,:));
end
for i=1:9
  plot3(nsbr_pred.point(i).profile(1,:),...
        nsbr_pred.point(i).profile(2,:),...
        nsbr_pred.point(i).profile(3,:),'Color',cm(2,:));
end
title('Comparison between computed and predicted periodic orbits')
xlabel('x'); ylabel('y'); zlabel('z')
view(3); grid on
\end{lstlisting}
The resulting plot is shown in \cref{sm:fig:RH:compare_orbits}.
To compare the computed and predicted periods, we use the following code.
\begin{lstlisting}[style=customMatlab]
%% Compare computed and predicted periods
figure(10); clf; hold on;
% Plot computed periods on nsbr(1) and nsbr(2)
omegas1=arrayfun(@(p)p.period,nsbr.point);
omegas1_pred=arrayfun(@(p)p.period,nsbr_pred.point);
plot(getpars(nsbr,ind.Iapp),omegas1,'Color',cm(1,:));
plot(getpars(nsbr_pred,ind.Iapp),omegas1_pred,'.','Color',cm(1,:));
title('Compare computed and predicted periods')
xlabel('$I_{app}$','Interpreter','LaTex')
ylabel('$\omega$','Interpreter','LaTex')
axis([-19.1016  -18.6500    1.1144    1.1220])
legend('Computed period','Predicted period','Location','SouthEast')
\end{lstlisting}
The resulting plot is shown in \cref{sm:fig:RH:compare_periods}.

\begin{figure}
\centering
\ifcompileimages
\subfloat[]{%
  \tikzsetnextfilename{RH_compare_orbits}%
  \input{tikz/RH_compare_orbits}%
\label{sm:fig:RH:compare_orbits}} \hfill
\subfloat[]{%
  \tikzsetnextfilename{RH_compare_periods}%
  \input{tikz/RH_compare_periods}%
\label{sm:fig:RH:compare_periods}}
\else
\subfloat[]{\includegraphics{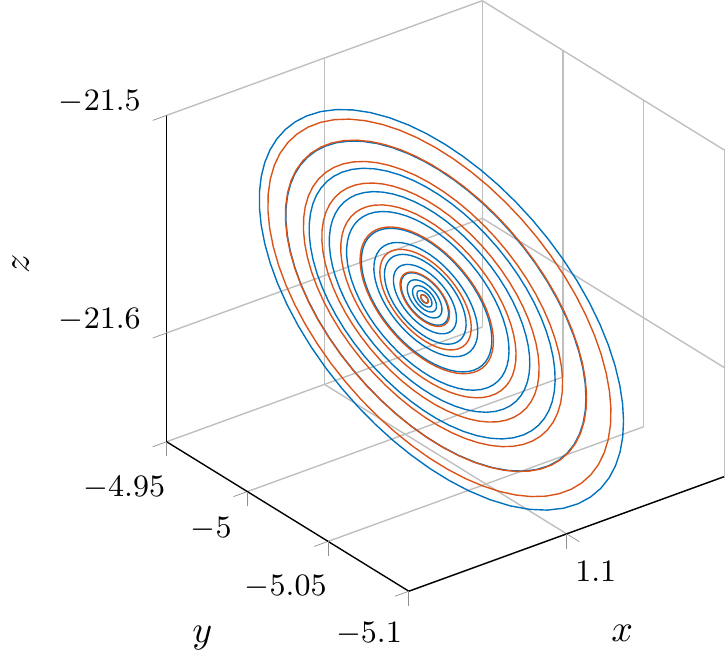}\label{sm:fig:RH:compare_orbits}} \hfill
\subfloat[]{\includegraphics{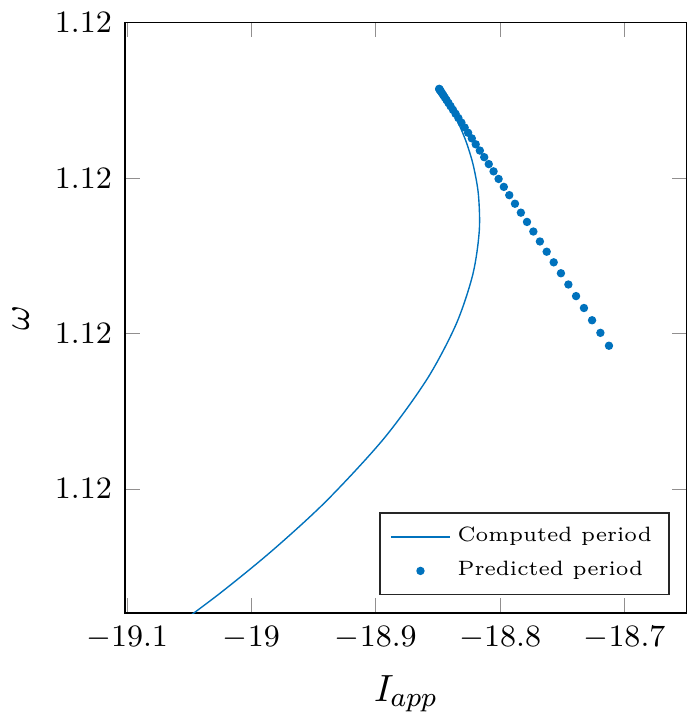}\label{sm:fig:RH:compare_periods}}
\fi
\caption{In \textup{(a)} the computed periodic orbits (blue) and predicted periodic orbits (red) are compared. In \textup{(b)} the computed predicted periods of the cycles are compared. We see that both are in good agreement.}
\end{figure}

\subsection{Simulation with \PYDELAY}
\begin{figure}[!th]
\centering
\ifcompileimages
  \tikzsetnextfilename{RH_simulation1}%
  \input{tikz/RH_simulation1}%
 \\
  \tikzsetnextfilename{RH_simulation2}%
  \input{tikz/RH_simulation2}%
 \\[0.6cm]
  \tikzsetnextfilename{RH_simulation3}%
  \input{tikz/RH_simulation3}%
 \\
  \tikzsetnextfilename{RH_simulation4}%
  \input{tikz/RH_simulation4}%

\else
\includegraphics{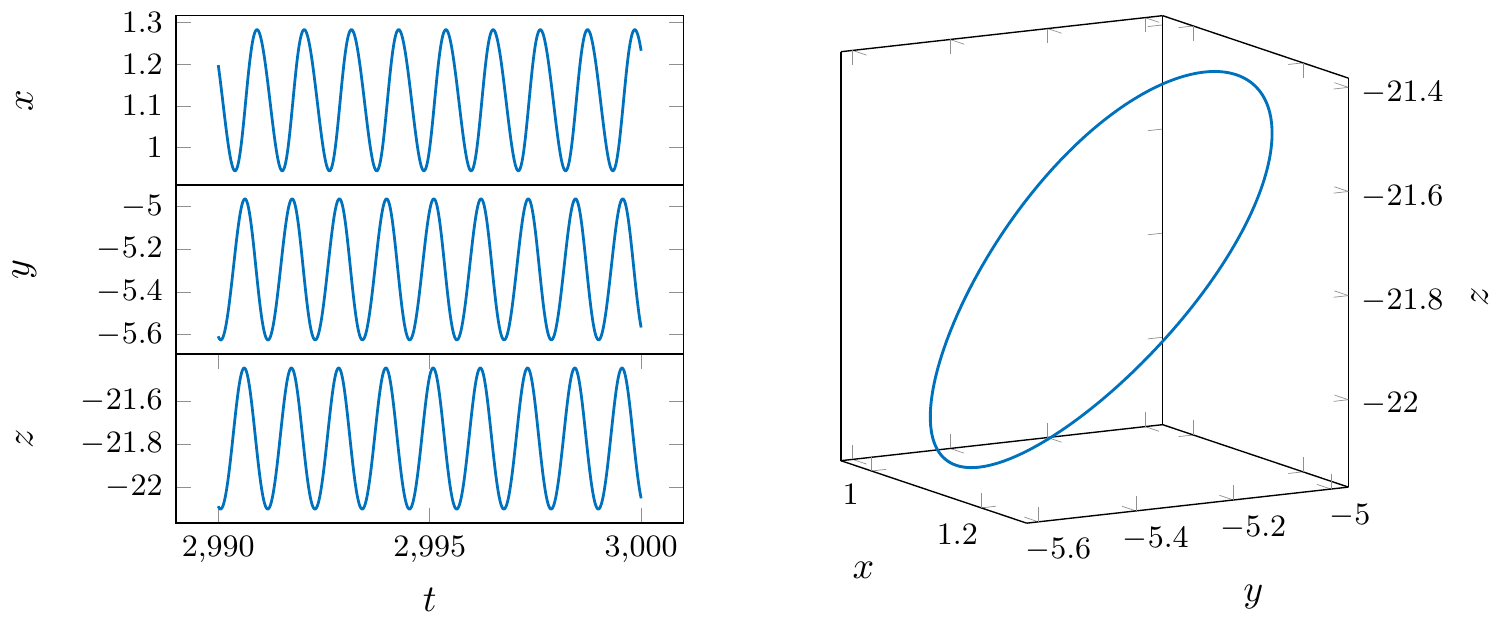}
\includegraphics{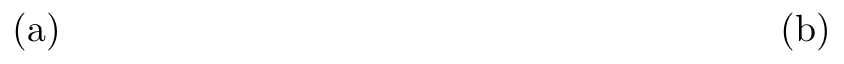}
\includegraphics{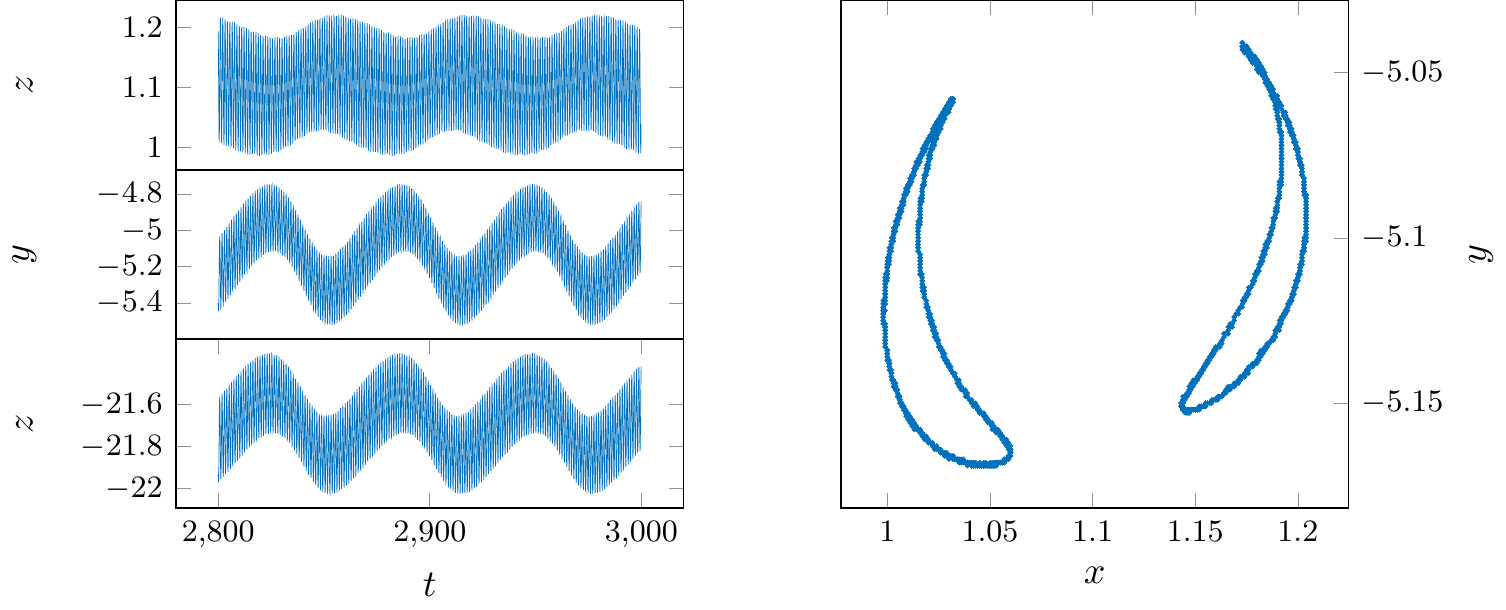}
\includegraphics{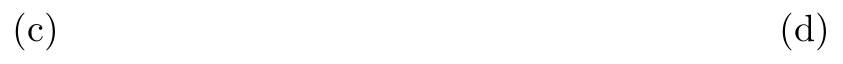}
\fi
\caption{Simulation near fold-Hopf point in \cref{eq:Rose-Hindmarsh}. In \textup{(a)} and \textup{(b)}, a stable periodic solution is shown. In \textup{(c)} and \textup{(d)}, a stable torus is shown. In \textup{(c)}, the time series is plotted while in \textup{(d)}, the cross-section defined by $z(t)=-21.75$ in the phase-space $(x,y,z)$ is taken.}
\label{fig:Rose-Hindmarsh-simulation}
\end{figure}

In this section we simulate the dynamics near the fold-Hopf point. The following code, from the file \lstinline|RH_simulation_torus.py|, uses the Python package \PYDELAY \cite{Flunkert2009Flunkert}.
\begin{lstlisting}[language=Python,escapechar=!]
import numpy as np
from pydelay import dde23
import matplotlib.pyplot as plt
from mpl_toolkits.mplot3d import Axes3D
import matplotlib.pyplot as plt
from matplotlib import colors

# Number of time units
tfinal_cycle = 10000
tfinal_torus = 20000

# Define DDE
eqns = {
    'x' : 'y-a*pow(x,3)+b*pow(x(t-tau),2)-c*z+Iapp',
    'y' : 'c-d*pow(x,2)-y',
    'z' : 'r*(S*(x-chi)-z)'
}

# Set parameters for torus
params_torus = {
    'a':1,'b':3,'c':1,'d':5.0,'chi':-1.6,'r':1.4,
    'tau':0.940246941050084,
    'Iapp':-18.902420391705071,
    'S':-8.045234985422740
}

# Set parameters for torus
params_cycle = params_torus.copy()
params_cycle['Iapp']=-18.886177304147466
params_cycle['S']=-8.044197084548104

# Set number of timesteps from the end to plot
timesteps_torus=300
timesteps_cycle=10

# Select periodic orbit or torus
tfinal,params,timesteps=tfinal_cycle,params_cycle,timesteps_cycle
#tfinal,params,timesteps=tfinal_torus,params_torus,timesteps_torus

# Solve DDE
dde = dde23(eqns=eqns, params=params)
dde.set_sim_params(tfinal=tfinal)
dde.set_sim_params(tfinal=tfinal, dtmax=0.001)
histfunc = {
    'x': lambda t: 1.097167540709727, 
    'y': lambda t: -5.018883061935152,
    'z': lambda t: -21.577340325677817
}
dde.hist_from_funcs(histfunc, 51)
dde.run()

# Subtract solution components
sol0 = dde.sample(tfinal-timesteps,tfinal, 0.001)
t = sol0['t']
x = sol0['x']
y = sol0['y']
z = sol0['z']

# Plot time series
fig = plt.figure()
plt.figure(1)
plt.subplot(311)
plt.xlabel('$t$')
plt.ylabel('$x(t)$')
plt.plot(t, x, c='royalblue')

plt.subplot(312)
plt.xlabel('$t$')
plt.ylabel('$y(t)$')
plt.plot(t, y, c='royalblue')

plt.subplot(313)
plt.xlabel('$t$')
plt.ylabel('$z(t)$')
plt.plot(t, z, c='royalblue')
fig.set_size_inches(18.5, 10.5)
plt.show()

# Plot the solution in phase-space
fig = plt.figure()
ax = fig.gca(projection='3d')
ax.w_xaxis.set_pane_color((1.0, 1.0, 1.0, 1.0))
ax.w_yaxis.set_pane_color((1.0, 1.0, 1.0, 1.0))
ax.w_zaxis.set_pane_color((1.0, 1.0, 1.0, 1.0))
ax.plot(x, y, z, c='royalblue')

ax.set_xlabel('$x(t)$')
ax.set_ylabel('$y(t)$')
ax.set_zlabel('$z(t)$')

# Fix to get z_label inside the figure
from matplotlib import rcParams
rcParams.update({'figure.autolayout': True})

ax.view_init(9,-44)
plt.show()

# Subtract solution components for Poincar!{\color{comment}\'e}! section
sol0 = dde.sample(tfinal-1000, tfinal, 0.001)
t = sol0['t']
x = sol0['x']
y = sol0['y']
z = sol0['z']

# Poincar!{\color{comment}\'e}! section
def poincaresection(x1, x2, x3,x1_label, x2_label, val):
  zero_cross = np.where(np.diff(np.sign(x3-val)))
  plt.figure(1)
  plt.xlabel(x1_label)
  plt.ylabel(x2_label)
  plt.plot(x1[zero_cross], x2[zero_cross],'.', c='royalblue')
  plt.show()
  return

x1_label='$x(t)$'
x2_label='$y(t)$'
poincaresection(x,y,z,x1_label,x2_label,-21.75)

\end{lstlisting}
We will simulate the dynamics in region I and II, see \cref{sm:fig:RH_bifurcation_diagram_II}, where a stable periodic orbit and stable two-dimensional torus, respectively, should be present. For a point in region I, we will take the unfolding parameter values 
\[
(I_{app},S)=(-18.886177304147466,-8.044197084548104).
\]
As an initial condition, we use the constant history function with values of the location of the fold-Hopf point. We integrate the DDE on the time interval $t\in[0,10.000]$ with \PYDELAY. In  \cref{fig:Rose-Hindmarsh-simulation} the time series of the components $x$,$y$ and $z$ and the orbit in $(x,y,z)$-space of the last 10 time steps are shown, clearly indicating a stable orbit.

Next, we simulate the dynamics in region II. Therefore, we adjust the unfolding parameter values to 
\[
(I_{app},S)=(-18.902420391705071,-8.045234985422740).
\]
Furthermore, we increase the integration interval to $t\in[0,20.000]$. Keeping the history function the same, we plot the last 1000 time steps, see \cref{fig:Rose-Hindmarsh-simulation}. We conclude that the dynamics near the fold-Hopf point are as predicted in \cite{Kuznetsov2004}.

\section{Hopf-Hopf and generalized Hopf bifurcations in Active control system}
In \cref{sec:acs_example} we considered the following active control system
\begin{equation}
\begin{cases}
\begin{aligned}
\dot{x}(t)&=\tau y(t),\\
\dot{y}(t)&=\tau\left(-x(t)-g_{u}x(t-1)-2\zeta y(t)-g_{v}y(t-1)+f(t)\right),
\end{aligned}
\end{cases}\label{sm:eq:acs}
\end{equation}
which  is used to control the response of structures to internal or external excitation, see \cite{Peng2013}. The function $f$ is substituted by $\beta x^{3}(t-1)$ and the parameters
\[
g_{u}=0.1,\quad g_{v}=0.52,\quad\beta=0.1,
\]
are fixed. The control parameters are $\zeta$ and $\tau$. 

\begin{remark}
	This demonstration can be found in the directory \lstinline|demos/tutorial/VII/acs| relative to the main directory of the \DDEBIFTOOL package. Here, we omit the code to generate a system file. The system file \lstinline|sym_acs_mf.m| has been generated with the script \lstinline|gen_acs.m|. Also, we assume that the \DDEBIFTOOL package has been loaded as in \cref{sm:lst:searchpath}.
	The code in \crefrange{sm:sec:acs:pars_and_funcs}{sm:sec:acs_bifurcation_diagram} highlights the important parts of the file \lstinline|acs.m|.
\end{remark}

\subsection{Set parameter names and funcs structure} \label{sm:sec:acs:pars_and_funcs}
We set the parameter names and define the \lstinline|funcs| structure.
\begin{lstlisting}[style=customMatlab,escapechar=!]
%% Set parameter names
parnames={'zeta','tau','tau_scaled'};
cind=[parnames;num2cell(1:length(parnames))];
ind=struct(cind{:});
!\matlabrule!
%% Set funcs structure
% We load the precalculated multilinear forms. These have been
% generated with the file gen_sym_acs.m.
funcs=set_symfuncs(@sym_acs_mf,'sys_tau',@()ind.tau_scaled);
\end{lstlisting}

\subsection{Stability and normal form coefficients of the Hopf-Hopf point} We construct a steady-state at the Hopf-Hopf point and calculate its stability.
\begin{lstlisting}[style=customMatlab]
%% Hopf-Hopf point
% Construct steady-state point
stst.kind='stst';
stst.x=[0;0];
stst.parameter(ind.zeta)=-0.016225;
stst.parameter(ind.tau)=5.89802;
stst.parameter(ind.tau_scaled)=1;
% Calculate stability
method=df_mthod(funcs,'stst');
stst.stability=p_stabil(funcs,stst,method.stability);
stst.stability.l1
\end{lstlisting}
The \MATLAB console shows the following output.
\begin{lstlisting}[style=matlabConsole]
ans =

    0.0000 + 4.5275i
    0.0000 - 4.5275i
   -0.0000 + 7.6449i
   -0.0000 - 7.6449i
\end{lstlisting}
The eigenvalues confirm that the point under consideration is indeed a Hopf-Hopf point. Furthermore, the remaining eigenvalues have negative real parts. Next, we calculate the normal form coefficients and the transformation to the center manifold with the function \lstinline|nmfm_hoho|, which implements the coefficients as derived in \cref{sec:HH_coef}. For this we need to set the argument \lstinline|free_pars| to the unfolding parameter $(\zeta,\tau)$. These coefficients will be used to start the continuation of the various branches emanating from the Hopf-Hopf point.
\begin{lstlisting}[style=customMatlab]
%% Calculate coefficients of the parameter dependent normal form
hopf=p_tohopf(funcs,stst);
method=df_mthod(funcs,'hopf');
hopf.stability=p_stabil(funcs,hopf,method.stability);
hoho=p_tohoho(hopf);
unfolding_parameters=[ind.zeta, ind.tau];
hoho=nmfm_hoho(funcs,hoho,'free_pars',unfolding_parameters);
hoho.nmfm
\end{lstlisting}
The \MATLAB console shows the following output.
\begin{lstlisting}[style=matlabConsole,keepspaces=true]
ans = 

?struct? with fields:

    g2100: -0.0915 + 0.1214i
    g1011: -0.3084 + 0.4096i
    g1110: 0.2151 + 0.3876i
    g0021: 0.1813 + 0.3268i
    theta: -1.7009
    delta: -2.3517
        b: [2x2 double]
    h0011: [1x1 struct]
    h0020: [1x1 struct]
    h2000: [1x1 struct]
        K: [2x2 double]
    h0000: [1x2 struct]
\end{lstlisting}
We conclude that this Hopf-Hopf bifurcation is of `difficult' type, since
\[
( \Re g_{2100} )( \Re g_{0021} ) = -0.0166 < 0,
\]
see \cite{Kuznetsov2004}. Furthermore, the quantities
\[
\theta=\frac{\text{Re }g_{1011}}{\text{Re }g_{0021}}=-1.7009,\qquad\delta=\frac{\text{Re }g_{1101}}{\text{Re }g_{2100}}=-2.3517
\]
are such that $\theta<0,\,\delta<0,\,\theta\delta>0$. It follows that we are in case VI.

\subsection{Set bifurcation parameter range and step size bounds}
Before continuing the various branches emanating from the transcritical-Hopf point, we create the variable \lstinline|brpars| containing parameter bounds and maximal stepsizes.
\begin{lstlisting}[style=customMatlab]
%% Set bifurcation parameter range and step size bounds
brpars={'max_bound',[ind.tau 16],...
        'min_bound',[ind.tau 5],...
        'max_step', [ind.zeta 0.04; ind.tau 0.04]};
\end{lstlisting}

\subsection{Continuing Hopf and Neimark-Sacker bifurcation curves}
\begin{figure}
\ifcompileimages
\subfloat[]{%
  \tikzsetnextfilename{acs_ns_branches}%
  \input{tikz/acs_ns_branches}%
}
\hfill
\subfloat[\label{sm:fig:acs_hopf_branches}]{%
  \tikzsetnextfilename{acs_hopf_branches}%
  \input{tikz/acs_hopf_branches}%
}
\else
\subfloat[]{\includegraphics{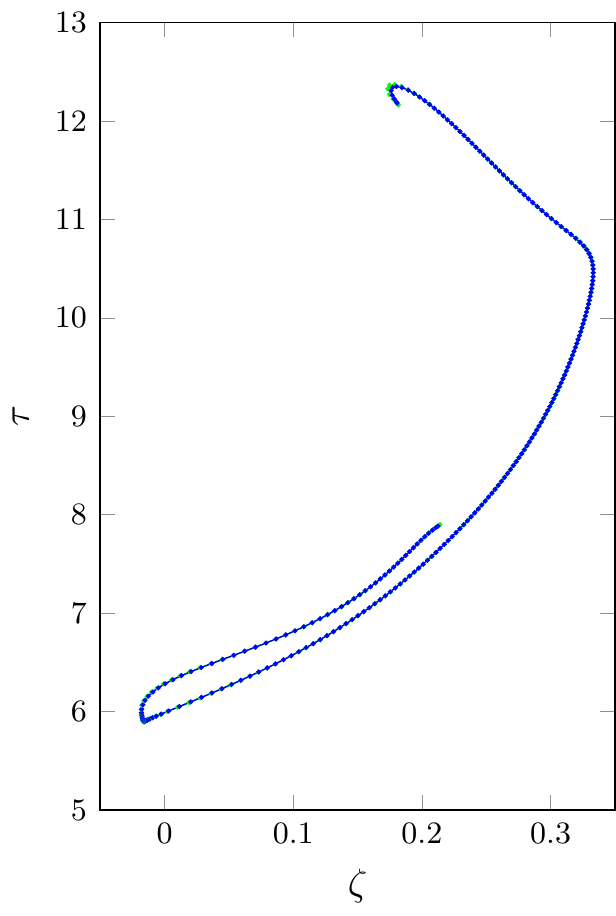}}
\hfill
\subfloat[]{\includegraphics{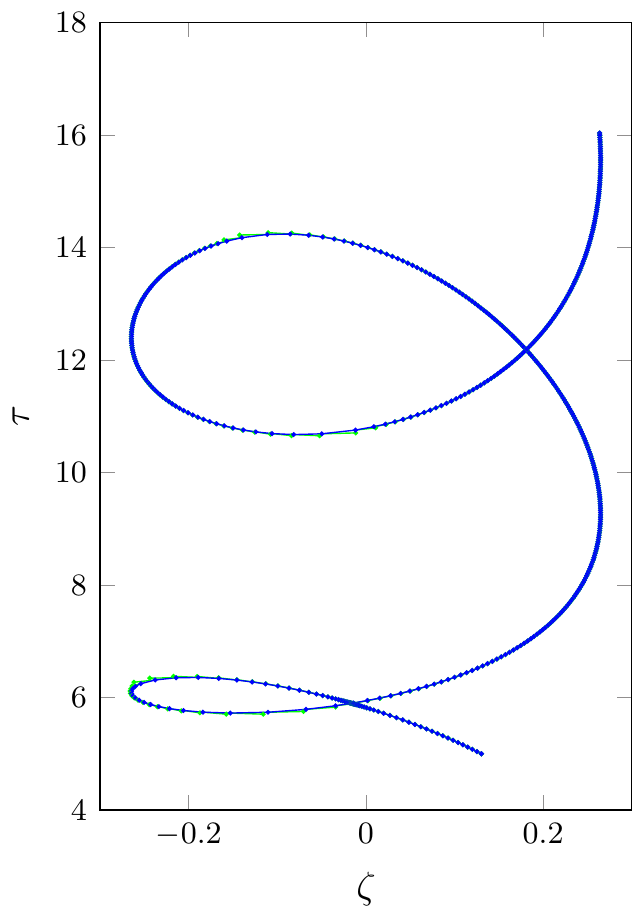}\label{sm:fig:acs_hopf_branches}}
\fi
\caption{In \textup{(a)}, the continued branches \lstinline|ns1_br| and \lstinline|ns2_br| are plotted. In \textup{(b)}, the continued Hopf branch \lstinline|hopf_br1| is plotted.}
\label{sm:fig:continued_branches}
\end{figure}

We use the function \lstinline|C1branch_from_C2point| to start to continue the branches emanating from the Hopf-Hopf point.
\begin{lstlisting}[style=customMatlab]
%% Continue Neimark-Sacker curves emanating from Hopf-Hopf point
[trfuncs,nsbr,suc]=C1branch_from_C2point(funcs,hoho,...
    unfolding_parameters,'codim2','hoho','codim1',....
    'TorusBifurcation',brpars{:},'step',1e-01,'plot',0);
assert(all(suc(:)>0))
ntrsteps=186; [nsbr(1),suc]=br_contn(trfuncs,nsbr(1),ntrsteps);
assert(suc>0)
ntrsteps=61;  [nsbr(2),suc]=br_contn(trfuncs,nsbr(2),ntrsteps);
assert(suc>0)
\end{lstlisting}
\begin{lstlisting}[style=customMatlab]
%% Continue Hopf curve emanating from Hopf-Hopf point
[~,hbr,suc]=C1branch_from_C2point(funcs,hoho,...
    unfolding_parameters,'codim2','hoho','codim1','hopf'....
    ,brpars{:},'step',1e-03,'plot',0);
assert(all(suc(:)>0))
nop=1000; [hbr(1),suc]=br_contn(funcs,hbr(1),nop);assert(suc>0)
hbr(1)=br_rvers(hbr(1));
[hbr(1),suc]=br_contn(funcs,hbr(1),nop); assert(suc>0)
nop=10; [hbr(2),suc]=br_contn(funcs,hbr(2),nop); assert(suc>0)
hbr(2)=br_rvers(hbr(2));
[hbr(2),suc]=br_contn(funcs,hbr(2),nop); assert(suc>0)
\end{lstlisting}
In \cref{sm:fig:continued_branches}, the computed branches are shown. We see in \cref{sm:fig:acs_hopf_branches} that it is redundant to continue the second Hopf branch emanating from the Hopf-Hopf point. Indeed, the first Hopf branch connects the Hopf-Hopf point to itself. Here, we verified that the underlying points coincide.

\subsection{Predictors}
For comparison in the bifurcation diagram, we obtain predictors for the various branches by setting the argument \lstinline|predictor| to 1 and \lstinline|step| to a range of $\varepsilon$-values when calling the function \lstinline|C1branch_from_C2point|.
\begin{lstlisting}[style=customMatlab]
%% Predictors for Neimark-Sacker and Hopf curves
[trfuncs,nsbr_pred,suc]=C1branch_from_C2point(funcs,hoho,...
    unfolding_parameters,'codim2','hoho','codim1',...
    'TorusBifurcation','step',linspace(1e-03,2,40),'predictor',1);
[~,hbr_pred,suc]=C1branch_from_C2point(funcs,hoho,...
    unfolding_parameters,'codim2','hoho','codim1','hopf',...
    'step',linspace(-2e-01,2e-01,30),'predictor',1);
\end{lstlisting}

\subsection{Bifurcation diagram}
We plot the computed curves and the predictors for the Neimark-Sacker and Hopf curves with the following code.
\begin{lstlisting}[style=customMatlab]
%% Close-up near Hopf Hopf point in parameter space with predictors
figure(2); clf; hold on;
hbr1_pred_pm  = [getpars(hbr_pred(1), ind.zeta)
                 getpars(hbr_pred(1), ind.tau)];
hbr2_pred_pm  = [getpars(hbr_pred(2), ind.zeta)
                 getpars(hbr_pred(2), ind.tau)];
nsbr1_pred_pm = [getpars(nsbr_pred(1),ind.zeta)
                 getpars(nsbr_pred(1),ind.tau)];
nsbr2_pred_pm = [getpars(nsbr_pred(2),ind.zeta)
                 getpars(nsbr_pred(2),ind.tau)];
plot(hbr_pm(1,:),hbr_pm(2,:),'Color',cm(1,:),...
    'DisplayName','Hopf branches')
plot(nsbr1_pm(1,:),nsbr1_pm(2,:),'Color',cm(2,:),...
    'DisplayName','Neimark-Sacker branches')
plot(nsbr2_pm(1,:),nsbr2_pm(2,:),'Color',cm(2,:),...
    'HandleVisibility','off')
plot(hbr1_pred_pm(1,:), hbr1_pred_pm(2,:) ,'.',...
    'Color',cm(1,:),'DisplayName','Hopf predictors')
plot(hbr2_pred_pm(1,:), hbr2_pred_pm(2,:) ,'.',...
    'Color',cm(1,:),'HandleVisibility','off')
plot(nsbr1_pred_pm(1,:),nsbr1_pred_pm(2,:),'.',...
    'Color',cm(2,:),'DisplayName','Neimark-Sacker predictors')
plot(nsbr2_pred_pm(1,:),nsbr2_pred_pm(2,:),'.',...
    'Color',cm(2,:),'HandleVisibility','off')
plot(hoho.parameter(ind.zeta),hoho.parameter(ind.tau),'k.',...
    'MarkerSize',12,'DisplayName','Hopf Hopf point')
title(['Close-up near Hopf Hopf point '...
       'in parameter space with predictors'])
xlabel('$\zeta$','Interpreter','LaTex');
ylabel('$\tau$','Interpreter','LaTex');
legend('Location','NorthWest')
axis([-0.0562    0.0109    5.7098    6.1757])
box on
\end{lstlisting}
In \cref{sm:fig:acs_hoho_predictors} the predictors in parameter space are compared.

\subsection{Plot comparing computed and predicted periods}
To compare the computed and predicted periods, we use the following code.
\begin{lstlisting}[style=customMatlab]
%% Compare computed and predicted periods
figure(3); clf; hold on;
% Plot computed periods on nsbr(1) and nsbr(2)
omegas1=arrayfun(@(p)p.period,nsbr(1).point);
omegas2=arrayfun(@(p)p.period,nsbr(2).point);
omegas1_pred=arrayfun(@(p)p.period,nsbr_pred(1).point);
omegas2_pred=arrayfun(@(p)p.period,nsbr_pred(2).point);
plot(getpars(nsbr(1),ind.zeta),omegas1,'Color',cm(1,:));
plot(getpars(nsbr(2),ind.zeta),omegas2,'Color',cm(2,:));
plot(getpars(nsbr_pred(1),ind.zeta),...
     omegas1_pred,'.','Color',cm(1,:));
plot(getpars(nsbr_pred(2),ind.zeta),...
     omegas2_pred,'.','Color',cm(2,:));
title('Compare computed and predicted periods')
xlabel('$\zeta$','Interpreter','LaTex');
ylabel('$\omega$','Interpreter','LaTex');
\end{lstlisting}
In \cref{sm:fig:acs_compare_periods}, the computed and predicted periods of the cycles are compared using the formulas as given in \cref{eq:HH_period_predictors}.

\begin{figure}
\centering
\ifcompileimages
\subfloat[]{%
  \tikzsetnextfilename{acs_hoho_predictors_for_sm}%
  \input{tikz/acs_hoho_predictors_for_sm}%
\label{sm:fig:acs_hoho_predictors}} \hfill
\subfloat[]{%
  \tikzsetnextfilename{acs_compare_omegas}%
  \input{tikz/acs_compare_omegas}%
\label{sm:fig:acs_compare_periods}}
\else
\subfloat[]{\includegraphics{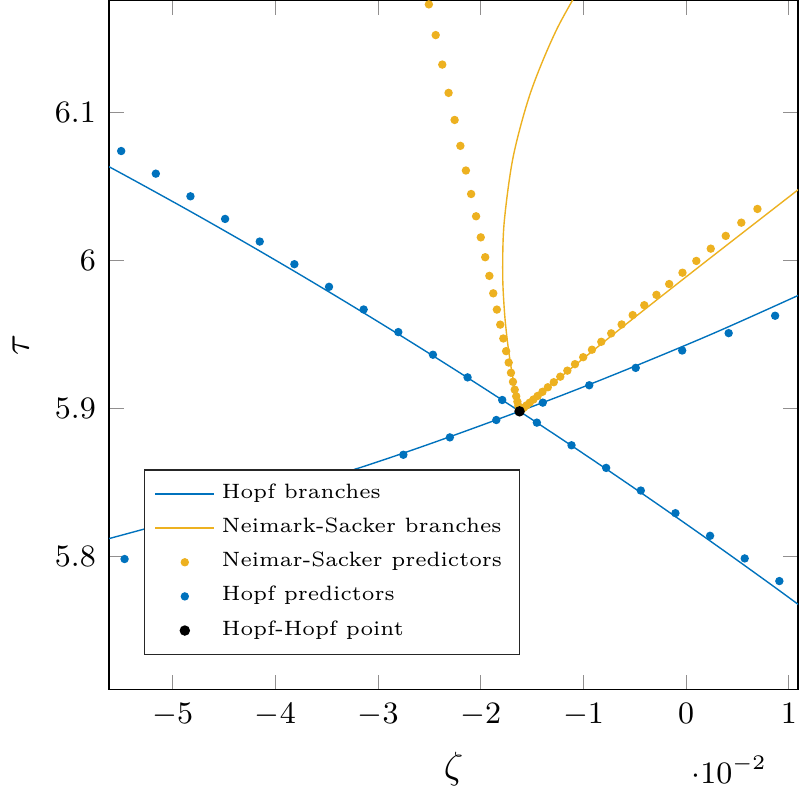}\label{sm:fig:acs_hoho_predictors}} \hfill
\subfloat[]{\includegraphics{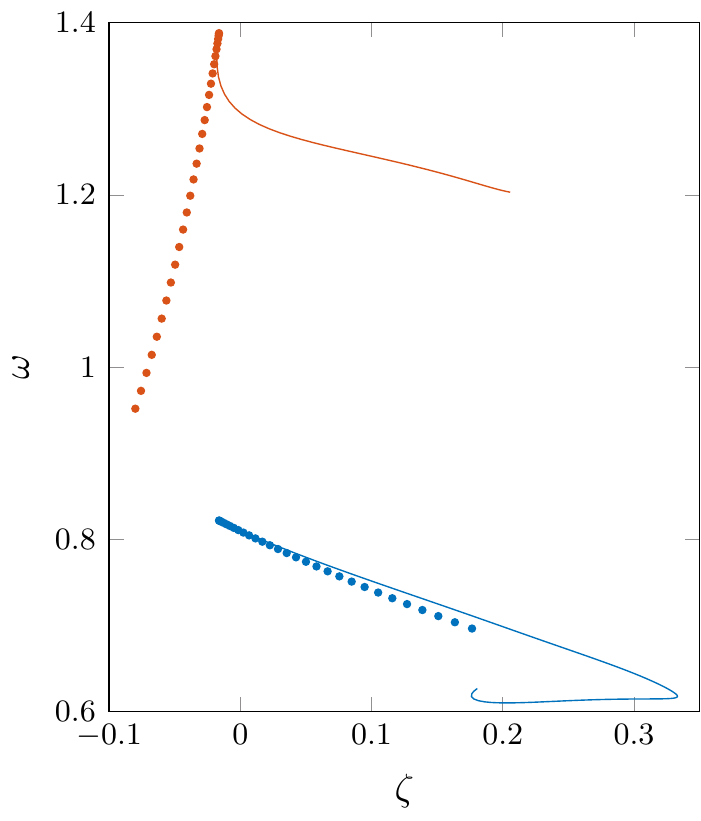}\label{sm:fig:acs_compare_periods}}
\fi
\caption{In \textup{(a)}, the computed curves are compared with the predicted curves for the Hopf-Hopf point \lstinline|hoho|. In \textup{(b)}, the computed and predicted periods of the cycles are compared.}
\end{figure}

\subsection{Detect special points on the Hopf branch}
Using the detection capabilities from \DDEBIFTOOL via the function \lstinline|LocateSpecialPoints|, we detect one additional Hopf-Hopf point and
three generalized Hopf points.
\begin{lstlisting}[style=customMatlab]
%% Detect codimension two points on hopf_br1
[hopf_br_wbifs,hopftests,hc2_indices,hc2_types]=...
    LocateSpecialPoints(funcs,hbr(1));
\end{lstlisting}
The \MATLAB console shows the following output.
\begin{lstlisting}[style=matlabConsole,keepspaces=true]
HopfCodimension2: calculate stability if not yet present
HopfCodimension2: calculate L1 coefficients
HopfCodimension2: (provisional) 3 gen. Hopf 4 Hopf-Hopf  detected.
br_insert: detected 1 of 7: hoho. Normalform:
     g2100: 0.1813 + 0.3268i
     g1011: 0.2151 + 0.3876i
     g1110: -0.3084 + 0.4096i
     g0021: -0.0915 + 0.1214i
     theta: -2.3517
     delta: -1.7009

br_insert: detected 2 of 7: genh. Normalform:
     L2: -0.0458
     L1: -1.6706e-11

br_insert: detected 3 of 7: hoho. Normalform:
     g2100: 0.1813 + 0.3268i
     g1011: 0.2151 + 0.3876i
     g1110: -0.3084 + 0.4096i
     g0021: -0.0915 + 0.1214i
     theta: -2.3517
     delta: -1.7009

br_insert: detected 4 of 7: hoho. Normalform:
    g2100: -0.0034 + 0.3217i
    g1011: -0.0046 + 0.4380i
    g1110: -0.1494 + 0.4498i
    g0021: -0.0509 + 0.1531i
    theta: 0.0912
    delta: 43.8639

br_insert: detected 5 of 7: genh. Normalform:
    L2: 0.0070
    L1: -5.9006e-14

br_insert: detected 6 of 7: genh. Normalform:
    L2: 0.0075
    L1: 2.2535e-14

br_insert: detected 7 of 7: hoho. Normalform:
    g2100: -0.0034 + 0.3217i
    g1011: -0.0046 + 0.4380i
    g1110: -0.1494 + 0.4498i
    g0021: -0.0509 + 0.1531i
    theta: 0.0912
    delta: 43.8639
\end{lstlisting}
There are seven bifurcations detected on the Hopf branch \lstinline|hbr(1)|: four Hopf-Hopf bifurcations and three generalized Hopf bifurcations. However, by inspecting the parameters of the detected Hopf Hopf points suggest that there are only two distinct Hopf-Hopf points,  which are connected with the same Hopf branch. By comparing the location of the underlying points confirms the premise. We are therefore not interested in continuing the Hopf curves emanating from these points. What we are interested in is continuing the Neimark-Sacker and limit point of cycle curves emanating from the second Hopf-Hopf point and the three generalized Hopf points, respectively.

The normal form coefficients of the second Hopf-Hopf point are such that
\[
( \Re g_{2100} )( \Re g_{0021} ) = 1.7331e - 04 > 0
\]
and
\[
\theta\geq \delta > 0, \qquad \delta\theta=4>1.
\]
We conclude that we are in case I of the `simple' type, see \cite[page 360]{Kuznetsov2008}. Therefore, no stable invariant two-dimensional torus is predicted for nearby parameter values. We only expect to find two stable periodic orbits.

\subsection{Continuing third Neimark-Sacker bifurcation and limit point of cycle bifurcation curves}
It turns out that we only need to continue one of the Neimark-Sacker bifurcation curves emanating from the second
Hopf-Hopf point. Indeed, the other Neimark-Sacker bifurcation curve is given
by \lstinline|ns_br(1)|.
\begin{lstlisting}[style=customMatlab,escapechar=!]
%% Subtract generalized Hopf points
genh_indices=hc2_indices(strcmp(hc2_types,'genh'));
genhpts=hopf_br_wbifs.point(genh_indices);
!\matlabrule!
%% Continue limit point of cycles emanating from the generalized Hopf points
[~,lpc_br1,suc]=C1branch_from_C2point(funcs,genhpts(1),...
unfolding_parameters,'codim2','genh','codim1','POfold',...
brpars{:},'step',1e-01,'plot',0);
assert(all(suc>0))
[poffuncs,lpc_br2,suc]=C1branch_from_C2point(funcs,genhpts(3),...
unfolding_parameters,'codim2','genh','codim1','POfold',...
brpars{:},'step',5e-03,'plot',0);
assert(all(suc>0))
ntrsteps=193; [lpc_br1,suc]=br_contn(poffuncs,lpc_br1,ntrsteps);
assert(suc>0)
ntrsteps=220; [lpc_br2,suc]=br_contn(poffuncs,lpc_br2,ntrsteps);
assert(suc>0)
!\matlabrule!
%% Subtract second Hopf-Hopf point on the Hopf branch
hoho_indices=hc2_indices(strcmp(hc2_types,'hoho'));
hoho2=hopf_br_wbifs.point(hoho_indices(3));
!\matlabrule!
%% Continue Neimark-Sacker curves emanating from Hopf-Hopf point
[trfuncs,nsbr34,suc]=C1branch_from_C2point(funcs,hoho2,...
    unfolding_parameters,'codim2','hoho','codim1',...
    'TorusBifurcation',brpars{:},'step',5e-03,'plot',0);
assert(all(suc(:)>0))
ntrsteps=200; [nsbr34(2),suc,fail,rjct]=...
    br_contn(trfuncs,nsbr34(2),ntrsteps);
\end{lstlisting}

\subsection{Bifurcation diagram} \label{sm:sec:acs_bifurcation_diagram}
\begin{figure}
\centering
\ifcompileimages
  \tikzsetnextfilename{acs_hoho_genh_lpc_ns}%
  \input{tikz/acs_hoho_genh_lpc_ns}%

\else
\includegraphics{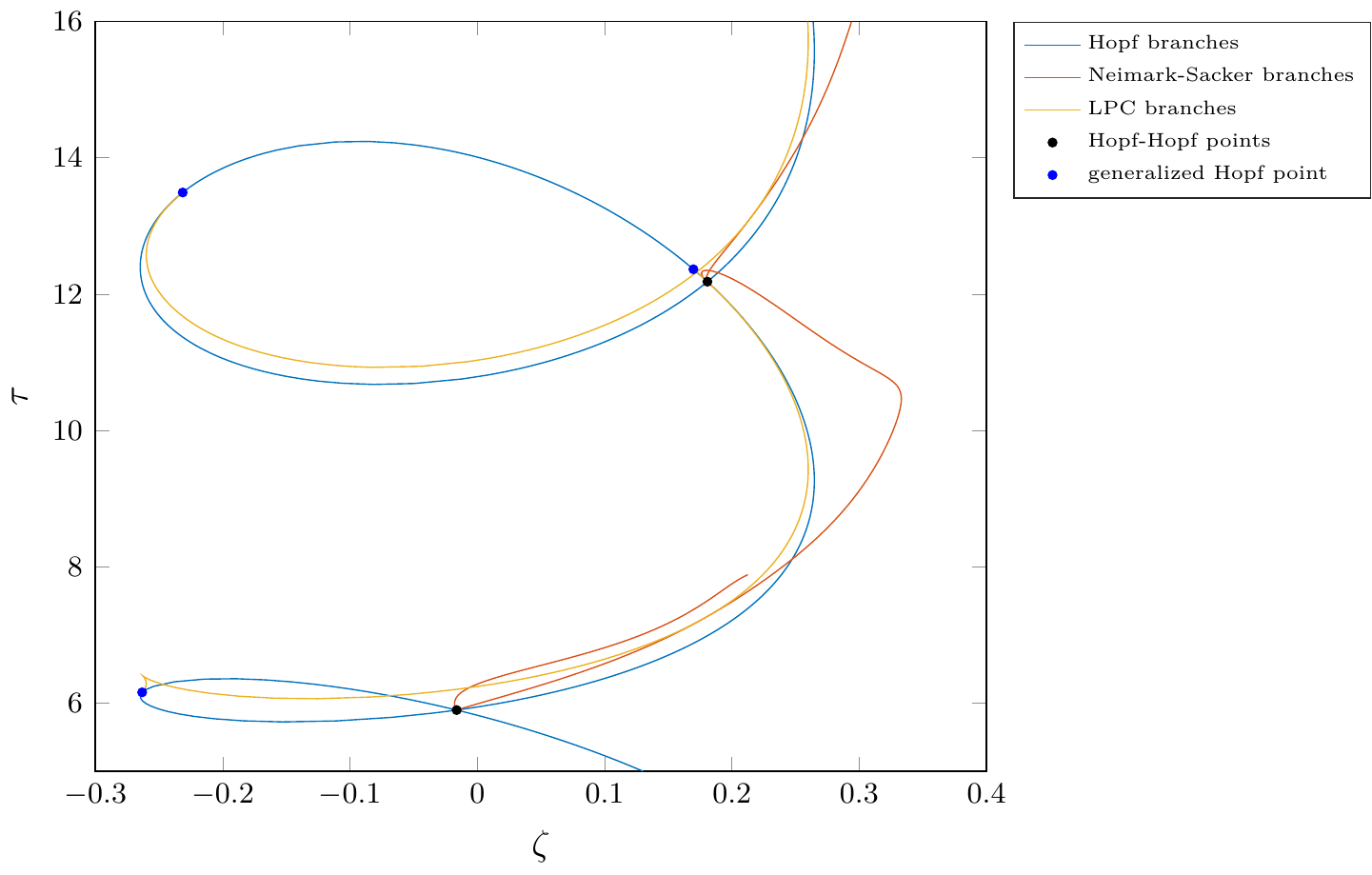}
\fi
\caption{Unfolding from multiple codimension two points detected in the DDE \cref{sm:eq:acs}.}
\label{sm:fig:acs_hoho_genh_lpc_ns}
\end{figure}

We plot the computed degenerate Hopf points, the limit point of limit cycle curves, and the Neimark-Sacker curves with the following code.
\begin{lstlisting}[style=customMatlab]
%% Bifurcation diagram in $(\zeta,\tau)$
figure(4); clf; hold on;
% Subtract paramater values from the branches
lpc_br1_pm= [getpars(lpc_br1,   ind.zeta)
             getpars(lpc_br1,   ind.tau)];
lpc_br2_pm= [getpars(lpc_br2,   ind.zeta)
             getpars(lpc_br2,   ind.tau)];
nsbr3_pm  = [getpars(nsbr34(2), ind.zeta)
             getpars(nsbr34(2) ,ind.tau)];
% Plot curves
plot(hbr_pm(1,:),    hbr_pm(2,:),    'Color',cm(1,:),...
    'DisplayName','Hopf branches');
plot(nsbr1_pm(1,:),  nsbr1_pm(2,:),  'Color',cm(2,:),...
    'DisplayName','Neimark-Sacker branches');
plot(nsbr2_pm(1,:),  nsbr2_pm(2,:),  'Color',cm(2,:),...
    'HandleVisibility','off');
plot(nsbr3_pm(1,:),  nsbr3_pm(2,:),  'Color',cm(2,:),...
    'HandleVisibility','off');
plot(lpc_br1_pm(1,:),lpc_br1_pm(2,:),'Color',cm(3,:),...
    'DisplayName','Generalized Hopf');
plot(lpc_br2_pm(1,:),lpc_br2_pm(2,:),'Color',cm(3,:),...
    'HandleVisibility','off');
% Add bifurcation points
plot(hoho.parameter(ind.zeta), hoho.parameter(ind.tau),...
    'k.','MarkerSize',8,'DisplayName','Hopf Hopf point')
plot(hoho2.parameter(ind.zeta),hoho2.parameter(ind.tau),...
    'k.','MarkerSize',8,'HandleVisibility','off')
plot(genhpts(1).parameter(ind.zeta),...
    genhpts(1).parameter(ind.tau),'b.','MarkerSize',8,...
    'DisplayName','Generelized Hopf point')
plot(genhpts(2).parameter(ind.zeta),...
     genhpts(2).parameter(ind.tau),'b.','MarkerSize',8,...
     'HandleVisibility','off')
plot(genhpts(3).parameter(ind.zeta),
     genhpts(3).parameter(ind.tau),'b.','MarkerSize',8,...
     'HandleVisibility','off')
title('Bifurcation diagram in (\zeta,\tau)')
xlabel('$\zeta$','Interpreter','LaTex')
ylabel('$\tau$','Interpreter','LaTex')
axis([-0.3    0.4    5    16])
legend()
\end{lstlisting}
\cref{sm:fig:acs_hoho_genh_lpc_ns} shows the resulting bifurcation diagram. There, we see that two generalized Hopf points are connected by the same limit point of cycles bifurcation curve.

\subsection{Simulation near Hopf-Hopf point with \texttt{pydelay}}
\begin{figure}[ht]
\subfloat[\label{fig:-Bifurcation-diagram-hoho}Bifurcation diagram near Hopf-Hopf point]{
\centering
\ifcompileimages
  \tikzsetnextfilename{acs_hoho_zoom}%
  \input{tikz/acs_hoho_zoom}%

\else
\includegraphics{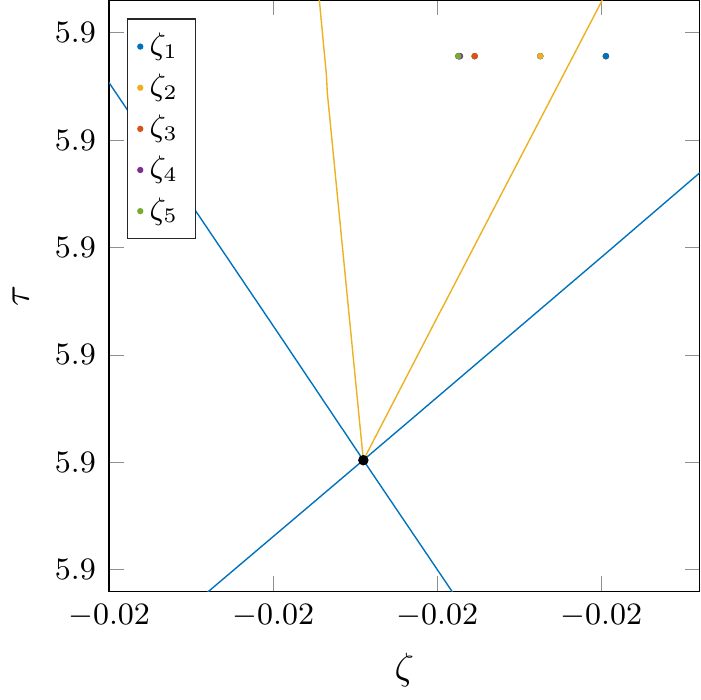}
\fi
}
\hspace*{\fill}
\subfloat[\label{fig:acs_zeta1}$\zeta_{1}$ Stable periodic orbit]{
\centering
\ifcompileimages
  \tikzsetnextfilename{acs_simulation_zeta1}%
  \input{tikz/acs_simulation_zeta1}%

\else
\includegraphics{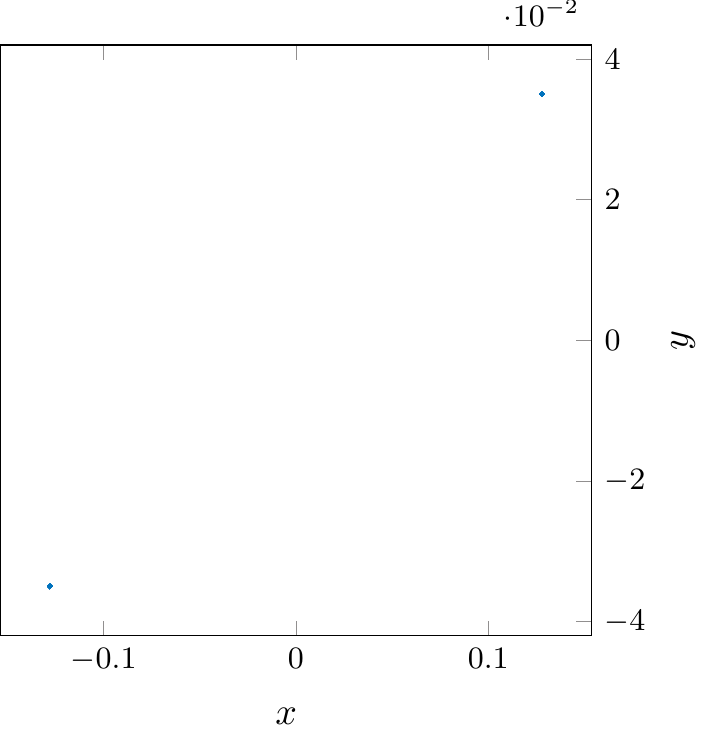}
\fi
}
\caption{In \textup{(a)}, the points $(\zeta_i,\tau)$ for $i=1,\dots,5$ of the parameter values where the simulation is performed are plotted. Note that $\zeta_4$ and $\zeta_5$ are almost indistinguishable. In \textup{(b)}, the Poincar\'e section of the simulation using \PYDELAY at the parameters ($\zeta_2,\tau)=(-0.015685728828307,5.901783308978358)$ shows a stable periodic orbit.}
\end{figure}

\begin{figure}[ht]
\centering
\subfloat[\label{fig:acs_zeta2}$\zeta_{2}$ Stable two-dimensional torus]{
\ifcompileimages
  \tikzsetnextfilename{acs_simulation_zeta2}%
  \input{tikz/acs_simulation_zeta2}%

\else
\includegraphics{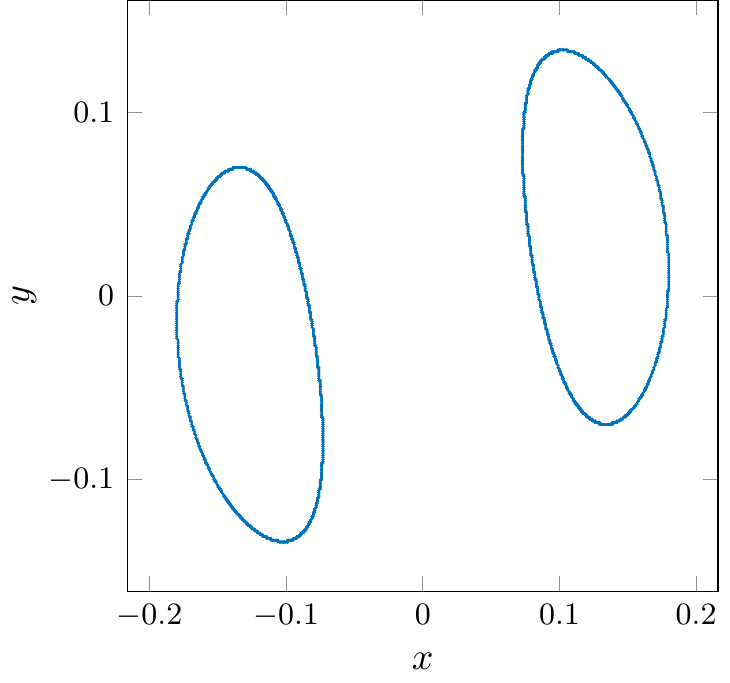}
\fi
}
\hspace*{\fill}
\subfloat[\label{fig:acs_zeta3}$\zeta_{3}$ Stable two-dimensional torus]{
\ifcompileimages
  \tikzsetnextfilename{acs_simulation_zeta3}%
  \input{tikz/acs_simulation_zeta3}%

\else
\includegraphics{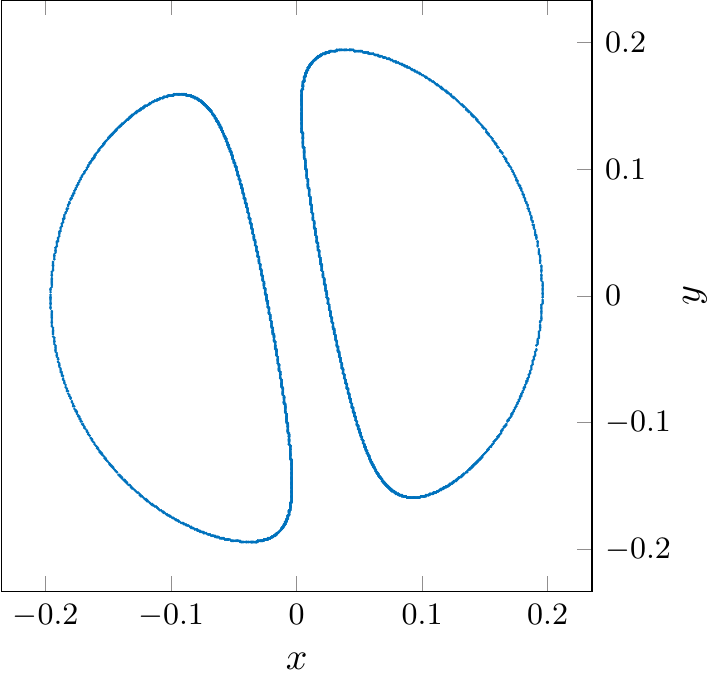}
\fi
} \\
\subfloat[\label{fig:zeta4}$\zeta_{4}$ Stable three-dimensional torus]{
\ifcompileimages
  \tikzsetnextfilename{acs_simulation_zeta4}%
  \input{tikz/acs_simulation_zeta4}%

\else
\includegraphics{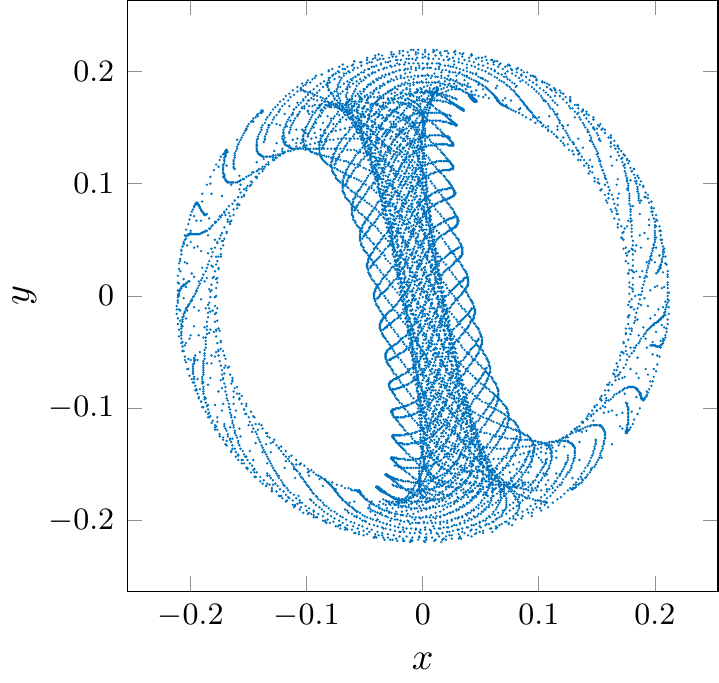}
\fi
}
\hspace*{\fill}
\subfloat[\label{fig:acs_zeta5}$\zeta_{5}$ Three-dimensional torus near blow-up]{
\ifcompileimages
  \tikzsetnextfilename{acs_simulation_zeta5}%
  \input{tikz/acs_simulation_zeta5}%

\else
\includegraphics{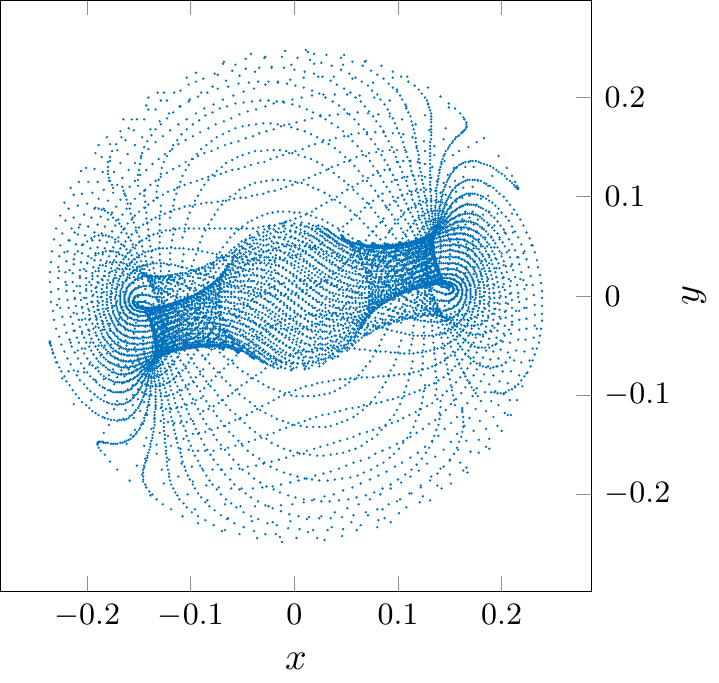}
\fi
}
\caption{Simulation with \texttt{pydelay} illustrating the branching of a three-dimensional torus from a two-dimensional torus. We refer to the text for further description.}
\end{figure}
In this last Section, we simulate the dynamics near the manually constructed Hopf-Hopf point at parameter values \cref{eq:acs-HH-pm}. As remarked before, the unfolding of the Hopf-Hopf point \blist{hoho} is of `difficult' type case VI. The normal form coefficients predict a stable invariant two-dimensional torus. Furthermore, this torus undergoes a bifurcation in which a three-dimensional torus is born. Since the DDE under investigation \cref{eq:acs3-1} does not contain any terms of order higher than three, the results for the normal form remain valid for the system. To confirm the unfolding in  \cref{fig:acs_hoho_predictors}, we fix the delay $\tau=5.901783308978358$ and take for $\zeta$ consecutive the values
\begin{align*}
\zeta_{1} & =-0.015485728828307,\\
\zeta_{2} & =\zeta_{1}-0.0002,\\
\zeta_{3} & =\zeta_{1}-0.0004,\\
\zeta_{4} & =\zeta_{1}-0.000445,\\
\zeta_{5} & =\zeta_{1}-0.0004496,
\end{align*}
see \cref{fig:-Bifurcation-diagram-hoho}. The cross-sections in \cref{fig:acs_zeta1,fig:acs_zeta5} are generated with the following Python code using \PYDELAY, see also the file \lstinline|acs_simulation.py|
\begin{lstlisting}[language=Python]
import numpy as np
import pylab as pl
from pydelay import dde23
from matplotlib import colors

# Number of time units
tfinal = 90000

# Define DDE
eqns = {
    'x':'tau*y',
    'y':'tau*(-x-0.1*x(t-1)-2*zeta*y' \
         '-0.52*y(t-1)+0.1*pow(x(t-1),3))'
}

# Set parameters
tau=5.901783308978358
zeta1=-0.015485728828307 # periodic orbit
zeta2=zeta1-0.0002       # torus
zeta3=zeta1-0.0004       # torus near bifurcation to 3d torus
zeta4=zeta1-0.000445     # 3d torus
zeta5=zeta1-0.0004496    # 3d torus near destruction

# Solve DDE
dde = dde23(eqns=eqns, params={'zeta':zeta2, 'tau':tau})
dde.set_sim_params(tfinal=tfinal, dtmax=0.1, AbsTol=1e-08, 
	RelTol=1e-06)
histfunc = {'x': lambda t: 0.01, 'y': lambda t: 0 }
dde.hist_from_funcs(histfunc, 51)
dde.run()

# Subtract solution components
a1=1;
dt=1e-04
sol = dde.sample(tfinal-tfinal/10, tfinal, dt)
soldelayed = dde.sample(tfinal-tfinal/10-a1, tfinal-a1, dt)
t = sol['t']
x = sol['x']
y = sol['y']
xdelayed = soldelayed['x']
zero_crossings = np.where(np.diff(np.sign(xdelayed)))[0]

# Scatter plot
x_cros=x[zero_crossings]
y_cros=y[zero_crossings]
params = {
    'figure.figsize': (15, 15),
    'axes.labelsize': 'x-large',
    'xtick.labelsize':'x-large',
    'ytick.labelsize':'x-large'
}
pl.rcParams.update(params)
fig = pl.figure()
pl.figure(1)
pl.scatter(x_cros,y_cros,s=0.8,color='royalblue')
pl.xlabel('$x$')
pl.ylabel('$y$')
pl.show()
\end{lstlisting}

In \cref{fig:acs_zeta1} there are two dots, corresponding to a stable period orbit. Crossing the curve \blist{ns1\_br}, a stable two-dimensional torus branches off, see \cref{fig:acs_zeta2}. The torus still exists at $\zeta=\zeta_{3}$, as seen in  \cref{fig:acs_zeta3}. Then, at $\zeta=\zeta_{4}$, only slightly smaller than $\zeta_{3}$, a three-dimensional torus is observed, see \cref{fig:zeta4}. Lastly, \cref{fig:acs_zeta5} shows the three-dimensional torus near the curve where the torus blows up.

\section{Transcritical-Hopf bifurcation in Van der Pol's oscillator with delayed position and velocity feedback}
Consider the generalized van der Pol's oscillator with delayed feedback\begin{equation}
\begin{cases}
\begin{aligned}
\dot{x}(t)&=\left(\tau_0+\mu_2\right)y(t),\\[0.5em]
\dot{y}(t)&=\left(\tau_0+\mu_{2}\right)\big[-x(t)-\varepsilon(x^2(t)-1)y(t)+(1+\mu_1)x(t-\tau)-0.2y(t-1)\\[0.5em]
&\qquad-0.2x^2(t-1)-0.2x(t-\tau)y(t-1)-0.2 y^2(t-1)+0.5x^3(t-1)\big],
\end{aligned}
\end{cases}\label{sm:eq:vdp}
\end{equation}
see \cref{sec:HT_example} and \cite{Bramburger2014}. The parameter $\varepsilon=0.3$ is fixed. For $\tau_{0}\approx1.757290761249588$ a transcritical-Hopf bifurcation is located at $(\mu_1,\mu_2)=(0,0)$.

\begin{remark}
	This demonstration can be found in the directory \lstinline|demos/tutorial/VII/vdpo| relative to the main directory of the \DDEBIFTOOL package. Here, we omit the code to generate a system file. The system file \lstinline|sym_vdpo_mf.m| has been generated with the script \lstinline|gen_sym_vdpo.m|. Also, we assume that the \DDEBIFTOOL package has been loaded as in \cref{sm:lst:searchpath}. The code in \crefrange{sm:sec:vdpo:pars_and_funcs}{sm:sec:vdpo:comparing_period_orbits} highlights the important parts of the file \lstinline|vdpo.m|.
\end{remark}

\subsection{Set parameter names and funcs structure} \label{sm:sec:vdpo:pars_and_funcs}
We set the parameter names and define the \lstinline|funcs| structure.
\begin{lstlisting}[style=customMatlab]
%% Set parameter names
parnames={'mu1','mu2','tau'};
cind=[parnames;num2cell(1:length(parnames))];
ind=struct(cind{:});
%% Set funcs structure
% We load the precalculated multilinear forms. These have been
% generated with the file gen_sym_vdpo.m.
funcs=set_symfuncs(@sym_vdpo_mf,'sys_tau',@()ind.tau);
\end{lstlisting}

\subsection{Stability and normal form coefficients of the transcritical-Hopf point} We construct a steady-state at the transcritical-Hopf point and calculate its stability.
\begin{lstlisting}[style=customMatlab]
%% Construct transcritical-Hopf bifucation point
stst=dde_stst_create('x',[0;0],'parameter',[0 0 1]);
% Calculate stability
method=df_mthod(funcs,'stst');
stst.stability=p_stabil(funcs,stst,method.stability);
stst.stability.l1
\end{lstlisting}
The \MATLAB console shows the following output.
\begin{lstlisting}[style=matlabConsole]
ans =

    0.0000 + 0.0000i
   -0.0000 + 2.4539i
   -0.0000 - 2.4539i
\end{lstlisting}
We have a zero eigenvalue and a pair of purely imaginary eigenvalues. Furthermore, the remaining eigenvalues have negative real parts. Next, we calculate the normal form coefficients and the transformation to the center manifold with the function \lstinline|nmfm_zeho|, which implements the coefficients as derived in \cref{sec:fold-Hopf,sec:transcritical-Hopf}. For this we need to set the argument \lstinline|free_pars| to the unfolding parameter $(\mu_1,\mu_2)$. These coefficients will be used to start the continuation of the various branches emanating from the transcritical-Hopf point.
\begin{lstlisting}[style=customMatlab]
%% Coefficients of the parameter dependent normal form
ht=p_tohopf(funcs,stst);
ht=p_tozeho(ht);
unfolding_pars=[ind.mu1, ind.mu2];
ht=nmfm_zeho(funcs,ht,'transcritical',1,'free_pars',unfolding_pars);
ht.nmfm
\end{lstlisting}
The \MATLAB console shows the following output.
\begin{lstlisting}[style=matlabConsole,keepspaces=true]
ans = 

?struct? with fields:

             g200: 0.2121
             g110: -0.1337 + 0.2672i
             g011: 0.4241
             g300: 0.4935
             g111: 1.0243
             g210: -0.8178 - 0.4283i
             g021: -0.3302 - 0.1646i
                b: 0.2121
                c: 0.4241
                d: -0.1337 - 5.4430i
                e: -0.2435
                s: 0.0899
            theta: -0.6303
    transcritical: 1
             h200: [1x1 struct]
             h011: [1x1 struct]
             h020: [1x1 struct]
             h110: [1x1 struct]
                K: [2x2 double]
           omega1: 0.4644
           omega2: 1.2768
\end{lstlisting}
The normal form coefficients are such that
\[
g_{011} \Re ( g_{110}) = 0.4241 \Re (- 0.1337 + 0.2672i ) < 0.
\]
Therefore, there are two Neimark-Sacker bifurcation curves predicted, see \cref{sec:HT_predictors}.

\subsection{Set bifurcation parameter range and step size bounds}
Before continuing the various branches emanating from the transcritical-Hopf point, we create the variable \lstinline|brpars| containing parameter bounds and maximal stepsizes.
\begin{lstlisting}[style=customMatlab]
%% Set bifurcation parameter range and step size bounds
brpars={'max_bound',[ind.mu1 0.0139; ind.mu2 0.0094],...
        'min_bound',[ind.mu1 -0.0190; ind.mu2 -0.0069 ],...
        'max_step' ,[ind.mu1  1.0e-02; ind.mu2 1.0e-02]};
\end{lstlisting}

\subsection{Continuing Hopf, transcritical and Neimark-Sacker bifurcation curves}
We use the function \lstinline|C1branch_from_C2point| to continue the various branches emanating from the singularity.
\begin{lstlisting}[style=customMatlab,escapechar=!]
%% Continue Neimark-Sacker curves emanating from 
%  the transcritical-Hopf point
[trfuncs,nsbr,suc]=C1branch_from_C2point(funcs,ht,unfolding_pars,...
    'codim2','zeho','codim1','TorusBifurcation',...
    'step',1e-04,'plot',0,brpars{:});
assert(all(suc(:)>0))
ntrsteps=27; [nsbr(1),suc]=br_contn(trfuncs,nsbr(1),ntrsteps);
assert(suc>0)
ntrsteps=30; [nsbr(2),suc]=br_contn(trfuncs,nsbr(2),ntrsteps);
assert(suc>0)
!\matlabrule!
%% Continue Hopf curves emanating from fold-Hopf point
[~,hbr,suc]=C1branch_from_C2point(funcs,ht,unfolding_pars,...
    'codim2','zeho','codim1','hopf',brpars{:},'step',1e-05,'plot',0);
assert(all(suc(:)>0))
for i=2:-1:1
nop=1000; hbr(i)=br_contn(funcs,hbr(i),nop);
hbr(i)=br_rvers(hbr(i));
hbr(i)=br_contn(funcs,hbr(i),nop);
end
!\matlabrule!
%% Continue transcritical curve emanating from fold-Hopf point
[~,tcbr]=C1branch_from_C2point(funcs,ht,unfolding_pars,...
    'codim2','zeho','codim1','fold',brpars{1:4},....
    'step',linspace(-8.0e-03,8.0e-03,10),'plot',0);
\end{lstlisting}

\subsection{Detect special points on the Hopf branches}
We continue with detecting bifurcations on the Hopf branches. The Hopf points on the branch \lstinline|hbr_wbifs(i)|$(i=1,2)$ will contain the normal form coefficients \lstinline|L1| and \lstinline|L2|. These will be used to visualize the criticality of the Hopf points (sub or super) in the bifurcation diagram.
\begin{lstlisting}[style=customMatlab]
%% Detect special points on the Hopf branches
for i=2:-1:1
  [hbr_wbifs(i),hopftests(i),hc2_indices,hc2_types]=...
  LocateSpecialPoints(funcs,hbr(i));
  al{i}=arrayfun(@(x)x.parameter(ind.mu1),hbr_wbifs(i).point);
  figure(i); clf;
  plot(al{i},hopftests(i).zeho(1,:),'.-',al{i},zeros(size(al{i})));
  xlabel('$\mu_1$','Interpreter','LaTex');
  ylabel('First Lyapunov coefficient (L1)')
  title('Criticality along Hopf bifurcation curve')
end
\end{lstlisting}

\subsection{Predictors}
For comparison in the bifurcation diagram we obtain predictors for the various branches by setting the argument \lstinline|predictor| to 1 and \lstinline|step| to a range of $\varepsilon$-values when calling the function \lstinline|C1branch_from_C2point|.
\begin{lstlisting}[style=customMatlab]
%% Predictors for Neimark-Sacker, Hopf and transcritical curves
[trfuncs,nsbr_pred]=C1branch_from_C2point(funcs,ht,...
    unfolding_pars,'codim2','zeho','codim1','TorusBifurcation',...
    'step',linspace(1e-05,4e-02,20),'predictor',true);
[~,tcbr_pred]=C1branch_from_C2point(funcs,ht,unfolding_pars,...
    'codim2','zeho','codim1','fold',...
    'step',linspace(-8.0e-03,8.0e-03,10),'predictor',true);
[~,hbrsub_pred]=C1branch_from_C2point(funcs,ht,unfolding_pars,...
    'codim2','zeho','codim1','hopf',brpars{:},...
    'step',linspace(-1e-01,0,40),'predictor',1);
[~,hbrsup_pred]=C1branch_from_C2point(funcs,ht,unfolding_pars,...
    'codim2','zeho','codim1','hopf',brpars{:},...
    'step',linspace(0,1e-01,40),'predictor',1);
% Correct super- and subcritical hopf branches for the first curve
tempbr=hbrsub_pred(1);
hbrsub_pred(1)=hbrsuper_pred(1);
hbrsuper_pred(1)=tempbr;
\end{lstlisting}

\subsection{Bifurcation diagram}
\begin{figure}[ht]
\centering
\ifcompileimages
  \tikzsetnextfilename{VDPO_bifurcation_diagram}%
  \input{tikz/VDPO_bifurcation_diagram}%

\else
\includegraphics{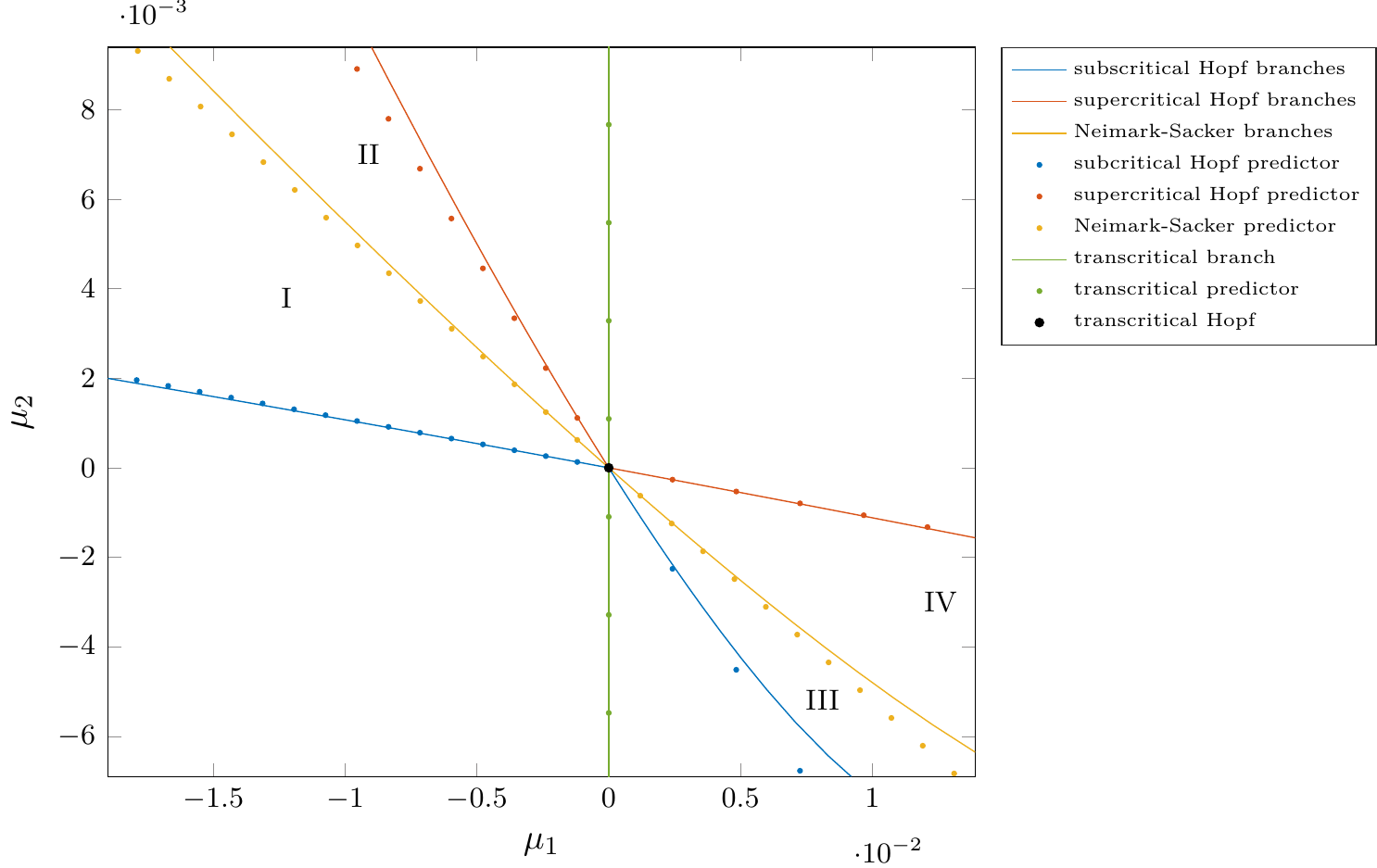}
\fi
\caption{\label{sm:fig:HT_bifurcation_diagram} Bifurcation diagram near the transcritical-Hopf bifurcation in the delayed Van der Pol's oscillator given by \cref{eq:vdp}. There are two supercritical Hopf curves (blue), two subcritical Hopf curves (red), two Neimark-Sacker curves (yellow), and one transcritical curve (green).  We see that the predictors (dotted) give a good approximation for nearby values.}
\end{figure}

We plot the obtained curves and the predictors for the Neimark-Sacker, Hopf, and transcritical curves with the following code.
\begin{lstlisting}[style=customMatlab]
%% Plot comparing computed and predicted curves
figure(4); clf; hold on;
plot(ht.parameter(ind.mu1),ht.parameter(ind.mu2),'k.'...
    ,'MarkerSize',12)
tcbr_pm_pred = [getpars(tcbr_pred,ind.mu1); ...
    getpars(tcbr_pred,ind.mu2)];
plot(tcbr_pm_pred(1,:),tcbr_pm_pred(2,:),'.','Color',cm(5,:));
plot(tcbr_pm(1,:),tcbr_pm(2,:),'Color',cm(5,:));
for i=2:-1:1
  nsbr_pm_pred{i} = [getpars(nsbr_pred(i),ind.mu1); ...
    getpars(nsbr_pred(i),ind.mu2)];
  hbrsub_pm_pred{i}  = [getpars(hbrsub_pred(i),ind.mu1); ...
    getpars(hbrsub_pred(i),ind.mu2)];
  hbrsup_pm_pred{i} = [getpars(hbrsup_pred(i),ind.mu1); ...
   getpars(hbrsup_pred(i),ind.mu2)];
  plot(nsbr_pm_pred{i}(1,:),nsbr_pm_pred{i}(2,:),...
    '.','Color',cm(3,:));
  plot(hbrsub_pm_pred{i}(1,:),hbrsub_pm_pred{i}(2,:),....
    '.','Color',cm(2,:));
  plot(hbrsup_pm_pred{i}(1,:),hbrsup_pm_pred{i}(2,:),...
    '.','Color',cm(1,:));
  plot(nsbr_pm{i}(1,:),nsbr_pm{i}(2,:),'Color',cm(3,:));
  plot(hbrsup_pm{i}(1,:),hbrsup_pm{i}(2,:),'Color',cm(2,:));
  plot(hbrsub_pm{i}(1,:),hbrsub_pm{i}(2,:),'Color',cm(1,:));
end
legend({'transcritical Hopf','transcritical predictor',...
'transcritical curve','Neimark-Sacker predictor',...
'subcritical Hopf predictor','supercritical Hopf predictor',...
'Neimark-Sacker branch','subcritical Hopf branch',...
'superscritical Hopf branch'})
title('Neimark-Sacker curve emanating from the transcrical-Hopf point')
axis([-0.0190    0.0139   -0.0069    0.0094])
xlabel('$\mu_1$','Interpreter','LaTex')
ylabel('$\mu_2$','Interpreter','LaTex')
text(-0.0129,0.0038,'I');  text(-0.0093,0.007,'II');
text(0.008,-0.0052,'III'); text(0.0115,-0.003,'IV');
legend('Location','NorthWest'); box on
% Reverse the stacking order of the graphics
chi=get(gca,'Children'); set(gca,'Children',flipud(chi));
\end{lstlisting}
\cref{sm:fig:HT_bifurcation_diagram} shows the resulting bifurcation diagram.

\subsection{Plot comparing computed and predicted periodic orbits} \label{sm:sec:vdpo:comparing_period_orbits}
\begin{figure}
	\centering
	\includegraphics{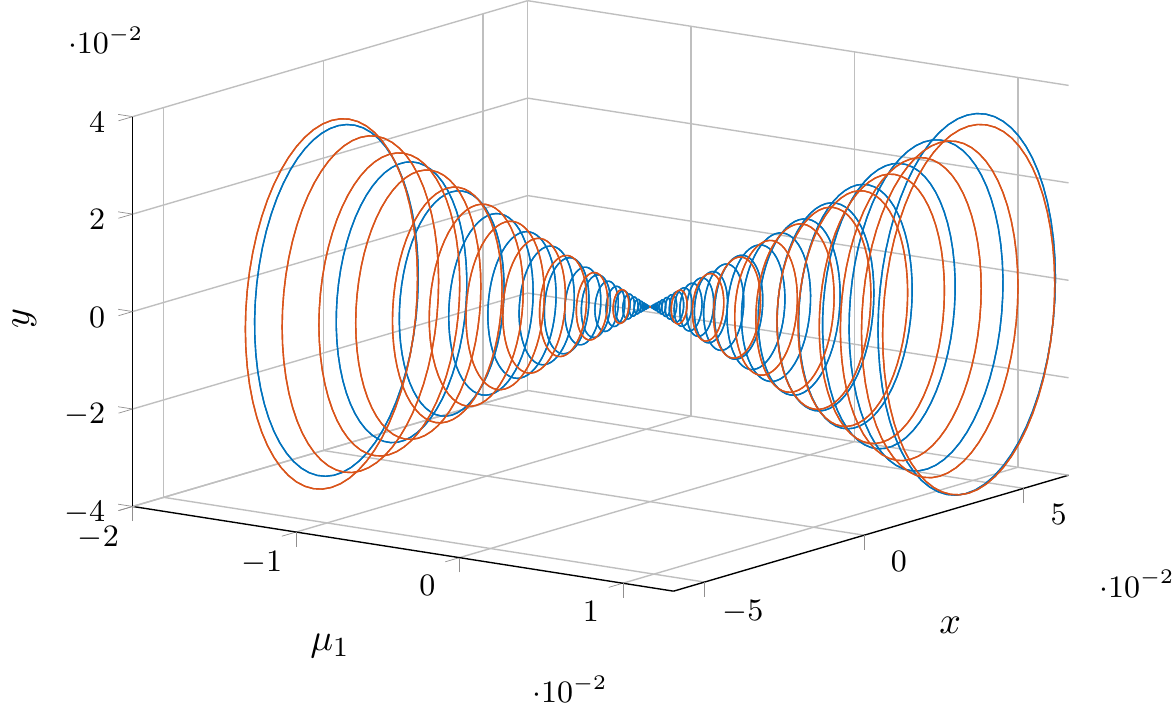}
	\caption{Comparison between predicted periodic orbits (red) and computed periodic orbits (blue) emanating from the transcritical-Hopf bifurcation}
	\label{sm:VDPO_compare_periodic_orbits}
\end{figure}

Lastly, we create a plot to compare the computed and predicted periodic orbits.
\begin{lstlisting}[style=customMatlab]
%% Plot comparing computed and predicted periodic orbits
figure(5); clf; hold on;
genpars=@(br,i)ones(2,length(...
    br.point(1).profile(1,:))).*br.point(i).parameter(ind.mu1);
for i=1:23
  plot3(genpars(nsbr(1),i),nsbr(1).point(i).profile(1,:),...
    nsbr(1).point(i).profile(2,:),'Color',cm(1,:));
end
for i=1:12
  plot3(genpars(nsbr_pred(1),i),...
    nsbr_pred(1).point(i).profile(1,:),...
    nsbr_pred(1).point(i).profile(2,:),'Color',cm(2,:));
end
for i=1:20
  plot3(genpars(nsbr(2),i),nsbr(2).point(i).profile(1,:),...
    nsbr(2).point(i).profile(2,:),'Color',cm(1,:));
end
for i=1:12
  plot3(genpars(nsbr_pred(2),i),...
    nsbr_pred(2).point(i).profile(1,:),...
    nsbr_pred(2).point(i).profile(2,:),'Color',cm(2,:));
end
title('Comparison between computed and predicted periodic orbits')
xlabel('\mu_1'); ylabel('x'); zlabel('y');
view(3)
\end{lstlisting}
The resulting plot is show in \cref{sm:VDPO_compare_periodic_orbits}.

\subsection{Simulation near transcritical-Hopf point with \texttt{pydelay}}
We simulate the dynamics in regions $\text{III}$ and $\text{IV}$ of  \cref{sm:fig:HT_bifurcation_diagram}. Since the critical normal form coefficients are such that
\[
s=1,\qquad \theta<0, \qquad e<0,
\]
a stable cycle and stable torus should be present. 
The simulation in regions $\text{I}$ and $\text{II}$ have also been carried out, but have been omitted here. The following code can be found in the file \lstinline|vdpo_simulation.py|.
\begin{lstlisting}[language=Python,escapechar=!]
import numpy as np
import pylab as pl
from pydelay import dde23
from matplotlib import colors

# Number of time units
tfinal = 30000

# Define DDE
eqns = {
    'x' : '(tau0+mu2)*y',
    'y' : '(tau0+mu2)*((1+mu1)*x(t-tau)-0.2*y(t-tau)'
           '-0.2*pow(x(t-tau),2)'
           '-0.2*x(t-tau)*y(t-tau)-0.2*pow(y(t-tau),2)'
           '-epsilon*(pow(x,2)-1)*y-x+0.5*pow(x,3))'
}

# Set parameters
# Period cycle
tau0=1.757290761249588
params1 = {
    'mu1':0.0049,
    'mu2':-0.0031,
    'tau0':tau0,
    'epsilon':0.3,
    'tau':1
}aH
/tru
# Torus
params2 = {
    'mu1':-0.006871405962603,
    'mu2':0.003871232826592+0.00001,
    'tau0':1.757290761249588,
    'epsilon':0.3,
    'tau':1
}

# Solve DDE
dde = dde23(eqns=eqns, params=params1)
dde.set_sim_params(tfinal=tfinal, dtmax=0.1, AbsTol=1e-8, RelTol=1e-6)
histfunc = {'x': lambda t: 0, 'y': lambda t: -0.2 }
dde.hist_from_funcs(histfunc, 51)
dde.run()

# Subtract solution components
M=600
sol1 = dde.sample(tfinal-M, tfinal, 0.1)
t = sol1['t']
x = sol1['x']
y = sol1['y']

# Plot the solution in phase-space
pl.plot(x, y)
pl.show()

# Plot the solution in phase-space
import matplotlib as mpl
from mpl_toolkits.mplot3d import Axes3D
import matplotlib.pyplot as plt
# Substract delayed solution component
soltau = dde.sample(tfinal-M-tau0,tfinal-tau0, 0.1)
xtau = soltau['y']

# Poincar!{\color{comment}\'e}! section
def poincaresection(x, xtau, y, x_label, y_label, val):
  zero_cross = np.where(np.diff(np.sign(xtau-val)))
  plt.figure(1)
  plt.xlabel(x_label)
  plt.ylabel(y_label)
  plt.plot(x[zero_cross], y[zero_cross],'.', c='royalblue')
  plt.show()
  return

x_label='$x(t)$'
y_label='$x(t-\tau)$'
poincaresection(x,xtau,y,x_label,y_label,-0.009)
\end{lstlisting}

\begin{figure}[ht]
\centering
\ifcompileimages
  \tikzsetnextfilename{vdpo_simulation1}%
  \input{tikz/vdpo_simulation1}%
 \\
  \tikzsetnextfilename{vdpo_simulation2}%
  \input{tikz/vdpo_simulation2}%
 \\[0.6cm]
  \tikzsetnextfilename{vdpo_simulation3}%
  \input{tikz/vdpo_simulation3}%
 \\
  \tikzsetnextfilename{vdpo_simulation4}%
  \input{tikz/vdpo_simulation4}%

\else
\includegraphics{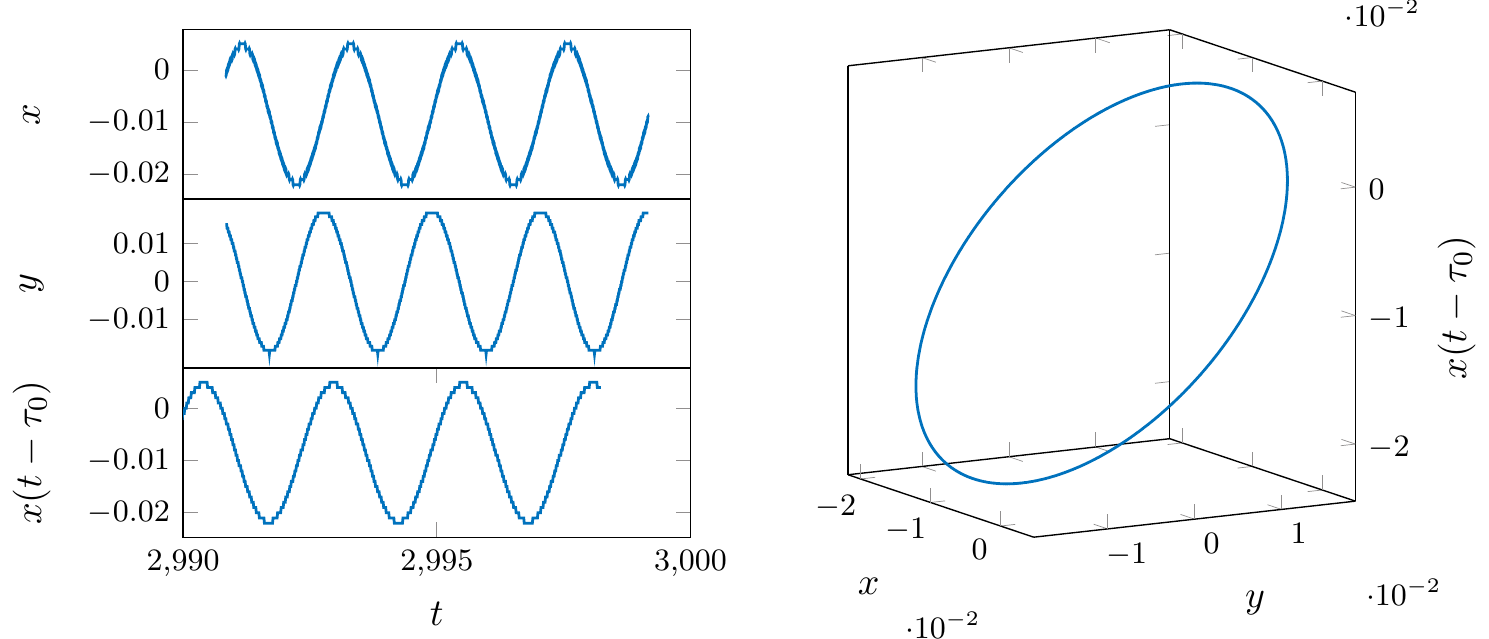}
\includegraphics{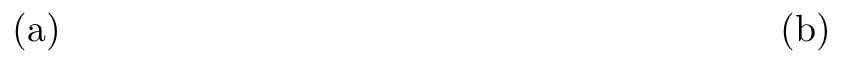}
\includegraphics{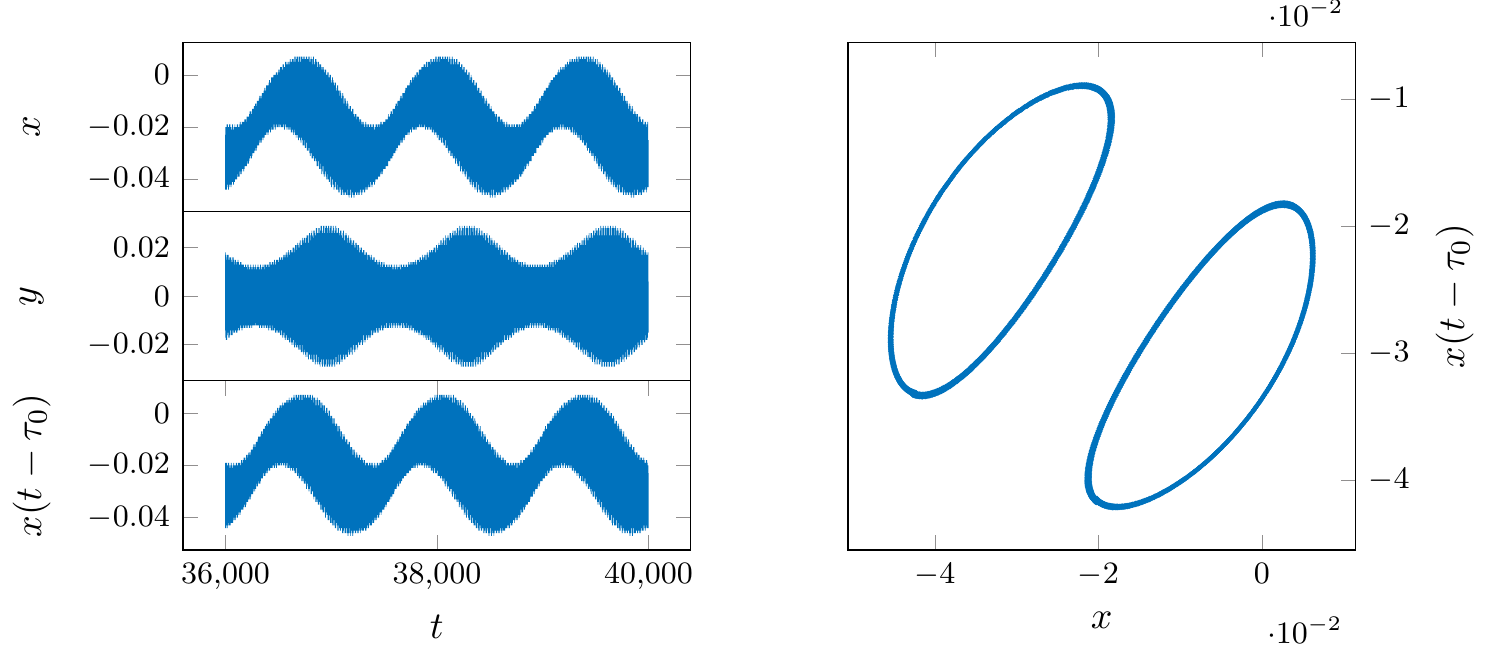}
\includegraphics{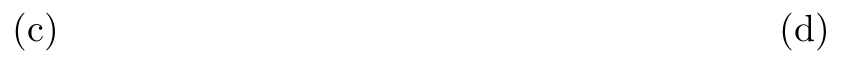}
\fi
\caption{In \textup{(a)} and \textup{(b)}, the periodic solution is shown to be present for parameter values in region $\text{III}$ of  \cref{sm:fig:HT_bifurcation_diagram}. In \textup{(c)}, the torus present in $\text{IV}$ of  \cref{sm:fig:HT_bifurcation_diagram} is shown. In \textup{(d)}, a cross-section of the torus with $y(t)=-0.009$ in the phase-space $(x,x(t-\tau_0),y)$ is taken.}
\end{figure}

\ifsiam
\bibliographystyle{siamplain}
\bibliography{references}

\end{document}
\fi


\fi

\ifarxiv
\clearpage
\setcounter{page}{1}
\renewcommand{\thepage}{\roman{page}}
\pdfbookmark[0]{References}{references}
\renewcommand\refname{References for main text and supplement}
\fi

\bibliographystyle{siamplain}
\bibliography{references}

\end{document}